\documentclass[12pt]{article}
\RequirePackage{amsthm} 
\RequirePackage{amsmath} 
\RequirePackage{amssymb}
\usepackage{amssymb,amsbsy,amsmath,amsfonts,amssymb,amscd}
\usepackage{latexsym}
\usepackage[frenchb,english]{babel}
\theoremstyle{plain}
\newtheorem{thm}{Th\'eor\`eme}[section] 

\newtheorem*{thmppal}{Th\'eor\`eme principal}
\newtheorem{lem}[thm]{Lemme} 
\newtheorem{coro}[thm]{Corollaire}
\newtheorem{prop}[thm]{Proposition}
\newtheorem*{propA}{Proposition \ref{inter}}
\newtheorem*{propB}{Proposition \ref{corr-config}}
\theoremstyle{definition}
\newtheorem{defn}[thm]{D\'efinition} 
\newtheorem{rem}[thm]{Remarque}
\newtheorem{ex}[thm]{Exemple}

\def\deuxind#1#2{{\buildrel{\scriptstyle #1}\over{\scriptstyle #2}}}
\newcommand{\Coess}{C\!oess}
\newcommand{\C}{\mathbb{C}} 
\newcommand{\N}{\mathbb{N}}

\newcommand{\K}{\mathbb{K}}
\renewcommand{\P}{\mathbb{P}}
\newcommand{\xinf}{x_{-\infty}}
\newcommand{\Xinf}{x_{\infty}}
\newcommand{\yinf}{y_{-\infty}}
\newcommand{\Yinf}{y_{\infty}}
\renewcommand{\a}{\alpha}
\renewcommand{\b}{\beta}
\newcommand{\g}{\gamma}
\renewcommand{\d}{\delta}
\renewcommand{\t}{\theta}

\newcommand{\II}{\mathcal I\!\mathcal I}
\newcommand{\winv}{w^{-1}}
\newcommand{\zinv}{z^{-1}}
\newcommand{\tinv}{\tau^{-1}}
\newcommand{\EE}{\mathcal E}
\newcommand{\Reg}{\mbox{Reg}}
\newcommand{\Ex}{\mbox{Ex}}
\begin{document}
\thispagestyle{empty}
\noindent {\it To appear in} Advances in Maths.
\vspace*{0.7cm}
\begin{center}
\LARGE Singularit\'es g\'en\'eriques et quasi-r\'esolutions des vari\'et\'es de
  Schubert pour le groupe lin\'eaire\\

\vspace*{0.8cm}

\large Aur\'elie Cortez\\
\normalsize 28 juin 2002
\end{center}
\selectlanguage{english}
\begin{abstract}
We determine explicitly the irreducible components of the singular
locus of any Schubert variety for $GL_n(\K),\ \K$ being an
algebraically closed field of arbitrary characteristic. We also
describe the generic singularities along each of them.

The case of covexillary Schubert varieties was solved in an earlier
work of the author [Ann. Inst. Fourier {\bf 51} 2 (2001),
375-393]. Here, we first exhibit some irreducible components of the
singular locus of $X_w,$ by describing the generic singularity along 
each of them. Let $\Sigma_w$ be the union of these components. As
mentioned above, the equality $\Sigma_w=\mbox{ Sing }X_w$ is known for
covexillary varieties, and we base our proof of the general case on
this result. More
precisely, we study the exceptional locus of certain 
quasi-resolutions of a non-covexillary
Schubert variety $X_w,$ and we relate the intersection of these loci
to $\Sigma_w.$ Then,
by induction on the dimension, we can establish the equality. 
\end{abstract}

\selectlanguage{frenchb}
\renewcommand\proofname{Preuve}
\renewcommand\refname{Bibliographie}

\vspace{1cm}
\small 
\noindent Aur\'elie Cortez\\
D\'epartement de math\'ematiques\\
Institut Galil\'ee\\
Universit\'e Paris 13\\
99 avenue J.B. Cl\'ement\\
93430 Villetaneuse\\ 
France\\
e-mail : cortez@math.univ-paris13.fr
\normalsize
\renewcommand{\thefootnote}{\fnsymbol{footnote}}
\footnotetext[1]{Mots-Cl\'es : Schubert varieties, singular loci,
  generic singularities.} 
\footnotetext[1]{Classification Math. : 14M15, 20F55.}
\newpage
%
\section*{Introduction}
\setlength{\unitlength}{0.75cm}
L'objet de ce travail est, \`a la suite de \cite{covex}, de d\'ecrire
explicitement les  
composantes irr\'eductibles du lieu singulier d'une vari\'et\'e de Schubert
arbitraire pour $G\!L_n(\K),$ ainsi que la singularit\'e le long de
chacune d'entre elles. Signalons d\`es maintenant que des r\'esultats
analogues ont \'et\'e obtenus de fa\c con concomitante, par des m\'ethodes
compl\`etement diff\'erentes, par L. Manivel d'une part,
S. Billey et G. Warrington d'autre part, et enfin
C. Kassel, A. Lascoux et C. Reutenauer (voir \`a la fin de
cette introduction). 

Notons $G=G\!L_n(\K),$ et $B$ le sous-groupe form\'e des matrices
triangulaires sup\'erieures. Les vari\'et\'es de Schubert dans $G/B$ sont
param\'etr\'ees par le groupe sym\'etrique d'ordre $n,$ not\'e
$\mathfrak S_n.$ Pour $w\in \mathfrak S_n,$ on note $e_w$ le point
$wB$ de $G/B,\ C_w=Be_w$ la cellule de Schubert, et
$X_w=\overline{C_w}$ la vari\'et\'e de Schubert associ\'es. 
Le lieu singulier d'une vari\'et\'e de Schubert est une r\'eunion de vari\'et\'es de
Schubert ; si $X_v\subseteq X_w,$ le point $e_v$ a un voisinage ouvert
dans $X_w$ qui se d\'ecompose comme le produit $C_v\times
\mathcal{N}_{v,w},$ pour une certaine vari\'et\'e  
$\mathcal{N}_{v,w}$ appel\'ee transversale. On a alors le fait
suivant, sur lequel repose notre d\'emarche : $X_v$ est une composante
irr\'eductible du lieu 
singulier Sing $X_w$ si et seulement si la transversale
$\mathcal{N}_{v,w}$ a un unique point singulier. 
L'\'etude des
transversales permet ainsi, du m\^eme coup, d'identifier des composantes du
lieu singulier, et de d\'ecrire la singularit\'e.  

Voici un aper\c cu des r\'esultats obtenus. D'apr\`es \cite{L-Sa}, la
vari\'et\'e de Schubert 
$X_w$ associ\'ee \`a la permutation $w$ de $\mathfrak S_n$ est
singuli\`ere si et seulement s'il existe des entiers $a<b<c<d$ dans
$[1,n]$ v\'erifiant : $w(d)<w(b)<w(c)<w(a)$ -- on dira que ces entiers
forment une configuration (4231) de $w$ -- ou bien
$w(c)<w(d)<w(a)<w(b)$ -- configuration (3412). On d\'emontre que les
composantes irr\'eductibles de Sing $X_w,$ sont param\'etr\'ees par des 
configurations de points du graphe de la permutation $w,$ 
qui sont des raffinements des notions pr\'ec\'edentes. Plus pr\'ecis\'ement, 
on appelle {\it configuration I} de $w$ un ensemble de points du
graphe de $w,$  
$$\mathcal I=\left\{(\Xinf,\Yinf),\ (\xinf,\yinf)\}\cup\{
  (x_i,y_i),\ i\in [-t,-1]\cup [1,s]\right\},$$ 
avec $s,\ t\ge 0,$ comme repr\'esent\'e sur le premier diagramme, tel que
la zone hachur\'ee ne contienne aucun point du graphe. On associe \`a une
telle configuration une permutation $\tau({\scriptstyle\mathcal I})$
en faisant agir un cycle sur ses points, comme repr\'esent\'e sur le
deuxi\`eme diagramme (les points du graphe de
$\tau({\scriptstyle\mathcal I})$ sont repr\'esent\'es par les $\times$).

\vspace*{0.5cm}
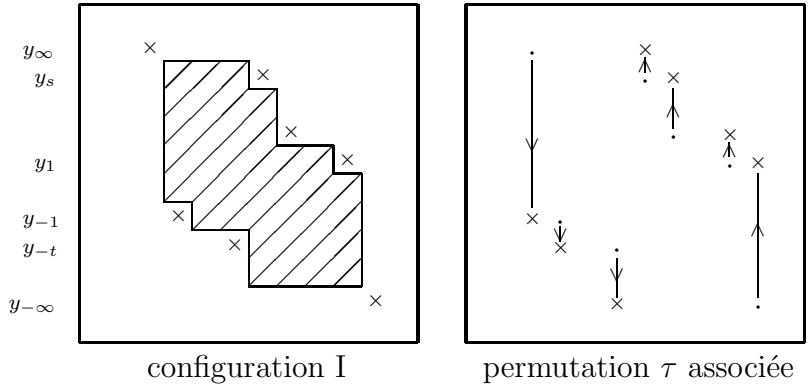
\begin{center}
\begin{minipage}{5cm}
\begin{picture}(6,6)

\thicklines
\put(0,0){\line(1,0){6}}
\put(0,0){\line(0,1){6}}
\put(0,6){\line(1,0){6}}
\put(6,0){\line(0,1){6}}

\thinlines

\put(1.1,5.15){$\scriptstyle\times$}
\put(1.6,2.15){$\scriptstyle\times$}
\put(2.6,1.65){$\scriptstyle\times$}
\put(3.1,4.65){$\scriptstyle\times$}
\put(3.6,3.65){$\scriptstyle\times$}
\put(4.6,3.15){$\scriptstyle\times$}
\put(5.1,0.65){$\scriptstyle\times$}

\put(-1,5.1){$\scriptstyle\Yinf$}
\put(-0.8,4.6){$\scriptstyle y_s $}
\put(-0.8,3.1){$\scriptstyle y_1 $}
\put(-1,2.1){$\scriptstyle y_{-1}$}
\put(-1,1.6){$\scriptstyle y_{-t}$}
\put(-1.25,0.6){$\scriptstyle\yinf$}

\put(1.5,2.5){\line(1,0){0.5}}
\put(2,2){\line(0,1){0.5}}
\put(2,2){\line(1,0){1}}
\put(3,1){\line(0,1){1}}
\put(3,1){\line(1,0){2}}
\put(5,1){\line(0,1){2}}
\put(4.5,3){\line(1,0){0.5}}
\put(4.5,3){\line(0,1){0.5}}
\put(3.5,3.5){\line(1,0){1}}
\put(3.5,3.5){\line(0,1){1}}
\put(3,4.5){\line(1,0){0.5}}
\put(3,4.5){\line(0,1){0.5}}
\put(1.5,5){\line(1,0){1.5}}
\put(1.5,2.5){\line(0,1){2.5}}

\put(1.5,4.5){\line(1,1){0.5}}
\put(1.5,4){\line(1,1){1}}
\put(1.5,3.5){\line(1,1){1.5}}
\put(1.5,3){\line(1,1){1.5}}
\put(1.5,2.5){\line(1,1){2}}
\put(2,2.5){\line(1,1){1.5}}
\multiput(2,2)(0.5,0){3}{\line(1,1){1.5}}
\put(3,1.5){\line(1,1){1.5}}
\put(3,1){\line(1,1){2}}
\put(3.5,1){\line(1,1){1.5}}
\put(4,1){\line(1,1){1}}
\put(4.5,1){\line(1,1){0.5}}

\put(1.2,-0.6){configuration I}
\end{picture}
\end{minipage}
\begin{minipage}{5cm}
\begin{picture}(6,6)

\thicklines
\put(0,0){\line(1,0){6}}
\put(0,0){\line(0,1){6}}
\put(0,6){\line(1,0){6}}
\put(6,0){\line(0,1){6}}

\thinlines

\put(1.03,2.1){$\scriptstyle\times$}
\put(1.53,1.6){$\scriptstyle\times$}
\put(2.53,0.6){$\scriptstyle\times$}
\put(3.03,5.1){$\scriptstyle\times$}
\put(3.53,4.6){$\scriptstyle\times$}
\put(4.53,3.6){$\scriptstyle\times$}
\put(5.03,3.1){$\scriptstyle\times$}

\put(1.1,5){$\cdot$}
\put(1.6,2){$\cdot$}
\put(2.6,1.5){$\cdot$}
\put(3.1,4.5){$\cdot$}
\put(3.6,3.5){$\cdot$}
\put(4.6,3){$\cdot$}
\put(5.1,0.5){$\cdot$}
\linethickness{0.01cm}
\put(1.17,2.4){\line(0,1){2.6}}
\put(1.04,3.4){$\scriptstyle\vee$}
\put(1.67,1.8){\line(0,1){0.25}}
\put(1.54,1.8){$\scriptstyle\vee$}
\put(2.67,0.8){\line(0,1){0.7}}
\put(2.54,1.1){$\scriptstyle\vee$}
\put(5.17,0.8){\line(0,1){2.2}}
\put(5.04,1.9){$\scriptstyle\wedge$}
\put(4.67,3.3){\line(0,1){0.25}}
\put(4.54,3.3){$\scriptstyle\wedge$}
\put(3.67,3.8){\line(0,1){0.7}}
\put(3.54,4){$\scriptstyle\wedge$}
\put(3.17,4.8){\line(0,1){0.25}}
\put(3.04,4.8){$\scriptstyle\wedge$}

\put(0.3,-0.6){permutation $\tau$ associ\'ee}
\end{picture}

\end{minipage}
\end{center}

\vspace*{0.8cm}

On d\'efinit aussi les {\it configurations II} comme des ensembles de
points du graphe
$$\displaylines{\hspace*{1cm}\II=\left\{(a,w(a)),(b,w(b)),(c,w(c)),(d,w(d))\right\}\cup
\hfill\cr\hfill\left\{(c_i,d_i),
i\in[1,r]\right\}\cup \left\{(x_i,y_i),\ i\in [-t,-1]\cup [1,s]\right\},\hspace*{1cm}}$$ 
avec $r,s,t\ge  0,$ comme repr\'esent\'e sur le premier diagramme
ci-dessous, tel que
la zone hachur\'ee ne contienne pas d'autre point du graphe de $w$ que
les $(c_i,d_i),$ et que ceux-ci soient en position relative
Nord-Ouest/Sud-Est. On associe \`a une telle configuration une permutation
$\sigma({\scriptstyle\II})$ en faisant agir un cycle sur ses points,
comme repr\'esent\'e sur le deuxi\`eme diagramme.

\vspace*{0.5cm}
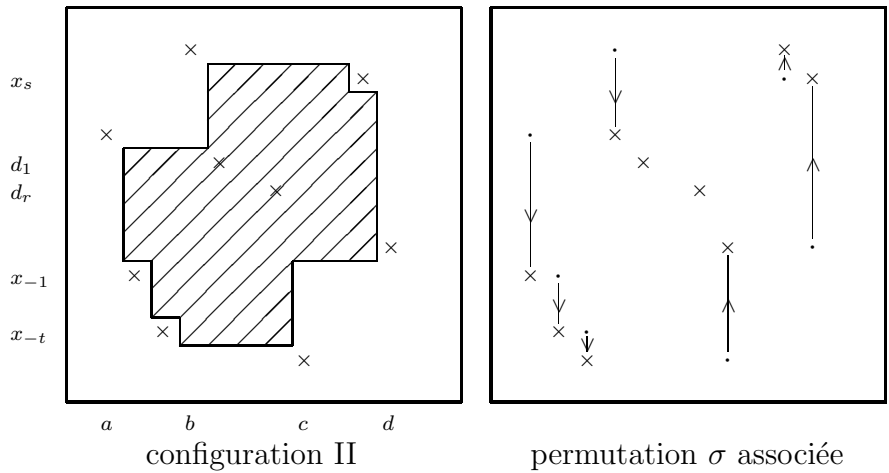
\begin{center}
\begin{minipage}{5.5cm}
\begin{picture}(7,7)

\thicklines
\put(0,0){\line(1,0){7}}
\put(0,0){\line(0,1){7}}
\put(0,7){\line(1,0){7}}
\put(7,0){\line(0,1){7}}

\thinlines
\put(0.55,4.65){$\scriptstyle\times$} 
\put(1.05,2.15){$\scriptstyle\times$} 
\put(1.55,1.15){$\scriptstyle\times$} 
\put(2.05,6.15){$\scriptstyle\times$} 
\put(2.55,4.15){$\scriptstyle\times$} 
\put(3.55,3.65){$\scriptstyle\times$} 
\put(4.05,0.65){$\scriptstyle\times$} 
\put(5.1,5.65){$\scriptstyle\times$} 
\put(5.6,2.65){$\scriptstyle\times$} 

\put(1,2.5){\line(1,0){0.5}}
\put(1.5,1.5){\line(0,1){1}}
\put(1.5,1.5){\line(1,0){0.5}}
\put(2,1){\line(0,1){0.5}}
\put(2,1){\line(1,0){2}}
\put(4,1){\line(0,1){1.5}}
\put(4,2.5){\line(1,0){1.5}}
\put(5.5,2.5){\line(0,1){3}}
\put(5,5.5){\line(1,0){0.5}}
\put(5,5.5){\line(0,1){0.5}}
\put(2.5,6){\line(1,0){2.5}}
\put(2.5,4.5){\line(0,1){1.5}}
\put(1,4.5){\line(1,0){1.5}}
\put(1,2.5){\line(0,1){2}}

\linethickness{0.01mm}
\put(1,4){\line(1,1){0.5}}
\put(1,3.5){\line(1,1){1}}
\put(1,3){\line(1,1){3}}
\put(1,2.5){\line(1,1){3.5}}
\put(1.5,2.5){\line(1,1){3.5}}
\put(1.5,2){\line(1,1){3.5}}
\put(1.5,1.5){\line(1,1){4}}
\put(2,1.5){\line(1,1){3.5}}
\put(2,1){\line(1,1){3.5}}
\put(2.5,1){\line(1,1){3}}
\put(3,1){\line(1,1){1}}
\put(3.5,1){\line(1,1){0.5}}
\put(2.5,5){\line(1,1){1}}
\put(2.5,5.5){\line(1,1){0.5}}
\put(4.5,2.5){\line(1,1){1}}
\put(5,2.5){\line(1,1){0.5}}

\put(0.6,-0.5){$\scriptstyle a$}
\put(2.1,-0.5){$\scriptstyle b$}
\put(4.1,-0.5){$\scriptstyle c$}
\put(5.6,-0.5){$\scriptstyle d$}

\put(-1,5.6){$\scriptstyle x_s$}
\put(-1,4.1){$\scriptstyle d_1$}
\put(-1,3.6){$\scriptstyle d_r$} 
\put(-1,2.1){$\scriptstyle x_{-1}$}
\put(-1,1.1){$\scriptstyle x_{-t}$}

\put(1.4,-1.1){configuration II}
\end{picture}
\end{minipage}
\begin{minipage}{5.5cm}
\begin{picture}(7,7)

\thicklines
\put(0,0){\line(1,0){7}}
\put(0,0){\line(0,1){7}}
\put(0,7){\line(1,0){7}}
\put(7,0){\line(0,1){7}}

\thinlines
\put(0.62,4.6){$\cdot$} 
\put(1.12,2.1){$\cdot$} 
\put(1.62,1.1){$\cdot$} 
\put(2.12,6.1){$\cdot$} 
\put(4.12,0.6){$\cdot$} 
\put(5.12,5.6){$\cdot$} 
\put(5.62,2.6){$\cdot$} 

\put(0.55,2.15){$\scriptstyle\times$} 
\put(1.05,1.15){$\scriptstyle\times$} 
\put(1.55,0.65){$\scriptstyle\times$} 
\put(2.05,4.65){$\scriptstyle\times$} 
\put(2.55,4.15){$\scriptstyle\times$} 
\put(3.55,3.65){$\scriptstyle\times$} 
\put(4.05,2.65){$\scriptstyle\times$} 
\put(5.05,6.15){$\scriptstyle\times$} 
\put(5.55,5.65){$\scriptstyle\times$} 

\linethickness{0.01mm}
\put(0.7,2.4){\line(0,1){2.2}}
\put(0.58,3.2){$\scriptstyle\vee$}
\put(1.2,1.4){\line(0,1){0.7}}
\put(1.08,1.6){$\scriptstyle\vee$}
\put(1.7,0.9){\line(0,1){0.25}}
\put(1.58,0.9){$\scriptstyle\vee$}
\put(4.2,0.9){\line(0,1){1.7}}
\put(4.08,1.6){$\scriptstyle\wedge$}
\put(5.7,2.9){\line(0,1){2.7}}
\put(5.58,4.1){$\scriptstyle\wedge$}
\put(5.2,5.9){\line(0,1){0.25}}
\put(5.08,5.9){$\scriptstyle\wedge$}
\put(2.2,4.9){\line(0,1){1.2}}
\put(2.08,5.3){$\scriptstyle\vee$}

\put(0.7,-1.1){permutation $\sigma$ associ\'ee}
\end{picture}
\end{minipage}
\end{center} 

\vspace*{1cm}
On d\'emontre le 

\noindent {\bf Th\'eor\`eme principal.} {\it Les composantes
irr\'eductibles du lieu singulier de $X_w$ sont les 
$X_{\tau({\scriptscriptstyle\mathcal I})}$ associ\'ees aux configurations I de
param\`etres $s,t\ge 1,$ et les
$X_{\sigma({\scriptscriptstyle\II})}$ associ\'ees aux configurations II de
param\`etres $r,s,t$ tels que $r=0$ ou $s=t=0.$
 
La transversale en une composante
$X_{\tau({\scriptscriptstyle\mathcal I})}$ est isomorphe \`a la vari\'et\'e des
matrices de taille $(s+1,t+1)$ et de rang au plus 1 ; pour une
composante $X_{\sigma({\scriptscriptstyle\II})}$ telle que $r=0,$ la
transversale est isomorphe \`a la vari\'et\'e des
matrices de taille $(s+t+2,2)$ et de rang au plus 1 ; enfin, pour une
composante $X_{\sigma({\scriptscriptstyle\II})}$ telle que $s=t=0,$ la
transversale est un c\^one quadratique non-d\'eg\'en\'er\'e de
dimension $2r+3.$}

On \'etablit ce r\'esultat par une m\'ethode g\'eom\'etrique, poursuivant
le travail accompli dans \cite{covex}. 
On montre dans un premier temps
que les $X_{\tau({\scriptscriptstyle\mathcal I})}$ et les
$X_{\sigma({\scriptscriptstyle\II})}$ de l'\'enonc\'e sont des
composantes irr\'eductibles du lieu singulier, en \'etudiant les
transversales (section 3). Notons $\Sigma_w$ la r\'eunion de ces
composantes. Dans \cite{covex}, on a montr\'e que $\Sigma_w=\mbox{Sing }
X_w$ lorsque $w$ est covexillaire -- c'est-\`a-dire ne contient pas de
configuration (3412) -- en construisant une r\'esolution de $X_w$ qui induit un
isomorphisme au-dessus du compl\'ementaire de $\Sigma_w.$ Dans le cas
g\'en\'eral, on construit cette fois des ``quasi-r\'esolutions'' de $X_w,$
c'est-\`a-dire des vari\'et\'es, \'eventuellement singuli\`eres, qui se
projettent birationnellement sur $X_w.$ On d\'ecrit les lieux
exceptionnels de ces quasi-r\'esolutions, puis on les relie au lieu
singulier de $X_w$ (section 4). Cela permet finalement d'\'etablir
l'\'egalit\'e cherch\'ee par r\'ecurrence sur la dimension de $X_w$ (section
5). A la suite de ce r\'esultat, on obtient, pour chaque composante
irr\'eductible $X_v$ de Sing $X_w,$ le polyn\^ome de Kazhdan-Lusztig
$P_{v,w},$ ainsi que la multiplicit\'e de $X_w$ en $e_v.$

Par ailleurs, on \'etablit au passage un r\'esultat suppl\'ementaire, de
nature plus combinatoire ; rappelons que les vari\'et\'es de Schubert
sont des sous-vari\'et\'es de la vari\'et\'e des drapeaux d\'efinies
par des relations d'incidence. Dans \cite{covex}, on a montr\'e que
lorsque $w$ est covexillaire, les composantes irr\'eductibles de Sing
$X_w$ sont des sous-vari\'et\'es d\'efinies par le renforcement d'une relation
d'incidence. On s'est int\'eress\'e dans le cas g\'en\'eral \`a ce type de
sous-vari\'et\'es ; il s'av\`ere qu'elles ne sont pas n\'ecessairement
irr\'eductibles, et l'on donne la description de leurs composantes
irr\'eductibles. Cela est fait dans la section 2, la
premi\`ere section rassemblant les notations et certains
rappels. Les r\'esultats d\'emontr\'es dans cet article ont
\'et\'e annonc\'es dans \cite{Note}. 

A partir de la description de l'espace tangent obtenue par
V. Lakshmibai et C. S. Seshadri (\cite{L-Se}) et par K. Ryan
(\cite{Ry}), des r\'esultats concernant le lieu singulier 
d'une vari\'et\'e de Schubert ont \'et\'e obtenus, par des m\'ethodes
combinatoires, dans une p\'eriode r\'ecente. D'abord, V. Gasharov a \'etabli dans
\cite{Gash} une direction de la conjecture de Lakshmibai-Sandhya
formul\'ee dans \cite{L-Sa}, en \'etudiant la variation de la dimension
des espaces tangents. Plus r\'ecemment, et de mani\`ere concomitante \`a notre
travail, la description des composantes irr\'eductibles du lieu
singulier a \'et\'e obtenue, presque simultan\'ement, par L. Manivel
(\cite{Man2}) d'une part, S. Billey et G. Warrington (\cite{BW})
d'autre part, et enfin C. Kassel, A. Lascoux et C. Reutenauer
(\cite{KLR}), \'egalement par une \'etude combinatoire de la variation de
la dimension 
des espaces tangents. Suite \`a cela, Manivel a \'egalement donn\'e la 
description des singularit\'es g\'en\'eriques (\cite{Man3}). Soulignons n\'eanmoins que
l'approche g\'eom\'etrique d\'evelopp\'ee ici apporte un \'eclairage
nouveau et permet sans doute une compr\'ehension plus profonde. On peut
d'ailleurs raisonnablement esp\'erer que l'introduction des
quasi-r\'esolutions ait des applications, concernant par exemple le
calcul des polyn\^omes de Kazhdan-Lusztig pour une permutation
arbitraire. Par ailleurs, on peut penser qu'une partie des m\'ethodes introduites
ici puisse se g\'en\'eraliser au cas des autres groupes semi-simples,
alors qu'une difficult\'e se pr\'esente imm\'ediatement 
pour l'approche combinatoire, puisque l'on ne
dispose pas, pour les autres groupes, d'une description de l'espace
tangent aussi maniable que dans le cas du groupe lin\'eaire.  

Je tiens \`a remercier mon directeur de th\`ese, P. Polo, pour le
soutien constant qu'il m'a apport\'e durant l'\'elaboration de ce
travail. Je remercie aussi le rapporteur pour une suggestion qui a
permis d'all\'eger la preuve du th\'eor\`eme \ref{tau=max}.

%
\section{Notations et rappels}
\setlength{\unitlength}{1cm}
$\K$ est un corps alg\'ebriquement clos de caract\'eristique
arbitraire, $G=G\!L_n(\K)$, $B$ est le sous-groupe de Borel des
matrices triangulaires sup\'erieures, $U$ (resp. $U^-$) est le groupe
des matrices triangulaires sup\'erieures (resp. inf\'erieures)
unipotentes, $T$ est le tore maximal des matrices diagonales. 

Pour $w\in \mathfrak S_n$, on note $e_w$ le point $wB$ de $G/B$,
$C_w=B e_w$ la cellule de Schubert, et $X_w=\overline{C_w}$
la vari\'et\'e de Schubert associ\'es. On note $\ell(w)$
le nombre d'inversions de $w$ ; on rappelle que l'on a $\ell(w)=\dim
X_w.$ On d\'esigne par $\Gamma_w$ le graphe de $w.$ 

On note $\K^{\bullet}=\K^1 \subset \cdots \subset \K^n$ le drapeau
standard dans $G/B,$ correspondant au sous-groupe de Borel $B.$ Pour
$g\in G,$ on note $g\K^{\bullet}$ le drapeau associ\'e. 

Tous les intervalles consid\'er\'es ici sont des intervalles de nombres
entiers. Si $i$ et $j$ sont deux entiers distincts de
$[1,n],$ on note $(i,j)$ la transposition de support $\{i,j\},$ et pour $i\le  
n-1,$ on note $s_i$ la transposition simple $(i,i+1).$  
Si $I\subseteq [1,n-1],$ on note $P_I$ le sous-groupe parabolique
contenant $B$ associ\'e, et $\mathfrak S_I$ le sous-groupe parabolique
de $\mathfrak S_n$ correspondant. On note 
$^I\mathfrak S_{min}$ (resp. $^I\mathfrak S_{max}$) l'ensemble des
repr\'esentants minimaux (resp. maximaux) des classes \`a droite de
$\mathfrak S_n$ modulo $\mathfrak S_I.$ 

On rappelle le lemme suivant, qui se v\'erifie ais\'ement en recensant les
inversions de $v$ et $v'.$  

\begin{lem}\label{difflongtransp}
Si $v'=(i,j)\,v$ avec $i<j$ et $v^{-1}(i)>
v^{-1}(j),$ on a $$\ell(v')=\ell(v)-1-2\#\{i<k<j\ |\
v^{-1}(j)<v^{-1}(k)<v^{-1}(i)\}.$$
\end{lem}

L'inclusion des vari\'et\'es de Schubert induit l'ordre de
Bruhat-Chevalley sur le groupe sym\'etrique : pour $v,w \in
\mathfrak S_n,\; v\leq w \iff X_v\subseteq X_w$.
\subsection{Singularit\'es g\'en\'eriques}\label{singgen}

Les composantes irr\'eductibles du lieu singulier de la vari\'et\'e de
Schubert $X_w$ sont donn\'ees par les permutations
maximales $v$ telles que le point $e_v$ soit un point singulier de
$X_w$. 

Etant donn\'e $v\le w$, l'ensemble $v(U^-) e_v\cap X_w$ est le voisinage
standard de $e_v$ dans $X_w.$ D'apr\`es la d\'ecomposition de Bruhat,
il est isomorphe au produit $C_v\times \mathcal{N}_{v , w} $ o\`u 
$\mathcal{N}_{v , w}= [v(U^-)\cap U^-] e_v\cap X_w $ ({\it cf.}
\cite{K-L}, Lemma A4). On appelle $\mathcal{N}_{v , w}$ la
transversale \`a $C_v$ dans $X_w.$ La cellule de Schubert $C_v$ \'etant un espace
affine, on a en fait : $X_v$ est une composante irr\'eductible du
lieu singulier de $X_w$ si et seulement si $\mathcal{N}_{v , w}$ a
$e_v$ pour unique point singulier. Ce sont les singularit\'es
g\'en\'eriques que l'on va d\'ecrire. 

Les situations connues jusqu'ici ({\it cf.} \cite{BP}, 3.3 et 4.6,
voir aussi \cite{BV}) sugg\`erent
de d\'efinir les deux types de singularit\'e suivants. Soit $X_v$ une
composante irr\'eductible du lieu singulier de $X_w$ ; on dira qu'elle
est de type $S_1$ s'il existe des entiers $s$ et $t$ ($s,t\ge 2$) tels que 
$\mathcal{N}_{v,w}$ soit isomorphe \`a la vari\'et\'e $\mathcal
C_{s,t}$ des matrices de taille $(s,t)$ et de rang au plus 1, et l'on
dira que $X_v$ est de type $S_2$ si $\mathcal{N}_{v,w}$ est isomorphe
\`a un c\^one  quadratique non d\'eg\'en\'er\'e de dimension au moins
5. Remarquons que dans le  
second cas, l'anneau local de $\mathcal{N}_{v,w}$ en $e_v$ est
factoriel ({\it cf.} \cite{Mum}, III.7, Example J) alors qu'il ne
l'est pas dans le premier (ceci peut se d\'eduire de \cite{Mum},
III.9, Prop. 1, en consid\'erant la r\'esolution $Z=\{(\mathcal D, u)\ |\ 
\mathcal D\in \P^{s-1}, u : \K^t\longrightarrow \mathcal D\}$ de $\mathcal
C_{s,t}$).  


\subsection{Ordre de Bruhat-Chevalley}
L'ordre de Bruhat-Chevalley sur $\mathfrak S_n$ peut \^etre
d\'ecrit en termes de clef d'une permutation ({\it cf.} \cite{Man},
Prop. 2.1.11).  
Cette description permet de d\'emontrer le lemme suivant. Si $v$ et $w$ sont 
deux permutations de $\{1,\ldots,n\}$ qui co\"\i ncident sur $k$ places, elles 
d\'efinissent naturellement des permutations $\tilde{v}$ et $\tilde{w}$ de 
$\{1,\ldots,n-k\}$ comme suit : soient $ i_1<\cdots <i_k$ dans $[1,n]$ tels 
que $v(i)=w(i)$ pour $i=i_1, \ldots, i_k$. Soit $\varphi$ l'unique bijection 
croissante de $\{1,\ldots ,n\}\setminus \{i_1,\ldots ,i_k\}$ dans $\{1,\ldots
,n-k\}$, et soit $\psi$ l'unique bijection croissante de
$\{1,\ldots ,n\}\setminus \{v(i_1),\ldots ,v(i_k)\}$ dans $\{1,\ldots
,n-k\}$. Soient alors $\tilde{v}=\psi\circ v\circ \varphi^{-1}$ et  
$\tilde{w}=\psi\circ w\circ \varphi^{-1}$ ; ce sont des \'el\'ements de
$\mathfrak S_{n-k}$. On a alors le 
\begin{lem}\label{focalisation} $v\le w \iff \tilde{v}\le\tilde{w}$. 
\end{lem}

On peut aussi caract\'eriser l'ordre de Bruhat-Chevalley \`a l'aide de la fonction 
rang d'une permutation : pour $w\in \mathfrak S_n$ et $p,q\in [1,n],$ on
d\'efinit $$r_w(p,q)=\#\{i\ |\ i\le p \quad\mbox{et}\quad w(i)\le
q\}.$$ On a alors le lemme suivant ({\it cf.} \cite{Man}, Prop. 2.1.12),

\begin{lem}\label{rang} Pour $v,w \in \mathfrak S_n$
$$ v\leq w \iff
r_v(p,q)\geq r_w(p,q), \ \ \forall\, p,q\in [1,n].$$
\end{lem} 

D'autre part, A. Lascoux et M.-P. Sch\" utzenberger ont introduit une nouvelle
approche, en d\'efinissant les rectrices d'une permutation ({\it cf.} \cite
{LS-2}). Il est plus commode pour nos besoins de consid\'erer la notion duale de 
corectrice ({\it cf.} \cite{covex}).

On appelle cograssmanniennes les permutations n'ayant qu'une mont\'ee,
et cobigrassmanniennes les permutations cograssmanniennes dont l'inverse
est aussi cograssmannienne. Une cobigrassmannienne est d\'etermin\'ee par un 
quadruplet d'entiers $(n_0,n_1,n_2,n_3),$ avec $n_0,n_3\in \N,\ n_1,n_2\in 
\N^*$ et $\sum n_i=n$ : on coupe $(n,n-1,\ldots,1)$ en quatre blocs dont les 
cardinaux sont les $n_i,$ et on permute les deux blocs m\'edians. En
termes de graphe, les cobigrassmanniennes sont de la forme suivante :

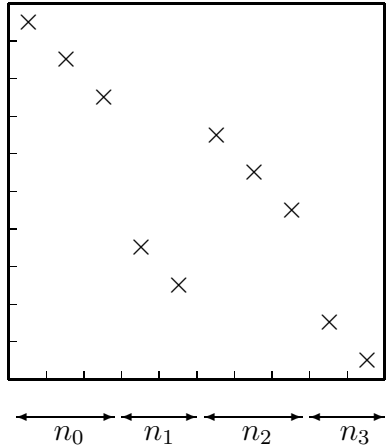
\begin{center}
\begin{picture}(5,6)
\thicklines
\put(0,0.8){\line(1,0){5}}
\put(0,0.8){\line(0,1){5}}
\put(0,5.8){\line(1,0){5}}
\put(5,0.8){\line(0,1){5}}
%
\linethickness{0.01mm}
\multiput(0.5,0.8)(0.5,0){9}{\line(0,1){0.1}}
\multiput(0,1.3)(0,0.5){9}{\line(1,0){0.1}}
\multiput(0.1,5.45)(0.5,-0.5){3}{$\times$}
\multiput(2.6,3.95)(0.5,-0.5){3}{$\times$}
\multiput(1.6,2.45)(0.5,-0.5){2}{$\times$}
\multiput(4.1,1.45)(0.5,-0.5){2}{$\times$}
\thinlines
\put(0.1,0.3){\vector(1,0){1.3}}
\put(1.4,0.3){\vector(-1,0){1.3}}
\put(0.6,0){$n_0$}
\put(1.5,0.3){\vector(1,0){1}}
\put(2.5,0.3){\vector(-1,0){1}}
\put(1.8,0){$n_1$}
\put(2.6,0.3){\vector(1,0){1.3}}
\put(3.9,0.3){\vector(-1,0){1.3}}
\put(3.1,0){$n_2$}
\put(4,0.3){\vector(1,0){1}}
\put(5,0.3){\vector(-1,0){1}}
\put(4.4,0){$n_3$}
\end{picture}
\end{center}

On dispose alors d'un crit\`ere simple pour comparer une permutation
arbitraire et une cobigrassmannienne, analogue \`a \cite{LS-2}, Lemme 4.3 : 
pour $w \in \mathfrak S_n$ et $c$ la cobigrassmannienne d\'efinie par
le quadruplet $(n_0,n_1,n_2,n_3),$ 
on a le 
\begin{lem}\label{w-inf-c}
$w \leq c$ si et seulement si l'ensemble $w([1,n_0+n_1])\cap [1,n_1+n_3]$
contient au moins $n_1$ \'el\'ements.
\end{lem}

Les cobigrassmanniennes permettent de d\'ecrire l'ordre de
Bruhat-\linebreak Chevalley de la
fa{\c c}on suivante : notant $\mathcal C$ l'ensemble des
cobigrassmanniennes, on munit l'ensemble des parties de $\mathcal
C$ de l'ordre inverse de l'inclusion. Alors, d'apr\`es \cite{LS-2}, l'application
qui \`a $w\in \mathfrak S_n$ associe l'ensemble $\{c\in\mathcal C \ |\ w\leq c\}$
induit un isomorphisme d'ensembles ordonn\'es de $\mathfrak S_n$ sur
son image dans $2^{\mathcal C}.$ Etant donn\'e une permutation $w$
de $\mathfrak S_n$, les \'el\'ements minimaux de $\{c\in
\mathcal C \ |\ w\leq c\}$ sont appel\'es corectrices de $w$ ; d'apr\`es
ce qui pr\'ec\`ede, leur donn\'ee d\'etermine enti\`erement $w$. Elles sont 
param\'etr\'ees par l'ensemble coessentiel de $w$, dual de l'ensemble essentiel 
d\'efini par Fulton ({\it cf.} \cite{F2}) : 
$$\Coess (w)=\left\{(p,q)\in [1,n]^2\left\vert 
\begin{array}{l}w(p-1)\le q < w(p)\\w^{-1}(q)\le
p-1<w^{-1}(q+1)\end{array}\right. \right\}.$$  

On notera $c_{p,q}$ la corectrice associ\'ee
au point coessentiel $(p,q)$ de $w$ (voir \cite{covex}, §2.3.2).

\begin{defn}Soit $c$ une cobigrassmannienne de quadruplet
$(n_0,n_1,n_2,n_3).$ On dit que $c$ est {\it it\'erable} si on a
$n_0,n_3 \ge 1,$ et on d\'efinit alors l'it\'er\'ee
$c ^1$ de $c$ par son quadruplet $(n_0-1,n_1+1,n_2+1,n_3-1)$. 
Elle v\'erifie l'in\'egalit\'e $c^1 \le c$.
\end{defn}
\subsection{Quadrants et rectangles}\label{quad}

Il est utile de revenir sur la d\'efinition des quadrants associ\'es \`a un
point $(p,q)$ du carr\'e $[1,n]^2$ donn\'ee dans \cite{covex}, pour
obtenir une notion plus satisfaisante dans le cas g\'en\'eral. On d\'efinit :
\begin{center}
\begin{picture}(13,6)
\put(0,4.4){$NO(p,q)=\{(i,j)\ |\ i\le p, j>q\}$}
\put(0,3.4){$SO(p,q)=\{(i,j)\ |\  i\le p, j\leq q\}$}        
\put(0,2.4){$N\!E(p,q)=\{(i,j)\ |\  i> p, j> q\}$}
\put(0,1.4){$S\hspace{-0.5mm}E(p,q)=\{(i,j)\ |\  i>p, j\le q\}$}
\thicklines
\put(7.5,0.4){\line(1,0){5}}
\put(7.5,0.4){\line(0,1){5}}
\put(7.5,5.4){\line(1,0){5}}
\put(12.5,0.4){\line(0,1){5}}
%
\linethickness{0.01mm}
\multiput(8,0.4)(0.5,0){9}{\line(0,1){0.1}}
\multiput(7.5,0.9)(0,0.5){9}{\line(1,0){0.1}}
\thinlines
\put(9.5,0.4){\line(0,1){5}}
\put(7.5,3.9){\line(1,0){5}}
\put(9.15,0){$p$}
\put(7.1,3.55){$q$}
\put(7.7,4.5){$NO(p,q)$}
\put(7.7,2.1){$SO(p,q)$}
\put(10.3,4.5){$N\!E(p,q)$}
\put(10.3,2.1){$S\hspace{-0.5mm}E(p,q)$}
\end{picture}
\end{center}

On consid\`ere $w\in \mathfrak S_n.$ Les quadrants associ\'es \`a 
$(p,q)$ d\'eterminent naturellement la partition suivante du graphe de $w$ : 
$$\begin{array}{l}NO_w(p,q)=\Gamma_w\cap NO(p,q),\\
SO_w(p,q)=\Gamma_w\cap SO(p,q),\\
NE_w(p,q)=\Gamma_w\cap NE(p,q),\\
SE_w(p,q)=\Gamma_w\cap SE(p,q).\end{array}$$

Il est \`a noter que si $(p+1,q)$ est un point coessentiel de $w,$ on
obtient la m\^eme partition du graphe de $w$ que celle d\'ecrite dans
\cite{covex}, §2.5. 

On note $\partial(NO_w(p,q))$ la fronti\`ere Sud-Est de $NO_w(p,q),$
c.-\`a-d. l'ensemble des points $(x,w(x))\in NO_w(p,q)$ tels que
$NO_w(p,q)$ ne rencontre pas $SE(x,w(x)).$ De m\^eme, on note
$\partial( SE_w(p,q))$) la fronti\`ere Nord-Ouest de $ SE_w(p,q).$
On les notera simplement $NO,$ $NO_w,\ \partial(NO_w),$ {\it etc.}
lorsqu'il n'y aura pas d'ambigu\"\i t\'e. 

\vspace*{1em}

On introduit aussi les notations suivantes ; soient $a$ et $b$ deux
entiers distincts de $[1,n].$ Ecrivant indiff\'eremment $[a,b]$ et
$[b,a]$ pour l'ensemble 
des entiers compris entre $a$ et $b,$ on note : 
$$\begin{array}{lcl}{\mathcal R}_{(a,b)}(w)
  &=&]\,a,b\,[\,\times\,]\,w(a),w(b)\,[\\ 
\overline{\mathcal R}_{(a,b)}(w)&=&[\,a,b\,]\times[\,w(a),w(b)\,]\\
{\mathcal R}^{(a,b)}(w)&=&]\,\winv(a),\winv(b)\,[\,\times\,]\,a,b\,[\\
\overline{\mathcal
  R}^{(a,b)}(w)&=&[\,\winv(a),\winv(b)\,]\times[\,a,b\,].\end{array}$$  
Soient maintenant $A$ et $B$ deux points distincts du graphe de $w,$
de coordonn\'ees respectives $(a,w(a))$ et $(b,w(b)),$ on d\'esigne aussi
par $\mathcal R_{(A,B)}(w)$ (resp. $\overline{\mathcal R}_{(A,B)}(w)$) 
le rectangle ${\mathcal R}_{(a,b)}(w)$ (resp. $\overline{\mathcal
  R}_{(a,b)}(w)$).

%
\section{Renforcement d'une condition d'incidence}\label{renfor}

Dans cette section, on \'etablit un r\'esultat de nature combinatoire,
ind\'ependant de l'\'etude des singularit\'es qui fait l'objet des
sections suivantes ; cette section n'est cependant pas \'etrang\`ere au
reste de ce travail, puisque la question trait\'ee ici a \'et\'e soulev\'ee par
la description donn\'ee dans \cite{covex} des composantes irr\'eductibles
du lieu singulier d'une vari\'et\'e de Schubert covexillaire, et que  
la description de certaines composantes
dans le cas g\'en\'eral fera appel \`a une construction d\'ecrite ici. 

La fonction rang d\'efinie en 1.2 permet de d\'ecrire 
les vari\'et\'es de Schubert en termes de relations d'incidence ({\it cf.} 
\cite{Man}, Prop. 3.6.4) : 
la vari\'et\'e de Schubert $X_w$ est l'ensemble des
drapeaux complets $V^{\bullet}$ de $\K^n$ qui v\'erifient pour tous 
$p,q \in [1,n],$ $$\dim (V^p \cap \K^q) \ge r_w(p,q).$$ 

Dans \cite{covex}, on a montr\'e que si $w$ est covexillaire, les
composantes irr\'eductibles de Sing $X_w$ sont des sous-vari\'et\'es de
$X_w$ d\'efinies par le renforcement d'une condition d'incidence,
c.-\`a-d. de la forme $$\{V^{\bullet}\in X_w\ |\ \dim (V^p \cap \K^q) \ge 
r_w(p,q)+1\}$$
pour certains entiers $p$ et $q.$ Nous nous int\'eressons ici \`a ce
type de sous-vari\'et\'es dans le cas g\'en\'eral : \'etant donn\'e
$(p,q)$ tel que $r_w(p,q)<\mbox{Min}(p,q),$ on consid\`ere
$$X_w^{(p,q)+}=\{V^{\bullet}\in X_w\ |\ \dim (V^p \cap \K^q) \ge 
r_w(p,q)+1\}.$$ 
La cobigrassmannienne $c$ de quadruplet d\'efini par 
$$\begin{array}{l}n_0=p-r_w(p,q)\\
n_1=r_w(p,q)\\
n_2=n-(p+q)+r_w(p,q)\\
n_3=q-r_w(p,q),\end{array}$$ est it\'erable, et d'apr\`es le lemme
\ref{w-inf-c}, elle majore $w,$ et l'on a 
$w\not\le c^1.$ On a $$X_w^{(p,q)+}=X_w\cap X_{c^1}.$$
R\'eciproquement, si $c$ est une cobigrassmannienne it\'erable, de quadruplet
$(n_0,n_1,n_2,n_3),$ telle que $w\le c$ et
$w\not\le c^1,$ posant $p=n_0+n_1$ et $q=n_1+n_3,$ on a $X_w\cap
X_{c^1}=X_w^{(p,q)+}.$ 

Nous allons d\'ecrire les composantes irr\'eductibles de ces
sous-vari\'et\'es, en g\'en\'eralisant la construction de
\cite{covex}. Cela revient \`a d\'ecrire les \'el\'ements maximaux de
$\Lambda(w,c^1),$ o\`u pour des permutations $w_1,\ldots, 
w_k,$ on note $\Lambda(w_1,\ldots, w_k)$
l'ensemble des permutations $z$ telles que $z\le w_i$ pour
$i=1,\ldots, k.$

\vspace*{1em}

Soit $c$ une cobigrassmannienne it\'erable de quadruplet
$(n_0,n_1,n_2,n_3)$ ; on note $$\Omega(c)=\{w\in \mathfrak S_n\ |\
w\le c\mbox{ et } w\not\le c^1\}.$$ On pose $p=n_0+n_1$ et
$q=n_1+n_3.$ On consid\`ere $w\in \Omega(c)$ ; d'apr\`es le lemme
\ref{w-inf-c}, on a 
$$\# \,SO_w=n_1,$$ et il vient alors aussi 
$$\begin{array}{l}\#\, NO_w=n_0,\\
\#\, NE_w=n_2,\\
\#\, SE_w=n_3.\end{array}$$

On voit facilement que la bigrassmannienne 
$$\displaylines{b=(1,\ldots, n_1,\
n_1+n_3+1,\ldots, n_1+n_3+n_0,\hfill\cr\hfill n_1+1, \ldots , n_1+n_3, \
n_1+n_3+n_0+1, \ldots ,n),}$$ qui, avec les notations de \cite{LS-2},
est la bigrassmannienne associ\'ee au quadruplet $(n_1,n_0,n_3,n_2),$ 
est le plus petit \'el\'ement de $\Omega(c).$ 

\begin{prop}\label{cas-b}
Les \'el\'ements maximaux de $\Lambda(b,c^1)$ sont exactement les $(i,j)\,b$
pour $n_1< i\le q,$ et $q<j\le q+n_0.$ 
\end{prop}

\begin{proof}[Preuve]
Comme on a $b(i)=i$ pour $i\in [1,n_1]\cup [n_1+n_3+n_0+1, n],$ toute
permutation $v\le b$ co\"\i ncide avec $b$ sur ces intervalles. On peut
donc supposer $n_1=n_2=0,$ c'est-\`a-dire que $b=(q+1, \ldots, n,\ 1,
\ldots, q),$ avec $n=p+q.$ La condition $v\le c^1$ est alors
$$\#SO_v\ge 1\ \ \ (\star).$$ On 
observe que l'on a $v\le b$ si et seulement si $$v(i)\le q+i \mbox{
  pour } i\le p,\mbox{ et } v(i)\ge i-p \mbox{ pour } i\ge p+1\ \ \ 
(\star\star).$$ Cela s'exprime aussi en termes de graphe, comme
repr\'esent\'e sur la figure suivante : le graphe de $v$ doit \^etre contenu dans la
zone non hachur\'ee 

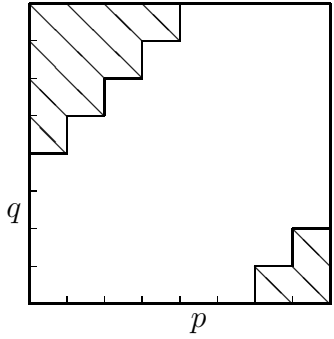
\begin{center}
\begin{picture}(4,4.5)
\thicklines
\put(0,0){\line(1,0){4}}
\put(0,0){\line(0,1){4}}
\put(0,4){\line(1,0){4}} 
\put(4,0){\line(0,1){4}}

\put(2.15,-0.3){$p$}
\put(-0.3,1.15){$q$}

\linethickness{0.01mm}
\multiput(0.5,0)(0.5,0){7}{\line(0,1){0.1}}
\multiput(0,0.5)(0,0.5){7}{\line(1,0){0.1}}
\thinlines
\put(3,0){\line(0,1){0.5}}
\put(3,0.5){\line(1,0){0.5}}
\put(3.5,0.5){\line(0,1){0.5}}
\put(3.5,1){\line(1,0){0.5}}

\multiput(3,0.5)(0.5,0){2}{\line(1,-1){0.5}}
\put(3.5,1){\line(1,-1){0.5}}

\multiput(0,2)(0.5,0.5){4}{\line(1,0){0.5}}
\multiput(0.5,2)(0.5,0.5){4}{\line(0,1){0.5}}

\multiput(0,2.5)(0,0.5){2}{\line(1,-1){0.5}}
\multiput(0,3.5)(0,0.5){2}{\line(1,-1){1}}
\put(0.5,4){\line(1,-1){1}}
\multiput(1,4)(0.5,0){2}{\line(1,-1){0.5}}
\end{picture} 
\end{center}
\vspace*{0.5cm}

On consid\`ere maintenant $v\in \Lambda(b,c^1)$ maximale. Soit
$r=\#\, SO_v=\#\,\{i\le p\ |\ v(i)\le q\}$ ; on a $r\ge 1$ puisque $v\le
c^1.$ On a aussi $r=\#\, NE_v=\#\,\{i>p\ |\ v(i)>q\}.$ Soient
$i_1<\cdots<i_r$ les $i\le p$ tels que $v(i)\le q,$ et
$i'_r<\cdots<i'_1$ les $i>p$ tels que $v(i)>q.$ La maximalit\'e de $v$
entra\^\i ne $v(i_1)>\cdots>v(i_r)$ et $v(i'_r)>\cdots>v(i'_1).$ On veut
montrer que $v(i)=b(i)$ pour $i\not\in \mathcal
A:=\{i_1,\ldots,i_r,i'_r,\ldots,i'_1\},$ et que $r=1.$ Or on a pour tout $i\in
[1,p],$
$$v(i)=b(i) \mbox{ ou } i\in\mathcal A\ \ (\ast).$$
En effet, c'est clair pour
$i=1$ car $b(1)=q+1$ ; soit $i\in [2,p],$ supposons $v(i)\not=b(i)$ et
$i\not\in \mathcal A.$ On a donc $q<v(i)<q+i.$ Or, comme $v\le b,$ on a
pour tout $i'\in [1,i-1], \ v(i')\le q+i'<q+i,$ il vient donc
$v^{-1}(q+i)>i.$ On a ainsi $(v(i),q+i)\,v>v,$ et on v\'erifie sans peine
que $(v(i),q+i)\,v\le b,c^1$ \`a l'aide de $(\star)$ et $(\star\star).$
Cela contredit la maximalit\'e de $v,$ et $(\ast)$ est d\'emontr\'ee. 

On d\'emontre de m\^eme que pour tout $i\in [p+1,n],$    
\begin{center}{$v(i)=b(i)$ ou $i\in\mathcal A.$}\end{center}

Cela entra\^\i ne, comme $q\ge v(i_1)>\cdots>v(i_r)$ et
$v(i'_r)>\cdots>v(i'_1)>q,$ que $v(i_1)=b(i'_1)=i'_1-p$ et
$v(i'_1)=b(i_1)=q+i_1.$ On obtient alors $r=1$ car si $r\ge 2,$ alors
$v<v(i_1,i'_1)\le b,c^1,$ exclu par maximalit\'e de $v.$ Cela montre que
$v=(i'_1-p,i_1+q)\,b,$ c'est-\`a-dire que $v$ est de la forme voulue.  

On remarque de plus que, si $i\le q,$ et $q<j,$ notant $D$ le rectangle
$[\,j-q,i+p-1\,]\times[\,i,j-1\,],$ on a $$r_{(i,j)b}=r_b+\chi_D,$$ o\`u
$\chi_D$ d\'esigne la fonction caract\'eristique de $D.$ Il r\'esulte alors
du lemme \ref{rang} que les permutations $(i,j)\,b$ pour $i\le q,$ et
$q<j,$ sont deux \`a deux incomparables. Cela ach\`eve la preuve de la
proposition \ref{cas-b}.   
\end{proof}

On consid\`ere maintenant un \'el\'ement arbitraire $w$ de $\Omega(c),$ et
l'on se donne $(P_+,P_-)\in \partial(NO_w)\times\partial(SE_w)$
(voir \ref{quad} pour ces notations). 
Notons $(\Xinf,y_{\infty})$ les coordonn\'ees de $P_+,$ et
$(\xinf,y_{-\infty})$ celles de $P_-.$ Si l'ensemble
$NE_w\,\cap\,{\mathcal R}_{(P_+, P_-)}(w)$ est non vide, sa
fronti\`ere Sud-Ouest constitue 
la suite NE associ\'ee \`a $(P_+,P_-).$ Soient $s$ la longueur de cette
suite ($s=0$ si $NE_w\,\cap\, {\mathcal R}_{(P_+, P_-)}(w)=\emptyset $), et 
$(x_i,y_i)_{1\le i\le s}$ les coordonn\'ees des points, index\'es de
sorte que $x_s<\cdots < x_1.$ On a alors $y_1<\cdots <y_s.$ 

De mani\`ere sym\'etrique, si l'ensemble $SO_w\,\cap\, {\mathcal
  R}_{(P_+, P_-)}(w)$ est non 
vide, sa fronti\`ere Nord-Est constitue la suite SO associ\'ee \`a
$(P_+,P_-).$ Soient $t$ la longueur de cette suite et
$(x_i,y_i)_{-t\le i\le -1}$ les coordonn\'ees des points, index\'es de
sorte que $x_{-1}<\cdots < x_{-t}.$ On a alors $y_{-t} <\cdots <y_{-1}.$

\begin{defn}\label{D1}
On notera $\mathcal {X}=\mathcal {X}_+\cup\mathcal {X}_-$
(resp. $\mathcal {Y}=\mathcal {Y}_+\cup\mathcal {Y}_-$)
l'ensemble des abscisses (resp. ordonn\'ees) ainsi distingu\'ees, o\`u 
$$\begin{array}{c}\mathcal {X}_+=\{x_i\ |\ i\in
  [1,s]\cup\{-\infty\}\},\quad \mathcal {Y}_+=\{y_i\ |\ i\in
  [1,s]\cup\{\infty\}\},\\ 
  \mathcal {X}_-=\{x_i\ |\ i\in
  [-t,-1]\cup\{\infty\}\},\quad
  \mathcal {Y}_-=\{y_i\ |\ i\in
  [-t,-1]\cup\{-\infty\}\}.
\end{array}$$

On d\'efinit alors le cycle $\g_{(P_+, P_-)}$ par 
$$\g_{(P_+, P_-)}=(\Yinf,\ y_{-1},\ \ldots,\ y_{-t},\ \yinf,\ y_1,\
\ldots,\ y_s),$$
et la permutation $\tau_{(P_+, P_-)}$ par 
$$\tau_{(P_+, P_-)} = \g_{(P_+, P_-)}\ w.$$
Il pourra \^etre n\'ecessaire par la suite de 
sp\'ecifier la permutation $w$ \`a laquelle $\tau_{(P_+, P_-)}$ est
associ\'ee, en \'ecrivant $\tau_{(P_+, P_-)}(w).$ En revanche, on la notera 
simplement $\tau$ lorsqu'il n'y aura pas d'ambigu\"it\'e sur $w$ et $(P_+, P_-).$

On note $D'$ la r\'egion de $[1,n]^2$ d\'efinie par :
$$D'=\overline{\mathcal R}_{(P_+,P_-)}(w)\setminus
\Big(\bigcup\limits_{i=1}^t SO(x_{-i}-1,y_{-i}-1)\cup
\bigcup\limits_{i=1}^s NE(x_i-1,y_i-1)\Big),$$
et $D$ la r\'egion obtenue en retirant \`a $D'$ la bande d'ordonn\'ee
$\Yinf,$ et la bande d'abscisse $\xinf.$
\end{defn}

On d\'emontre alors comme dans \cite{covex}, lemme 3.5, le

\begin{lem}\label{long}$($a$)\ r_{\tau}=r_w+\chi_D,\ \mbox{en
    particulier } \tau\le w,$\\
  $($b$)\ \ell(\tau)=\ell(w)-(s+t+1).$
\end{lem}

Nous allons d\'emontrer le
\begin{thm}\label{tau=max}
Les permutations $\tau_{(P_+, P_-)},$ pour $(P_+,
P_-)\in\partial(NO_w)\times\partial( SE_w),$ sont
exactement les \'el\'ements maximaux de $\Lambda(w, c^1).$
\end{thm}

\begin{proof}
Montrons pour commencer que les $\tau_{(P_+, P_-)},$ sont des
\'el\'ements maximaux de $\Lambda(w, c^1).$ On fixe $(P_+,
P_-)\in\partial(NO_w)\times\partial( SE_w),$ et l'on note simplement
$\tau$ la permutation $\tau_{(P_+, P_-)}.$ D'apr\`es l'assertion $(a)$
du lemme ci-dessus, $\tau \in \Lambda(w, c^1).$ Supposons que $\tau$
ne soit pas maximal dans $\Lambda(w, c^1)$ ; alors il existe une
transposition $t=(i,j)$ (avec $i<j$) telle que $t\tau$ soit encore
dans $\Lambda(w, c^1),$ avec $\ell(t\tau)=\ell(\tau)+1.$ On a
$r_{\tau}=r_{t\tau}+\chi_R,$ o\`u $R$ d\'esigne le rectangle
$[\tinv(i),\tinv(j)-1]\times [i,j-1].$ Comme $t\tau\le w,$ ce
rectangle doit \^etre contenu dans $D.$ Il vient alors $i\in \mathcal
{Y}_-$ et $j\in \mathcal {Y}_+,$ et donc le point $(p,q)$ est dans
$R.$ Ainsi l'on a $r_{t\tau}(p,q)=r_w(p,q),$ ce qui contredit
$t\tau\le c^1$ ; $\tau$ est donc maximal dans $\Lambda(w, c^1).$

On termine alors la preuve du th\'eor\`eme par r\'ecurrence sur la
longueur de $w\in \Omega(c)$ ; la proposition \ref{cas-b} donne le
r\'esultat pour l'\'el\'ement minimal de $\Omega(c).$   

On se donne maintenant $w\in \Omega(c)$ tel que $w>b.$ On suppose que
pour tout $z\in [b,c]$ tel que $\ell(z)<\ell(w),$ les \'el\'ements
maximaux de $\Lambda(z,c^1)$ sont exactement les $\tau_{(P_+,P_-)}(z)$
pour $(P_+,P_-)\in \partial(NO_z)\times\partial(SE_z).$ 
Compte-tenu de ce qui pr\'ec\`ede, il suffit de d\'emontrer que
$\Lambda(w,c^1)$ a au plus $\#\,\big(\partial(NO_w)\times\partial(SE_w)\big)$
\'el\'ements maximaux. 

Soit $y$ un \'el\'ement maximal de $\Lambda(w,c^1)$ ; comme $w>b,$
il existe soit un entier $j\not= q$ tel que $s_j w<w,$ soit un
entier $i\not=p$ tel que $w s_i<w.$ Pla\c cons-nous par exemple
dans le premier cas (l'autre situation se traite de mani\`ere
similaire). On remarque alors que $s_j
w\in \Omega(c).$ De plus, il r\'esulte de la ``propri\'et\'e Z'' de V. Deodhar
({\it cf.} \cite{propZ}, Theorem 1.1) que $s_j y\le w.$ D'autre part,
l'entier $j$ est choisi de sorte que $s_j c^1<c^1,$ et il vient donc aussi
$s_j y\le c^1.$ On en d\'eduit, par maximalit\'e de $y,$ que $s_j y<y.$ A
nouveau, par la propri\'et\'e Z, on obtient $s_j y\le s_j w.$ Il existe
donc un \'el\'ement maximal $\tau$ de $\Lambda(s_j w,c^1)$ tel que $s_j
y\le \tau.$ On montre alors, \`a l'aide de la ``propri\'et\'e Z'' que $y=s_j
* \tau,$ o\`u $*$ est l'unique loi associative sur $\mathfrak S_{n}$
telle que pour toute transposition simple $s$ et tout $v\in\mathfrak
S_{n}, \ s* v=\max (v, sv).$ On a donc obtenu 
$$\mbox {Max} \big(\Lambda(w,c^1)\big)\subseteq s_j *\mbox{ Max}\big(
\Lambda(s_j w,c^1)\big).$$   

De plus, si $\tau\in \Lambda(s_j w,c^1),$ alors on a $s_j * \tau\le
w,\ c^1.$ Il vient 
donc $$\mbox {Max}\big(\Lambda(w,c^1) \big)\subseteq\mbox{Max} \Big( s_j
*\mbox{ Max} \big(\Lambda(s_j w,c^1)\big)\Big)\hspace{3em}(\dagger).$$

Il nous suffit maintenant de 
montrer que $\#\,\big(\partial(NO_w)\times\partial( SE_w)\big)$ majore le
cardinal du terme de droite de $(\dagger).$ Par hypoth\`ese de r\'ecurrence, on a 
$$\#\,\mbox{Max}\Big(s_j *\mbox{ Max}\big(\Lambda(s_j w,c^1)
\big)\Big)\le \#\,\big(\partial(NO_{s_j w})\times\partial(SE_{s_j w})\big).$$
Supposons $j>q$ ; alors il est clair que $\partial(SE_{s_j
w})=\partial( SE_w),$ et l'on voit facilement que
$\#\,\partial(NO_{s_j w})=\#\,\partial(NO_w)\mbox{ ou } \#\,\partial(NO_w)+1.$
Plus pr\'ecis\'ement, on a $\#\,\partial(NO_{s_j w})=\#\,\partial(NO_w)+1$
si et seulement si $(w^{-1}(j),j)$ appartient \`a $\partial (NO_w),$ et
est l'unique 
point de $NO_w$ au Sud-Est de $(w^{-1}(j+1),j+1).$ C'est le seul cas \`a
consid\'erer, car si $\#\,\partial(NO_{s_j w})=\#\,\partial(NO_w),$ la majoration
cherch\'ee est \'etablie. On a alors $$\partial(NO_{s_j
w})=\partial(NO_w)\setminus\{(w^{-1}(j),j)\}\cup\{(w^{-1}(j),j+1),
(w^{-1}(j+1),j) \}.$$ Notons $P_+$ (resp. $P'_+$) le point de
$\Gamma_{s_j w}$ de coordonn\'ees 
$(w^{-1}(j+1),j)$ (resp. $(w^{-1}(j),j+1)$). 
Pour tout $P_-\in \partial(SE_w)=\partial( SE_{s_j w}),$ on a 
$$s_j * \tau_{(P_+,P_-)}(s_j w) =\tau_{(P_+,P_-)}(s_j w),$$ et $$s_j *
\tau_{(P'_+,P_-)}(s_j w) =s_j \tau_{(P'_+,P_-)}(s_j w).$$ Montrons
qu'alors $$s_j * \tau_{(P_+,P_-)}(s_j w) \le  s_j *
\tau_{(P'_+,P_-)}(s_j w) \ ;$$ on note $x_{-1}<\cdots
<x_{-t}$ les abscisses des points de la suite SO d\'efinie par
$(P_+,P_-).$ Soit $i$ le plus petit entier tel que $x_{-i}>w^{-1}(j)$
; alors les points de la suite SO d\'efinie par $(P'_+, P_-)$ ont pour
abscisses $x_{-i}<\cdots<x_{-t}.$ Les permutations $\tau_{
(P_+,P_-)}(s_j w)$ et $s_j \tau_{(P'_+,P_-)}(s_j w)$ ne diff\`erent qu'en les
abscisses $w^{-1}(j+1)<x_{-1}<\cdots<x_{-i+1}<w^{-1}(j),$ sur lesquelles $s_j
\tau_{(P'_+,P_-)}(s_j w)$ induit la permutation maximale. On
obtient ainsi $$s_j * \tau_{(P_+,P_-)}(s_j w)\le s_j *
\tau_{(P'_+,P_-)}(s_j w).$$ Il vient donc $$\#\,\mbox{Max}\Big(s_j
*\mbox{ Max} \big(\Lambda(s_j w,c^1)\big)\Big)\le \#\Big[\big(\partial(NO_{s_j
    w})\setminus\{P_+\}\big)\times\partial(SE_{s_j w})\Big],$$
c'est la majoration cherch\'ee. On traite de m\^eme le cas $j<q.$
On obtient ainsi le r\'esultat voulu pour $w,$ ce qui ach\`eve la
d\'emonstration du th\'eor\`eme.  
\end{proof}

\section{Des composantes de type $S_1$ et $S_2$}\label{comp types}
\subsection{Configurations I et II}\label{conf I et II}

Dans cette section, on d\'efinit les configurations I et II
d'une permutation $w,$ et on leur associe des permutations
inf\'erieures ou \'egales \`a 
$w.$ Les configurations I sont une version intrins\`eque de la
construction de la section pr\'ec\'edente (c'est-\`a-dire qui ne n\'ecessite
pas le choix pr\'ealable d'une cobigrassmannienne $c$ telle que $w\in
\Omega(c)$). 

\begin{defn}\label{confI}
$\bullet$ On appelle {\it configuration I} de $w$ un ensemble de points du
graphe de $w,$ 
$$\mathcal I=\left\{(\Xinf,\Yinf),\ (\xinf,\yinf)\}\cup\{
  (x_i,y_i),\ i\in [-t,-1]\cup [1,s]\right\},$$ 
avec $s,\ t\ge 0,$ v\'erifiant les in\'egalit\'es
$$\begin{array}{l}\Xinf<x_{-1}<\cdots<x_{-t}<x_s<\cdots<x_1<\xinf,\\
\yinf<y_{-t}<\cdots<y_{-1}<y_1<\cdots<y_s<\Yinf,\end{array}$$ 
et tels que, notant $\mathcal R=\mathcal R_{(\Xinf,\xinf)}(w),$ on ait : 
$$\Gamma_w\cap \mathcal R \subseteq\bigcup_{i=1}^t SO(x_{-i},y_{-i}) \cup
\bigcup_{i=1}^s NE(x_i-1,y_i-1)\hspace{3em}(\bigtriangleup).$$   

Graphiquement, une configuration I est un ensemble de points du graphe
comme repr\'esent\'e sur le diagramme suivant, tel que la zone
hachur\'ee ne contienne aucun point du graphe. 

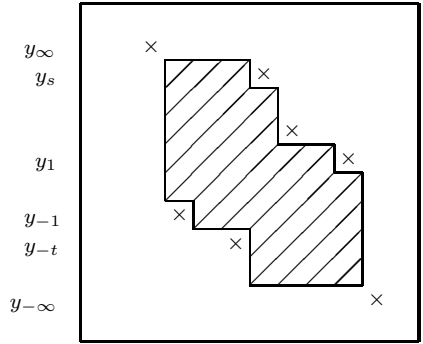
\setlength{\unitlength}{0.75cm}

\begin{center}
\begin{picture}(6,6)

\thicklines
\put(0,0){\line(1,0){6}}
\put(0,0){\line(0,1){6}}
\put(0,6){\line(1,0){6}}
\put(6,0){\line(0,1){6}}

\thinlines

\put(1.1,5.15){$\scriptstyle\times$}
\put(1.6,2.15){$\scriptstyle\times$}
\put(2.6,1.65){$\scriptstyle\times$}
\put(3.1,4.65){$\scriptstyle\times$}
\put(3.6,3.65){$\scriptstyle\times$}
\put(4.6,3.15){$\scriptstyle\times$}
\put(5.1,0.65){$\scriptstyle\times$}

\put(-1,5.1){$\scriptstyle\Yinf$}
\put(-0.8,4.6){$\scriptstyle y_s $}
\put(-0.8,3.1){$\scriptstyle y_1 $}
\put(-1,2.1){$\scriptstyle y_{-1}$}
\put(-1,1.6){$\scriptstyle y_{-t}$}
\put(-1.25,0.6){$\scriptstyle\yinf$}

\put(1.5,2.5){\line(1,0){0.5}}
\put(2,2){\line(0,1){0.5}}
\put(2,2){\line(1,0){1}}
\put(3,1){\line(0,1){1}}
\put(3,1){\line(1,0){2}}
\put(5,1){\line(0,1){2}}
\put(4.5,3){\line(1,0){0.5}}
\put(4.5,3){\line(0,1){0.5}}
\put(3.5,3.5){\line(1,0){1}}
\put(3.5,3.5){\line(0,1){1}}
\put(3,4.5){\line(1,0){0.5}}
\put(3,4.5){\line(0,1){0.5}}
\put(1.5,5){\line(1,0){1.5}}
\put(1.5,2.5){\line(0,1){2.5}}

\put(1.5,4.5){\line(1,1){0.5}}
\put(1.5,4){\line(1,1){1}}
\put(1.5,3.5){\line(1,1){1.5}}
\put(1.5,3){\line(1,1){1.5}}
\put(1.5,2.5){\line(1,1){2}}
\put(2,2.5){\line(1,1){1.5}}
\multiput(2,2)(0.5,0){3}{\line(1,1){1.5}}
\put(3,1.5){\line(1,1){1.5}}
\put(3,1){\line(1,1){2}}
\put(3.5,1){\line(1,1){1.5}}
\put(4,1){\line(1,1){1}}
\put(4.5,1){\line(1,1){0.5}}

\end{picture}
\end{center}

\setlength{\unitlength}{1cm}

\noindent $\bullet$ Si $st=0,$ on dit que la configuration est d\'eg\'en\'er\'ee, au
Nord-Est si $s=0,$ et au Sud-Ouest si $t=0.$ \\
$\bullet$ Si $\mathcal I$ est une configuration I de $w,$ on lui associe le
cycle $$\gamma({\scriptstyle\mathcal I})=(\Yinf,\ y_{-1},\ \ldots,\
y_{-t},\ \yinf,\ y_1,\ \ldots,\ y_s),$$ et l'on d\'efinit la permutation 
$$\tau({\scriptstyle\mathcal I})=\gamma({\scriptstyle\mathcal I})\,w.$$
On d\'efinit alors la r\'egion $D$ comme en \ref{D1}, et l'on
rappelle que l'on a : $r_{\tau({\scriptscriptstyle\mathcal I})}=r_w+\chi_D,$ en
  particulier $ \tau({\scriptstyle\mathcal I})\le w,$ et
  $\ell(\tau({\scriptstyle\mathcal I}))= \ell(w)-(s+t+1).$ 
\end{defn}

On obtient une param\'etrisation des configurations I de $w$ de la
mani\`ere suivante. 

\begin{defn}\label{PtCoBB}
  On dit qu'un point coessentiel $P=(p+1,q)$ de $w,$ est
  {\it bien bord\'e} si le graphe de $w$ rencontre les quadrants
  $NO(p,q)$ et $SE(p,q).$  Dans ce cas, on appelle {\it bordage
  minimal} de $P$ tout couple de points $(P_+, P_-)\in NO_w(p,q)\times
  SE_w(p,q)$
  tels que ${\mathcal R}_{(P_-, P_+)}(w)\cap NO(p,q)$ et ${\mathcal
  R}_{(P_-, P_+)}(w) \cap SE(p,q)$ ne 
  contiennent aucun point du graphe de $w.$
\end{defn}
On fixe un point coessentiel bien bord\'e $P=(p+1,q)$ de $w,$ et un
bordage minimal $(P_+, P_-)$ de $P.$ Alors on a
$(P_+,P_-)\in\partial(NO_w(p,q))\times\partial( SE_w(p,q)),$ et l'on consid\`ere
alors les suites NE et SO associ\'ees, not\'ees $(x_i,y_i)_{1\le i\le s}$
et $(x_i,y_i)_{-t\le i\le 
  -1}.$ Alors l'ensemble $$\mathcal I:=\left\{P_+,P_-\}\cup\{
  (x_i,y_i),\ i\in [-t,-1]\cup [1,s]\right\},$$ est une configuration I
  de $w.$ De plus, toutes les configurations I de $w$ sont obtenues de
  cette mani\`ere. On note $\mathcal {T}_w$ l'ensemble des triplets
  $(P,P_+,P_-),$ o\`u $P$ est un point coessentiel bien bord\'e de $w$ et
  $(P_+,P_-)$ un bordage minimal de $P.$ L'application qui \`a un tel
  triplet associe la configuration I d\'ecrite plus haut n'est en 
g\'en\'eral pas injective ; ces configurations sont param\'etr\'ees par les
classes d'\'equivalence de $\mathcal {T}_w$ pour la relation suivante : 
$(P,P_+,P_-)\sim (Q,Q_+,Q_-)$ si l'on a $P_+=Q_+,\ P_-=Q_-,$ et si 
les quadrants associ\'es aux points coessentiels $P$ et $Q$ d\'efinissent
la m\^eme partition de $\Gamma_w\cap \mathcal R_{(P_+,P_-)}(w).$  

\vspace*{1em}

On passe maintenant aux configurations II.

\begin{defn}$\bullet$ Une configuration (3412) de $w$ est la donn\'ee de quatre
  points du graphe, d'abscisses $a<b<c<d,$ tels que
  $w(c)<w(d)<w(a)<w(b).$ Par abus de langage, on assimilera cette
  donn\'ee \`a celle des abscisses.\\
$\bullet$ Une configuration (3412), d'abscisses $a<b<c<d,$ est dite
  {\it incompressible} s'il n'existe pas d'autre configuration (3412)
  d'abscisses $x<y<c<d$ telle que $(x,w(x))$ (resp. $(y,w(y))$) soit
  au Sud-Est de $(a,w(a))$ (resp. $(b,w(b))$), ni, de mani\`ere
  sym\'etrique, d'autre configuration (3412) d'abscisses $a<b<x<y$ telle
  que $(x,w(x))$ (resp. $(y,w(y))$) soit au Nord-Ouest de $(c,w(c))$
  (resp. $(d,w(d))).$ 
\end{defn} 

La donn\'ee d'une configuration (3412) d\'etermine une partition du
rectangle $\mathcal R_{II}=[\,a,d\,]\times [\,w(c),\,w(b)]$ en neuf
zones. D'une part les zones  
$$\begin{array}{lcl}NO_{II} & = & [\,a,b\,]\times[\,w(a),w(b)\,],\\ 
NE_{II} & = & [\,c,d\,] \times[\,w(a),w(b)\,],\\
SO_{II} & = & [\,a,b\,] \times[\,w(c),w(d)\,],\\
SE_{II} & = & [\,c,d\,] \times[\,w(c),w(d)\,],\end{array}$$ 
la zone centrale $$C=\ ]\,b,c\,[\,\times\, ]\,w(d),w(a)\,[,$$ et enfin les
zones m\'edianes  
$$\begin{array}{lcl}MN&=&]\,b,c\,[\,\times[\,w(a),w(b)\,],\\
MO& =&[\,a,b\,] \times\, ]\,w(d),w(a)\,[,\\
ME& =& [\,c,d\,] \times\, ]\,w(d),w(a)\,[,\\
MS& = &]\,b,c\,[\,\times[\,w(c),\,w(d)].\end{array}$$ 

\begin{ex}
On consid\`ere la permutation 
$$w=(11, 12, 17, 7, 3, 5, 16, 10, 1, 9, 2, 6, 15, 4, 18, 13, 8, 14)\in
\mathcal S_{18},$$  
et la configuration (3412) donn\'ee par les points d'abscisses 2, 7, 11,
17 (cette configuration est incompressible). Sur le diagramme suivant,
on a repr\'esent\'e par des $\otimes$ les points de la configuration, et
par des $\times$ les autres points du graphe de $w.$ On a repr\'esent\'e
le rectangle $\mathcal R_{II}$ associ\'e et les 9 zones d\'ecrites ci-dessus.

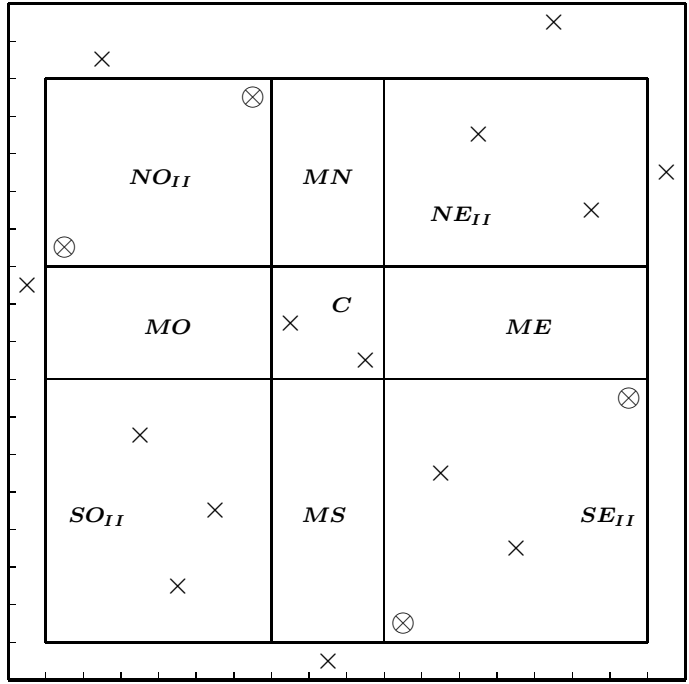
\begin{center}
\begin{picture}(10,10)

\put(0.49,5.55){$\otimes$}
\put(2.99,7.55){$\otimes$}
\put(4.99,0.55){$\otimes$}
\put(7.99,3.55){$\otimes$}

\put(-0.01,5.05){$\times$}
\put(0.99,8.05){$\times$}
\put(1.49,3.05){$\times$}
\put(1.99,1.05){$\times$}
\put(2.49,2.05){$\times$}
\put(3.49,4.55){$\times$}
\put(3.99,0.05){$\times$}
\put(4.49,4.05){$\times$}
\put(5.49,2.55){$\times$}
\put(5.99,7.05){$\times$}
\put(6.49,1.55){$\times$}
\put(6.99,8.55){$\times$}
\put(7.49,6.05){$\times$}
\put(8.49,6.55){$\times$}

\thicklines
\put(-0.1,-0.1){\line(1,0){9}}
\put(-0.1,-0.1){\line(0,1){9}}
\put(-0.1,8.9){\line(1,0){9}}
\put(8.9,-0.1){\line(0,1){9}}

\linethickness{0.01mm}
\multiput(0.4,-0.1)(0.5,0){17}{\line(0,1){0.1}}
\multiput(-0.1,0.4)(0,0.5){17}{\line(1,0){0.1}}

\thinlines
\put(0.4,0.4){\line(1,0){8}}
\put(0.4,3.9){\line(1,0){8}}
\put(0.4,5.4){\line(1,0){8}}
\put(0.4,7.9){\line(1,0){8}}

\put(0.4,0.4){\line(0,1){7.5}}
\put(3.4,0.4){\line(0,1){7.5}}
\put(4.9,0.4){\line(0,1){7.5}}
\put(8.4,0.4){\line(0,1){7.5}}

\boldmath
\put(1.5,6.5){$\scriptstyle{NO_{II}}$}
\put(3.8,6.5){$\scriptstyle{MN}$}
\put(5.5,6){$\scriptstyle{NE_{II}}$}
\put(1.7,4.5){$\scriptstyle{MO}$}
\put(4.2,4.8){$\scriptstyle{C}$}
\put(6.5,4.5){$\scriptstyle{ME}$}
\put(0.7,2){$\scriptstyle{SO_{II}}$}
\put(3.8,2){$\scriptstyle{MS}$}
\put(7.5,2){$\scriptstyle{SE_{II}}$}
\unboldmath

\end{picture}
\end{center}

\end{ex}

Soit maintenant $a<b<c<d$ une configuration (3412) incompressible. Alors,
d'une part, aucune des quatre zones m\'edianes ne contient de point du
graphe de $w,$ et, d'autre part, si 
$C\cap \Gamma_w$ est non vide, notant ses points $(c_i,d_i)_{1\le i\le r},$ 
avec $c_{1}<\cdots <c_{r},$ on a $d_1>\cdots >d_r.$ Cette suite de
points du graphe de $w$ est appel\'ee la suite centrale associ\'ee \`a la
configuration $a<b<c<d.$

On consid\`ere maintenant l'ensemble $NE_{II}\cap\Gamma_w.$ S'il n'est pas vide, 
sa fronti\`ere Sud-Ouest constitue la suite NE associ\'ee \`a la configuration 
$a<b<c<d.$ Soient $s$ la longueur de cette suite, et $(x_i,y_i)_{1\le i\le s}$ 
ses points, index\'es de sorte que $c< x_s <
\cdots < x_1 < d.$ On a alors $w(a) <y_1 < \cdots < y_s <w(b).$  

De mani\`ere sym\'etrique, on consid\`ere $SO_{II}\cap\Gamma_w.$ S'il est non vide, 
sa fronti\`ere Nord-Est constitue la suite SO associ\'ee \`a la configuration 
$a<b<c<d.$ Soient $t$ la longueur de cette suite et $(x_i,y_i)_{-t\le i\le -1}$
ses points, index\'es de sorte que $a< x_{-1}< \cdots <
x_{-t} < b.$ On a alors $w(c)<y_{-t}< \cdots <y_{-1} < w(d).$
 
\begin{defn}$\bullet$ On appelle {\it configuration II} de $w$ la donn\'ee
  d'une configuration (3412) incompressible et des trois suites de
  points associ\'ees comme ci-dessus. Une configuration II n'ayant pas
  de suite centrale sera dite {\it mixte}, et une configuration II
  avec une suite centrale mais n'ayant ni suite NE ni suite SO sera
  dite {\it pure}.\\ 
$\bullet$ On notera 
$$\begin{array}{c}\mathcal X_-=\{x_i\ |\  i\in [-t,-1]\},\quad \mathcal
  Y_-=\{y_i\ |\  i\in [-t,-1]\},\\
  \mathcal C=\{c_i\ |\  i\in [1,r]\},\quad \mathcal
  D=\{d_i\ |\ i\in [1,r]\},\\
  \mathcal X_+=\{x_i\ |\ i\in [1,s]\},\quad  
  \mathcal Y_+=\{y_i\ |\  i\in [1,s]\}.\\
  \end{array}$$ 

On d\'efinit alors le cycle $\g({\scriptstyle\II})$ par
$$\g({\scriptstyle\II})=\big(w(a),\ y_{-1},\ldots,\ y_{-t},\ w(c),\
w(d),\ y_1,\ldots,\ y_s,\ w(b)\big)$$
et la permutation $\sigma({\scriptstyle\II})$ par 
$$\sigma({\scriptstyle\II})=\g({\scriptstyle\II})\ w.$$

Soit $D'_{II}$ la r\'egion du carr\'e $[1,n]^2$ d\'efinie par 
$$\displaylines{D'_{II}=\mathcal R_{II} \setminus\Big(\bigcup\limits_{i=1}^t
SO(x_{-i}-1,y_{-i}-1)\cup \bigcup\limits_{i=1}^s
NE(x_i-1,y_i-1)\hfill\cr
\hfill\cup NO(b-1,w(a)-1)\cup SE(c-1,w(d)-1)
\Big),}$$
et soit $D_{II}$ la r\'egion obtenue en retirant \`a $D'_{II}$ la bande
d'ordonn\'ee $w(b)$ et la bande d'abscisse $d.$
\end{defn}

\begin{ex}
On reprend la permutation consid\'er\'ee dans l'exemple pr\'ec\'edent
$$w=(11, 12, 17, 7, 3, 5, 16, 10, 1, 9, 2, 6, 15, 4, 18, 13, 8, 14)\in
\mathcal S_{18},$$  
et la configuration (3412) incompressible donn\'ee par les points
d'abscisses 2, 7, 11, 17. Sur le diagramme suivant, on a repr\'esent\'e
par des \boldmath $\oplus$ \unboldmath les quatre points de la
configuration, par des $\oplus$ les points des trois suites
d\'ecrites ci-dessus, et par des $+$ les autres points du graphe de $w.$
 
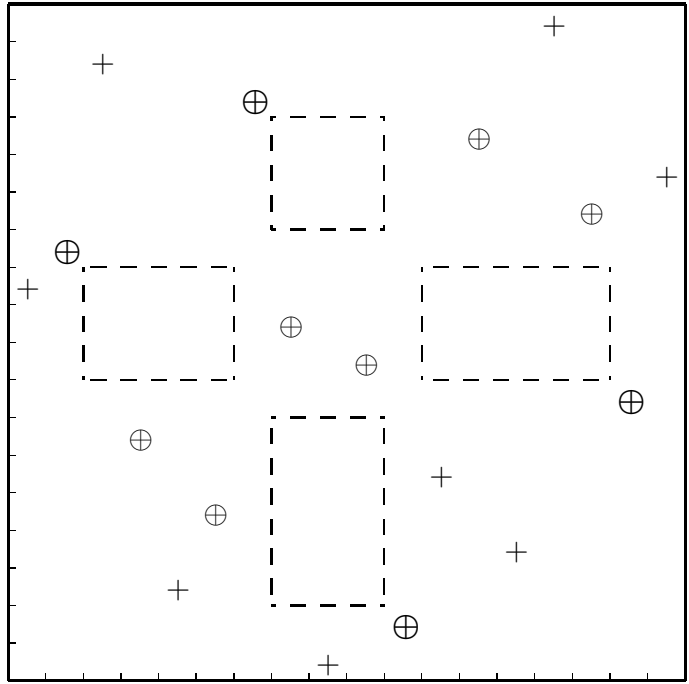
\begin{center}
\begin{picture}(10,10)
\boldmath
\put(0.5,5.5){$\oplus$}
\put(3,7.5){$\oplus$}
\put(5,0.5){$\oplus$}
\put(8,3.5){$\oplus$}
\unboldmath

\put(0,5){+}
\put(1,8){+}
\put(1.5,3){$\oplus$}
\put(2,1){+}
\put(2.5,2){$\oplus$}
\put(3.5,4.5){$\oplus$}
\put(4,0){+}
\put(4.5,4){$\oplus$}
\put(5.5,2.5){+}
\put(6,7){$\oplus$}
\put(6.5,1.5){+}
\put(7,8.5){+}
\put(7.5,6){$\oplus$}

\put(8.5,6.5){+}
\thicklines
\put(-0.1,-0.1){\line(1,0){9}}
\put(-0.1,-0.1){\line(0,1){9}}
\put(-0.1,8.9){\line(1,0){9}}
\put(8.9,-0.1){\line(0,1){9}}

\linethickness{0.01mm}
\multiput(0.4,-0.1)(0.5,0){17}{\line(0,1){0.1}}
\multiput(-0.1,0.4)(0,0.5){17}{\line(1,0){0.1}}

\thinlines
\multiput(1,3.9)(0.4,0){5}{\line(1,0){0.2}}
\put(0.9,3.9){\line(0,1){0.05}}
\multiput(0.9,4.15)(0,0.4){3}{\line(0,1){0.2}}
\put(0.9,5.35){\line(0,1){0.05}}
\multiput(1,5.4)(0.4,0){5}{\line(1,0){0.2}}
\put(2.9,3.9){\line(0,1){0.05}}
\multiput(2.9,4.15)(0,0.4){3}{\line(0,1){0.2}}
\put(2.9,5.35){\line(0,1){0.05}}

\put(3.4,5.9){\line(1,0){0.05}}
\multiput(3.65,5.9)(0.4,0){3}{\line(1,0){0.2}}
\put(4.85,5.9){\line(1,0){0.05}}

\put(3.4,5.9){\line(0,1){0.05}}
\multiput(3.4,6.15)(0,0.4){3}{\line(0,1){0.2}}
\put(3.4,7.35){\line(0,1){0.05}}

\put(4.9,5.9){\line(0,1){0.05}}
\multiput(4.9,6.15)(0,0.4){3}{\line(0,1){0.2}}
\put(4.9,7.35){\line(0,1){0.05}}

\put(3.4,7.4){\line(1,0){0.05}}
\multiput(3.65,7.4)(0.4,0){3}{\line(1,0){0.2}}
\put(4.85,7.4){\line(1,0){0.05}}

\put(3.4,0.9){\line(1,0){0.05}}
\multiput(3.65,0.9)(0.4,0){3}{\line(1,0){0.2}}
\put(4.85,0.9){\line(1,0){0.05}}

\put(3.4,3.4){\line(1,0){0.05}}
\multiput(3.65,3.4)(0.4,0){3}{\line(1,0){0.2}}
\put(4.85,3.4){\line(1,0){0.05}}

\multiput(3.4,1.05)(0,0.4){6}{\line(0,1){0.2}}

\multiput(4.9,1.05)(0,0.4){6}{\line(0,1){0.2}}


\multiput(5.55,3.9)(0.4,0){6}{\line(1,0){0.2}}
\multiput(5.55,5.4)(0.4,0){6}{\line(1,0){0.2}}

\put(5.4,3.9){\line(0,1){0.05}}
\multiput(5.4,4.15)(0,0.4){3}{\line(0,1){0.2}}
\put(5.4,5.35){\line(0,1){0.05}}

\put(7.9,3.9){\line(0,1){0.05}}
\multiput(7.9,4.15)(0,0.4){3}{\line(0,1){0.2}}
\put(7.9,5.35){\line(0,1){0.05}}

\end{picture}
\end{center}

La permutation $\sigma({\scriptstyle\II})$ est alors 
$$\sigma({\scriptstyle\II})=(11, 7, 17, 5, 3, 2, 12, 10, 1, 9, 8, 6,
16, 4, 18, 15, 13, 14).$$ 
Sur le diagramme suivant, on a repr\'esent\'e par des \boldmath $\oplus$
\unboldmath les points du graphe de $\sigma({\scriptstyle\II})$ dont
les ordonn\'ees sont 
dans $\{w(c),w(d),w(a),w(b)\},$ par des $\oplus$ ceux dont les
ordonn\'ees sont dans $\mathcal Y_-\cup \mathcal D\cup \mathcal Y_+,$ et
par des + les autres points du graphe de $\sigma({\scriptstyle\II}).$
On a repr\'esent\'e par 
des $\cdot$ les points du graphe de $w$ dont les abscisses sont dans
$\{a,b,c,d\}\cup \mathcal X_- \cup \mathcal X_+.$ Enfin, $D_{II}$ est la
r\'egion d\'elimit\'ee par les pointill\'es.

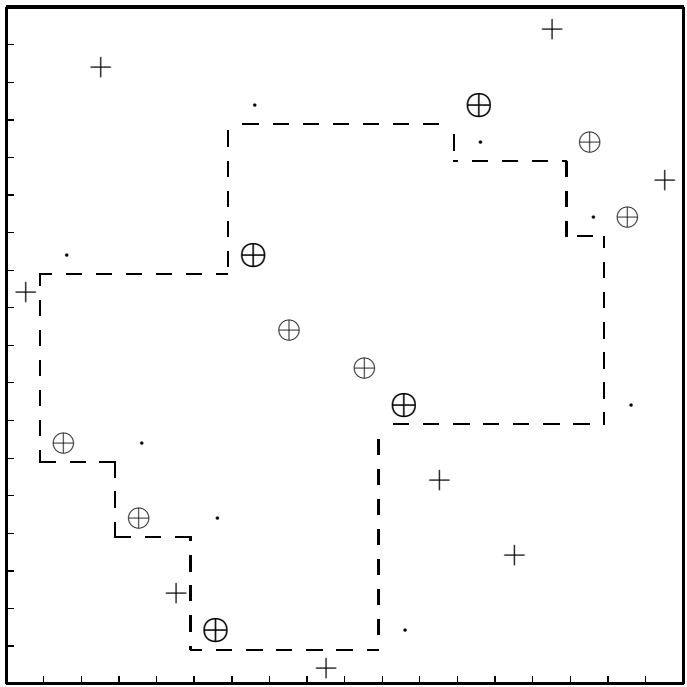
\begin{center}\begin{picture}(10,10)

\boldmath
\put(3,5.5){$\oplus$}
\put(6,7.5){$\oplus$}
\put(2.5,0.5){$\oplus$}
\put(5,3.5){$\oplus$}
\unboldmath

\put(0,5){+}
\put(1,8){+}
\put(0.5,3){$\oplus$}
\put(2,1){+}
\put(1.5,2){$\oplus$}
\put(3.5,4.5){$\oplus$}
\put(4,0){+}
\put(4.5,4){$\oplus$}
\put(5.5,2.5){+}
\put(7.5,7){$\oplus$}
\put(6.5,1.5){+}
\put(7,8.5){+}
\put(8,6){$\oplus$}

\put(8.5,6.5){+}

\thicklines
\put(-0.1,-0.1){\line(1,0){9}}
\put(-0.1,-0.1){\line(0,1){9}}
\put(-0.1,8.9){\line(1,0){9}}
\put(8.9,-0.1){\line(0,1){9}}

\linethickness{0.01mm}
\multiput(0.4,-0.1)(0.5,0){17}{\line(0,1){0.1}}
\multiput(-0.1,0.4)(0,0.5){17}{\line(1,0){0.1}}

\put(0.65,5.5){$\cdot$}
\put(3.15,7.5){$\cdot$}
\put(5.15,0.5){$\cdot$}
\put(8.15,3.5){$\cdot$}

\put(1.65,3){$\cdot$}
\put(2.65,2){$\cdot$}
\put(6.15,7){$\cdot$}
\put(7.65,6){$\cdot$}

\thinlines
\put(0.35,5.35){\line(1,0){0.15}}
\multiput(0.7,5.35)(0.4,0){5}{\line(1,0){0.2}}
\put(2.7,5.35){\line(1,0){0.15}}

\multiput(2.85,5.45)(0,0.4){5}{\line(0,1){0.2}}

\multiput(3.05,7.35)(0.4,0){7}{\line(1,0){0.2}}

\put(5.85,7.2){\line(0,-1){0.2}}

\put(5.85,6.85){\line(1,0){0.05}}
\multiput(6.1,6.85)(0.4,0){3}{\line(1,0){0.2}}
\put(7.3,6.85){\line(1,0){0.05}}

\multiput(7.35,6.85)(0,-0.4){3}{\line(0,-1){0.2}}

\put(7.5,5.85){\line(1,0){0.2}}

\put(7.85,5.85){\line(0,-1){0.15}}
\multiput(7.85,5.5)(0,-0.4){5}{\line(0,-1){0.2}}
\put(7.85,3.5){\line(0,-1){0.15}}

\multiput(7.65,3.35)(-0.4,0){7}{\line(-1,0){0.2}}

\multiput(4.85,3.15)(0,-0.4){7}{\line(0,-1){0.2}}

\put(4.85,0.35){\line(-1,0){0.15}}
\multiput(4.5,0.35)(-0.4,0){5}{\line(-1,0){0.2}}
\put(2.5,0.35){\line(-1,0){0.15}}

\put(2.35,0.35){\line(0,1){0.05}}
\multiput(2.35,0.6)(0,0.4){3}{\line(0,1){0.2}}
\put(2.35,1.8){\line(0,1){0.05}}

\multiput(2.35,1.85)(-0.4,0){3}{\line(-1,0){0.2}}

\multiput(1.35,1.85)(0,0.4){3}{\line(0,1){0.2}}

\multiput(1.35,2.85)(-0.4,0){3}{\line(-1,0){0.2}}

\put(0.35,2.85){\line(0,1){0.15}}
\multiput(0.35,3.2)(0,0.4){5}{\line(0,1){0.2}}
\put(0.35,5.2){\line(0,1){0.15}}
\end{picture}
\end{center}
\end{ex}

Le r\'esultat suivant est une extension du lemme 3.5 de \cite{covex}. 

\begin{lem}\label{longII}$($a$)\
  r_{\sigma({\scriptscriptstyle\II})}=r_w+\chi_{D_{II}}, \mbox{ en 
    particulier } \sigma({\scriptstyle\II}) \le w.$\\ 
  $($b$)\ \ell(\sigma({\scriptstyle\II}))=\ell(w)-(2r+s+t+3).$
\end{lem}

\begin{proof}
$(a)$ L'\'egalit\'e $r_{\sigma({\scriptscriptstyle\II})}=r_w+\chi_{D_{II}}$
r\'esulte de la construction de 
$\sigma({\scriptstyle\II}),$ et l'in\'e\-galit\'e
$\sigma({\scriptstyle\II}) \le w$ en d\'ecoule d'apr\`es le lemme 
\ref{rang}. 

$(b)$ On remarque que $\g({\scriptstyle\II})$ s'\'ecrit aussi 
$$\displaylines{\g({\scriptstyle\II}) =
  \Big[\big(w(d),y_1\big)\big(y_1,y_2\big)\cdots\big(y_{s-1},
   y_s\big)\big(y_s,w(b)\big)\Big]\hfill\cr\hfill   
  \Big[\big(w(a),y_{-1}\big)\big(y_{-1},y_{-2}\big)\cdots\big(
  y_{-t+1},y_{-t}\big)\big( 
  y_{-t},w(c)\big)\Big]\ \big(w(c),w(b)\big).}$$ 
En utilisant le lemme \ref{difflongtransp}, on obtient d'abord
$\ell((w(c),w(b))w)=\ell(w)-(2r+1).$ Ensuite, 
\`a chaque \'etape \'el\'ementaire $v'=(i,j)v$ du passage de $w$ \`a
$\sigma({\scriptstyle\II}),$ 
le rectangle de sommets $(v^{-1}(i),i)$ et $(v^{-1}(j),j)$ ne contient
pas d'autre point du graphe de $v,$ de sorte que l'on a
$\ell(v')=\ell(v)-1.$ On en d\'eduit l'\'egalit\'e annonc\'ee.\end{proof}
%
\subsection{Des composantes de type $S_1$ et $S_2$}\label{descomp}

Nous allons d\'eterminer, parmi les permutations associ\'ees aux
configurations I et II de $w,$ d\'efinies \`a la section pr\'ec\'edente,
celles qui donnent des composantes irr\'eductibles du lieu singulier de
$X_w.$ Pour cela, comme expliqu\'e en \ref{singgen}, nous d\'ecrivons les
transversales correspondantes dans $X_w.$ 

On consid\`ere d'abord une configuration I de $w,$ not\'ee $\mathcal I,$
et l'on note simplement $\tau$ la permutation associ\'ee. On
rappelle que l'on a d\'efini \`a la section \ref{renfor}, 
$$\begin{array}{c}\mathcal {Y}_+=\{y_i\ |\ i\in
  [1,s]\cup\{\infty\}\},\\ 
  \mathcal {X}_-=\{x_i\ |\ i\in
  [-t,-1]\cup\{\infty\}\}.
\end{array}$$

On d\'efinit encore
$$\mathcal{C}_{\mathcal{Y}_+,\mathcal{X}_-}=
  \left\{u=(u_{ij})\in GL_n\left\vert\begin{array}{l}
  u_{\tau(j)j}=1 \mbox{\; pour tout \;} j,\\
  u_{ij}=0 \mbox{\; si \;} i\notin \mathcal{Y}_+\mbox{\; ou \;}
  j\notin \mathcal{X}_-, \\
  \mbox{rg\;}(u_{ij})_{\deuxind{i\in \mathcal{Y}_+}{j\in
  \mathcal{X}_-}}\leq 1
\end{array}\right.\right\}.$$

Le tore $T$ agit sur $\mathcal{C}_{\mathcal{Y}_+,\mathcal{X}_-}$
par : $t\cdot u= t\,u\,(\overline{\tau}^{-1}\,t^{-1}\,\overline{\tau}),$ o\`u 
$\overline{\tau}$ d\'esigne la matrice de la permutation $\tau,$
c'est-\`a-dire la matrice dont les coefficients sont les
$\delta_{i\,\tau(j)}.$ La projection
de $\mathcal{M}_n$ (ensemble des matrices carr\'ees d'ordre $n$ \`a
coefficients dans $\K$) sur $\mathcal{M}_{s+1,t+1}$ (ensemble des
matrices de taille $(s+1,t+1)$ \`a coefficients dans $\K$) obtenue par omission
des lignes (resp. colonnes) d'indice n'appartenant pas \`a
$\mathcal{Y}_+$ (resp. $\mathcal{X}_-$) induit un isomorphisme
$T$-\'equivariant de $\mathcal{C}_{\mathcal{Y}_+,\mathcal{X}_-}$ sur  
$\mathcal{C}_{s+1,t+1}.$ 

On d\'emontre alors, exactement comme dans \cite{covex}, th\'eor\`eme 3.6, le 

\begin{thm}\label{gencomp}$($a$\,)$ L'application $u\longmapsto u\K^{\bullet}$
  induit un isomorphisme $T$-\'equivariant de
  $\mathcal{C}_{\mathcal{Y}_+,\mathcal{X}_-}$ sur
  $\mathcal{N}_{\tau , w}.$ \\
  $($b$\,)$ Par cons\'equent, si $st=0,\ e_{\tau}$ est un point lisse de
  $X_w,$ et si $st\not=0,\ X_{\tau}$ est une composante irr\'eductible
  de Sing $X_w,$ de type $S_1.$
\end{thm}

\vspace*{1em}

On consid\`ere maintenant une configuration II de $w,$ et la
permutation associ\'ee, not\'ee simplement $\sigma.$ En vue de la description de la
transversale $\mathcal{N}_{\sigma , w},$ il est utile de consid\'erer la
vari\'et\'e suivante : \'etant donn\'e trois entiers $i,j,k$ avec $i,k \ge 1$
et $j\ge 2,$ on d\'efinit 

$$\mathcal{N}_{i,j,k}=\left\{(M,N)\in \mathcal{C}_{j,k}\times
  \mathcal{C}_{i,j}\ |\ NM=0\right\}.$$

\begin{prop}\label{modele} $\mathcal{N}_{i,j,k}$ est une vari\'et\'e
  irr\'eductible de dimension $2j+i+k-3.$ Elle est singuli\`ere au point
  $(0,0).$ \\
  $($a$\,)$ Cette singularit\'e est isol\'ee si et seulement si
  $j=2$ ou $i=k=1.$\\
  $($b$\,)$ Si $i>1,\ j>2,$ et $k=1\ ($resp. $i=1,\ j>2,$ et $\ k>1)$
  alors Sing $\mathcal{N}_{i,j,1}=\mathcal{C}_{j,1}\times\{0\}\ 
  ($resp. Sing $\mathcal{N}_{1,j,k}=\{0\}\times\mathcal{C}_{1,j}).$\\
  $($c$\,)$ Si $\,i,k>1$ et $j>2,$ Sing $\mathcal{N}_{i,j,k}$ a deux
  composantes irr\'eductibles, $\mathcal{C}_{j,k}\times\{0\}$ et
  $\{0\}\times\mathcal{C}_{i,j}.$\\
  $($d$\,)$ On a :\\
  \hspace*{0.7cm}$\bullet\;\mathcal{N}_{i,2,k}$ est isomorphe \`a
  $\mathcal{C}_{i+k,2},$\\
  \hspace*{3ex} $\bullet\;\mathcal{N}_{1,j,1}$ est un c\^one quadratique
  non d\'eg\'en\'er\'e de dimension $2j-1.$\\
\end{prop}

\begin{proof} 
Les assertions $(a\,),\  (b\,)$ et $(c\,)$ de cette proposition sont
  un cas particulier du Th\'eor\`eme 1 de \cite{Gon} ; on en donne ici, dans
  ce cas particulier, une d\'emonstration directe, plus simple et plus
  g\'eom\'etrique.  

Notons $\P(\K^j)$ (resp. $\P^*(\K^j)$) l'espace projectif des droites
(resp. hyperplans) dans $\K^j.$ En consid\'erant la r\'esolution 
$$\displaylines{\hspace*{3em}\mathcal{Z}_{i,j,k}=\{(M, \mathcal{D}, \mathcal{H}, N)\in
  \mathcal{C}_{j,k}\times \P(\K^j)\times \P^*(\K^j)\times
  \mathcal{C}_{i,j}\  |\hfill\cr\hfill \mbox{Im } M\subseteq \mathcal{D}
  \subseteq \mathcal{H}\subseteq \mbox{Ker }N \},\hspace*{3em}}$$ on voit que
$\mathcal{N}_{i,j,k}$ est irr\'eductible et de dimension $2j+i+k-3.$

D'autre part, on voit facilement que l'espace tangent \`a
$\mathcal{N}_{i,j,k}$ en $(0,0)$ s'identifie \`a
$\mathcal{M}_{j,k}\times \mathcal{M}_{i,j}.$ Par cons\'equent, le point
$(0,0)$ est singulier dans $\mathcal{N}_{i,j,k}.$ Par 
ailleurs, notant $m_{pq}$ les coefficients de la matrice $M,$ on
d\'efinit les ouverts affines 
$$U_{pq}=\{M\in\mathcal{C}_{j,k}\ |\ m_{pq}\not=0\}$$
et $$\mathcal U_{pq}=\{(M,N)\in \mathcal{N}_{i,j,k}\ |\
m_{pq}\not=0\}.$$
On voit sans peine que, premi\`erement, l'application qui \`a $M$ associe
les coefficients $m_{pq},\ m_{rq},\ m_{ps}$ pour $r\not=p$ et
$s\not=q,$ induit un isomorphisme $U_{pq}\simeq \K^*\times \K^{j+k-2},$
et deuxi\`emement, l'application
$(M,N)\longmapsto (M,\overline{N}),$ o\`u $\overline{N}$ est la matrice
obtenue par omission de la colonne d'indice $p$ de $N,$ induit un
isomorphisme $\mathcal U_{pq}\simeq U_{pq}\times \mathcal
C_{i,j-1}.$ 

De m\^eme, on d\'efinit 
$$\mathcal V_{pq}=\{(M,N)\in \mathcal{N}_{i,j,k}\ |\ n_{pq}\not=0\}$$ et
$$V_{pq}=\{N\in\mathcal{C}_{i,j}\ |\ n_{pq}\not=0\}.$$ On obtient
comme pr\'ec\'edemment que $V_{pq}\simeq \K^*\times \K^{i+j-2}$ et que $\mathcal
V_{pq} \simeq \mathcal C_{j-1,k}\times V_{pq}.$

D'autre part, on observe que $\mathcal C_{r,s}$
est lisse si $r=1$ ou $s=1$ et a $0$ pour unique point singulier
sinon. Comme les ouverts $\mathcal U_{pq}$ et $\mathcal V_{pq}$
recouvrent $\mathcal{N}_{i,j,k}\setminus\{(0,0)\},$ on en d\'eduit les
assertions $(a\,),\ (b\,)$ et $(c\,).$ Il reste \`a \'etablir les deux
isomorphismes de $(d\,).$  

D'abord, il est clair que si $i=k=1,\;\mathcal{N}_{1,j,1}$ est un c\^one
quadratique non d\'eg\'en\'er\'e de dimension $2j-1.$ Supposons maintenant
$j=2$ ; alors  
$$\mathcal{N}_{i,2,k}=\{(M,N)\in
\mathcal C_{2,k}\times \mathcal C_{i,2}\ | \  N M=0 \}.$$
On v\'erifie sans difficult\'e que l'application $${\begin{array}[t]{ccl}
   \mathcal M_{2,k}\times\mathcal M_{i,2} & \longrightarrow &
   \mathcal M_{i+k,2}, \\ 
   (M,N)& \longmapsto & \left(\begin{array}{c}
   N\\\widetilde{M}\end{array}\right) 
\end{array}}$$ o\`u
$\widetilde{M}=\  ^tM\left(\begin{array}{cc}0&-1\\1&0\end{array}\right),$
induit un 
isomorphisme de $\mathcal{N}_{i,2,k}$ sur $\mathcal C_{i+k,2}.$ Cela
ach\`eve la preuve de la proposition.
\end{proof}

On passe maintenant \`a la description de la transversale
$\mathcal{N}_{\sigma , w}$ ; cela n\'ecessite encore quelques
notations. On d\'efinit 
$$\begin{array}{cc}\overline{\mathcal X_-}=\mathcal X_-
  \cup\{a\}, & \overline{\mathcal Y_+}=\mathcal Y_+\cup\{w(b)\},\\
  \overline{\mathcal C}=\mathcal C\cup\{b,c\},&
  \overline{\mathcal D}=\mathcal D \cup\{w(a),w(d)\},
\end{array}$$
et l'on d\'esigne par $\mathcal M_{\sigma}$ l'ensemble des $u\in GL_n$
tels que : \\
\hspace*{3ex}$\bullet \; u_{\sigma(j)j}=1$ pour tout $j,$ \\
\hspace*{3ex}$\bullet \; u_{ij}=0$ si $i\not=\sigma(j)$ et $(i,j)\not\in
  (\overline{\mathcal D}\times \overline{\mathcal X_-})\cup
  (\overline{\mathcal Y_+}\times \overline{\mathcal C}),$\\ 
\hspace*{3ex}$\bullet \;
(M,\overleftarrow{N})\in\mathcal{N}_{s+1,r+2,t+1},$ o\`u 
  $M$ et $N$ sont les matrices extraites de $u$ d\'efinies par
  $M=(u_{ij})_{(i,j)\in\overline{\mathcal D}\times\overline{\mathcal
  A_-}}$ et $N=(u_{ij})_{(i,j)\in\overline{\mathcal Y_+}\times
  \overline{\mathcal C}},$ et $\overleftarrow{N}$ est la matrice
  obtenue en lisant $N$ de droite \`a gauche. 

On fait agir le tore $T$ sur $\mathcal M_{\sigma}$ par $t\cdot u=
t\,u(\overline{\sigma}^{-1}\,t^{-1}\,\overline{\sigma}).$  

\begin{thm}\label{trans}L'application $u\longmapsto u\K^{\bullet}$
  induit un isomorphisme $T$-\'equiva\-riant de $\mathcal M_{\sigma}$
  sur $\mathcal{N}_{\sigma , w}.$
\end{thm} 

\begin{proof}
La preuve de ce th\'eor\`eme s'obtient en \'etendant les arguments de la
preuve du th\'eor\`eme 3.6 de \cite{covex}. 
Rappelons que $\mathcal{N}_{\sigma , w}=\big(\sigma(U^-)\cap U^- \big)
e_{\sigma} \bigcap X_w,$ et qu'il r\'esulte de la d\'ecomposition de
Bruhat que l'application 
$$\begin{array}{rccc}\phi : \;&\overline{\sigma} U^- \cap
  U^-\overline{\sigma} & 
\longrightarrow & \big(\sigma(U^-)\cap U^-\big) e_{\sigma}\\
 &u & \longmapsto & u\K^{\bullet}\end{array}$$
est un isomorphisme $T$-\'equivariant. Il s'agit donc de montrer que
$\mathcal M_{\sigma}=\phi^{-1}(\mathcal{N}_{\sigma , w}).$

On a d'abord $$\overline{\sigma} U^- \cap U^-\overline{\sigma}=
\left\{u\in G\!L_n\ \left\vert\begin{array}{l} 
u_{\sigma(j)j}=1 \mbox{\; pour tout \;} j,\\
u_{ij}=0 \mbox{\;si \;}i<\sigma(j) \mbox{\;ou \;}
j>\sigma^{-1}(i)\end{array}\right.\right\}.$$

Consid\'erons maintenant $u\in \phi^{-1}(\mathcal{N}_{\sigma , w}),$ et
montrons que pour $(i,j)\not\in (\overline{\mathcal D}\times 
\overline{\mathcal X_-})\cup (\overline{\mathcal Y_+}\times
\overline{\mathcal C}),$ avec $i>\sigma(j)$ et $j<\sigma^{-1}(i),$ on a
$u_{ij}=0.$ 

Pour commencer, si $j\not\in\overline{\mathcal
  X_-}\cup\overline{\mathcal C},$ le point $(j,\sigma(j))$ n'est pas
dans $D_{II},$ on a donc
  $r_{\sigma}\big(j,\sigma(j)\big)=r_w\big(j,\sigma(j)\big).$ L'espace
$u(\K^j)+\K^{\sigma(j)}$ contient la famille de vecteurs
$$\{e_1,\ldots,e_{\sigma(j)}\}\cup\{ue_p\ |\ p<j,
\sigma(p)>\sigma(j)\}\cup\{ue_j-e_{\sigma(j)}\}.$$
S'il existait $i>\sigma(j)$ tel que $u_{ij}\not=0,$ alors on aurait $i\notin
\{\sigma(p) \ |\ p<j\}$ car $u_{\sigma(p)j}=0$ pour $j<p,$ et donc cette
famille serait libre. Or son cardinal est 
$$\sigma(j)+j-r_{\sigma}\big(j,\sigma(j)\big)+1=
\sigma(j)+j-r_w\big(j,\sigma(j)\big)+1,$$
alors que $\dim(u(\K^j)+\K^{\sigma(j)})\le
\sigma(j)+j-r_w\big(j,\sigma(j)\big),$ une contradiction. Ainsi, si
$j\not\in\overline{\mathcal X_-}\cup\overline{\mathcal C},$ on a
$u_{ij}=0$ pour tout $i>\sigma(j).$ De m\^eme, on montre que si
$i\not\in \overline{\mathcal Y_+}\cup \overline{\mathcal D}$ alors
$u_{ij}=0$ pour tout $j<\sigma^{-1}(i).$ 

\smallskip

Remarquons de plus que si $(i,j)\in \overline{\mathcal
D}\times\overline{\mathcal C}, $ avec $\sigma(j)\not=i,$ alors on a
$i<\sigma(j)$ ou bien $j>\sigma^{-1}(i),$ car la restriction de
$\sigma$ \`a $\overline{\mathcal
C}\times\overline{\mathcal D}$ est l'\'el\'ement de plus grande longueur ; 
et donc $u_{ij}=0,$ d'apr\`es 
la description de $\overline{\sigma} U^- \cap U^-\overline{\sigma}.$
Par ailleurs, on voit que si $(i,j)\in \overline{\mathcal Y_+}\times 
\overline{\mathcal X_-},$ le point $(j,i-1)$ n'est pas dans $D_{II}$ ; en
consid\'erant l'espace $u(\K^j)+\K^{i-1},$ on montre alors comme
pr\'ec\'edemment que $u_{ij}=0.$ Cela prouve que $u_{ij}=0$ si
$i\not=\sigma(j)$ et $(i,j)\not\in (\overline{\mathcal D}\times 
\overline{\mathcal X_-})\cup (\overline{\mathcal Y_+}\times
\overline{\mathcal C}).$

Soient maintenant $\mu$ le rang de la matrice extraite 
$M=(u_{ij})_{(i,j)\in\overline{\mathcal D}\times 
\overline{\mathcal X_-}},$ et $\nu$ le rang de la matrice extraite
$N=(u_{ij})_{(i,j)\in 
\overline{\mathcal Y_+}\times \overline{\mathcal C}}.$ Soit $E$
l'espace  engendr\'e par la famille de vecteurs 
$$\{e_1,\ldots , e_{w(d)-1}\}\cup\{ue_j\ |\ j<b,
\sigma(j)>w(d)\}\cup\{ue_j\ |\ j\in \overline{\mathcal X_-}\}.$$
On voit que sa dimension est 
$$(w(d)-1)+(b-1)-r_{\sigma}\big(b-1,w(d)\big)+\mu.$$
Le point $(b-1, w(d))$ est dans $D_{II},$ donc on a 
$$\dim E=(w(d)-1)+(b-1)-r_w\big(b-1,w(d)\big)+\mu-1.$$
Or $E\subseteq u(\K^{b-1})+\K^{w(d)-1},$ donc $\dim E\le
(w(d)-1)+(b-1)-r_w\big(b-1,w(d)\big).$ Il en r\'esulte que $\mu \le 1.$
On montre de m\^eme que $\nu\le 1.$ 

Montrons maintenant que $\overleftarrow N M=0.$ Pour cela, on
introduit les notations suivantes : on pose $c_0=b,
c_{r+1}=c, d_0=w(a), \mbox{ et } d_{r+1}=w(d).$ Il
s'agit alors de montrer que pour tout $(i,j)\in \overline{\mathcal
  Y_+}\times \overline{\mathcal X_-}$ 
$$\sum_{k=0}^{r+1}u_{ic_k}u_{d_k j}=0.$$
Fixons $j\in\overline{\mathcal X_-}$ ; d'apr\`es ce qui pr\'ec\`ede,
$$ue_j=e_{\sigma(j)}+\sum_{k=0}^{r+1}u_{d_k j}e_{d_k},$$ et
pour tout $k\in [0,r+1],$ $$ue_{c_k}= e_{d_k}+
\sum_{i\in\overline{\mathcal Y_+}} u_{ic_k}e_i.$$
Soit $f_j=ue_j- e_{\sigma(j)}-\sum_{k=0}^{r+1}u_{d_k
  j}ue_{c_k}.$ Ce vecteur appartient \`a l'espace
$u(\K^c)+\K^{w(d)-1},$ et il s'\'ecrit en fait 
$$f_j=-\sum_{i\in\overline{\mathcal Y_+}}\left(\sum_{k=0}^{r+1}
  u_{ic_k}u_{d_k j}\right)e_i.$$
S'il existait $j\in\overline{\mathcal X_-}$ tel que $f_j\not=0,$ alors
la famille 
$$\{e_1,\ldots,e_{w(d)-1}\}\cup\{ue_p\ |\ p\le c,\ \sigma(p)\ge
w(d)\}\cup\{f_j\}$$ serait libre. Son cardinal est 
$w(d)-1+c-r_{\sigma}(c,w(d)-1)+1,$ et elle est contenue dans l'espace
$u(\K^c)+\K^{w(d)-1},$ dont la dimension est au plus
$w(d)-1+c-r_w(c,w(d)-1).$ Or le point $(c,w(d)-1)$ n'est pas
dans $D_{II}$ ; on a donc $r_{\sigma}(c,w(d)-1)=r_w(c,w(d)-1),$ une
contradiction.

Ainsi, on a montr\'e que $\phi^{-1}(\mathcal{N}_{\sigma ,
w})\subseteq\mathcal M_{\sigma}.$ Or d'apr\`es le lemme \ref{longII} et la
proposition \ref{modele}, ce sont deux vari\'et\'es
irr\'eductibles de dimension $2r+t+s+3$ ; on en d\'eduit donc
$\phi^{-1}(\mathcal{N}_{\sigma , w})=\mathcal M_{\sigma}.$ Cela prouve
le th\'eor\`eme \ref{trans}.
\end{proof}

\begin{coro}\label{CompII} $($a$\,)\  e_{\sigma}$ est un point singulier
  de $X_w.$\\ 
  $($b$\,)$ Si la configuration est mixte ({\it i.e.} si $r=0$),
  $X_{\sigma}$ est une composante 
  irr\'eductible de type $S_1$ de Sing $X_w,$ la transversale \'etant
  isomorphe \`a $\mathcal{C}_{s+t+2,2}.$\\
  $($c$\,)$ Si la configuration est pure ({\it i.e.} si $r\not =0$ et
  $s=t=0$), $X_{\sigma}$ est une composante
  irr\'eductible de type $S_2$ de Sing $X_w,$ la transversale \'etant
  isomorphe \`a $\mathcal K_{2r+3}.$\\
  $($d$\,)$ Si $rt\not=0$ et $s=0\ ($resp. $rs\not=0$ et $t=0)$ alors
  $X_{\sigma}$ est contenu dans exactement une composante irr\'eductible
  $X_{\tau}$ de Sing $X_w,$ associ\'ee \`a une configuration I de $w,$ et
  telle que $\mathcal N_{\tau,w}\simeq \mathcal C_{r+1,t+1}\ 
  ($resp. $\mathcal N_{\tau,w}\simeq \mathcal C_{s+1,r+1}).$\\
  $($e$\,)$ Si $rst\not=0,$ alors $X_{\sigma}$ est contenu dans exactement
  deux composantes irr\'eductibles de Sing $X_w,$ correspondant \`a
  des configurations I de $w,$ et dont les transversales sont
  isomorphes \`a $\mathcal C_{r+1,t+1}$ et $\mathcal C_{s+1,r+1}.$\\
\end{coro}

\begin{proof}
Les trois premiers points r\'esultent directement de la proposition
\ref{modele} et du th\'eor\`eme \ref{trans}. 

Supposons maintenant $rt\not=0.$ Alors les points $(a,w(a))$ et $(c,w(c)),$
la suite SO et la suite centrale forment une configuration I de
$w,$ non d\'eg\'en\'er\'ee. La permutation $\tau$ associ\'ee donne donc une
composante irr\'eductible du lieu singulier de $X_w.$ 

De plus, les graphes de $\tau$ et $\sigma$ ne diff\`erent qu'en les
points dont les abscisses sont dans $\overline{\mathcal C}\cup
\mathcal X_+\cup \{d\}$ ; leurs ordonn\'ees sont dans
$\overline{\mathcal D}\cup\overline{\mathcal Y_+}.$ Or sur
ces points, $\tau$ induit la permutation maximale telle que
$\tau(\overline{\mathcal C})\subseteq \overline{\mathcal D}\cup
\{w(b)\}.$ Comme $\sigma(\overline{\mathcal C})=\overline{\mathcal
  D},$ on obtient donc, en vertu du lemme \ref{focalisation}, $\sigma
\le \tau.$ On remarque de plus que $\ell(\tau)=\ell(w)-(r+t+1).$ 

De m\^eme, si $rs\not=0,$ les points $(b,w(b))$ et $(d,w(d)),$ la suite
centrale et la suite NE forment une configuration I non d\'eg\'en\'er\'ee de
$w.$ La permutation $\tau'$ associ\'ee donne une composante irr\'eductible
de Sing $X_w,$ et l'on a $\sigma\le \tau'.$ Elle v\'erifie de plus
$\ell(\tau')=\ell(w)-(r+s+1).$  

Par ailleurs, on observe que si $rst\not=0,$ les permutations $\tau$
et $\tau'$ sont distinctes. Or, d'apr\`es la proposition \ref{modele} et
le th\'eor\`eme \ref{trans}, $e_{\sigma}$ est contenu dans exactement une
composante irr\'eductible de Sing $X_w$ lorsque $rt\not=0$ et $s=0$ ou
lorsque $rs\not=0$ et $t=0,$ et dans exactement deux lorsque
$rst\not=0.$ Cela ach\`eve la d\'emonstration du corollaire. 
\end{proof}

\begin{rem}
Soit $\mathcal I$ une configuration I d\'eg\'en\'er\'ee de $w.$ D'apr\`es le
th\'eor\`eme \ref{gencomp}, le point $e_{\tau({\scriptscriptstyle\mathcal
I})}$ est lisse dans $X_w.$ On peut voir qu'il existe une configuration
(3412) incompressible $\mathcal I\!\mathcal I$ telle que la permutation
associ\'ee $\sigma({\scriptstyle\mathcal I\!\mathcal I})$ corresponde \`a
une composante irr\'eductible de Sing $X_w,$ et v\'erifie 
$\tau({\scriptstyle\mathcal I})> \sigma({\scriptstyle\mathcal
I\!\mathcal I}).$  
\end{rem}

%

\section{Quasi-r\'esolutions des vari\'et\'es de Schubert}\label{quasiresol}
\index{$\Sigma_w$}
\begin{defn}On note $\Sigma_w$ la r\'eunion des composantes du lieu
  singulier exhib\'ees \`a la section pr\'ec\'edente. 
\end{defn}

On consid\`ere une permutation $w$ non covexillaire. L'objet des deux
derni\`eres sections est de d\'emontrer que 
$\Sigma_w=\mbox{ Sing }X_w.$ A cet effet, nous allons introduire des
quasi-r\'esolutions des vari\'et\'es de Schubert non covexillaires,
c'est-\`a-dire des morphismes birationnels $P\times^Q X_y
\longrightarrow X_w$ pour certains sous-groupes paraboliques $P\supseteq Q$ de
$G$ et certaines vari\'et\'es de Schubert $X_y\subseteq X_w.$ (Rappelons
que, pour une vari\'et\'e alg\'ebrique $X,$ munie d'une action 
d'un sous-groupe ferm\'e $Q$ d'un groupe alg\'ebrique $P,$ on note
$P\times^Q X$ le quotient du produit $P\times X$ sous l'action
diagonale de $Q$ : $q\cdot (p,x)=(pq^{-1},qx)$ pour tous $q\in Q, p\in
P$ et $x\in X$ ; voir \cite{Ser}, et aussi \cite{Bia}. Une propri\'et\'e
cruciale de cette construction est que le morphisme $P\times^Q
X\longrightarrow P/Q$ est une fibration localement triviale, de fibre $X.$)  

On rappelle par ailleurs que, \'etant donn\'e un morphisme birationnel $\pi : X
\longrightarrow Y$ entre des vari\'et\'es irr\'eductibles, on
d\'efinit l'ouvert $\Reg(\pi)$ \index{$\Reg(\pi)$} de $X$ comme
l'ensemble des points 
admettant un voisinage ouvert $U$ tel que $\pi$ induise un
isomorphisme de $U$ sur $\pi(U).$ Le lieu exceptionnel de $\pi,\
\Ex(\pi),$ \index{$\Ex(\pi)$} est le ferm\'e compl\'ementaire de
$\Reg(\pi).$ D'autre part, 
on peut d\'efinir dans $Y$ l'ouvert $\Reg_{\pi}$ \index{$\Reg_{\pi}$}
comme l'ensemble des 
points admettant un voisinage $V$ tel que $\pi$ induise un
isomorphisme de $\pi^{-1}(V)$ sur $V.$ Alors $\pi$ induit un
isomorphisme de $\Reg(\pi)$ sur $\Reg_{\pi}.$ Si de plus $\pi$ est
surjectif, notant $\EE_{\pi}=\pi(\Ex(\pi)),\ \EE_{\pi}$
\index{$\EE_{\pi}$} est le ferm\'e compl\'ementaire de $\Reg_{\pi},$ et
l'on a aussi $\pi^{-1}(\EE_{\pi})=\Ex(\pi).$    
 
Nous \'etudierons les lieux exceptionnels des quasi-r\'esolutions, et
leurs images, que nous relierons finalement \`a $\Sigma_w.$

\begin{defn}Soit $\II$ une configuration (3412) de $w,$
  correspondant aux abscisses $a<b<c<d.$ Notons $\a=w(a), \b=w(b),
  \g=w(c), \d=w(d).$\\
\hspace*{3ex} $\bullet \;$ On associe \`a $\II$ la paire d'entiers suivante
  : sa hauteur $h(\II)=\a-\d$ et son amplitude $am(\II)=\b-\g.$\\
\hspace*{3ex} $\bullet \;$ On dit que $\II$ est
  {\it bien remplie} si $\winv\big(]\d,\a[\big)\subseteq\  ]b,c[.$
\end{defn}

Remarquons qu'il existe des configurations (3412) bien remplies : on
voit facilement qu'une configuration (3412) de hauteur minimale est
bien remplie. 

On fixe une configuration (3412) de $w,$ correspondant aux abscisses
$a<b<c<d,$ et aux ordonn\'ees $\g<\d<\a<\b,$ bien remplie
et d'amplitude minimale (parmi les configurations bien remplies). On
v\'erifie alors sans peine qu'elle est incompressible. On note
simplement $h$ sa hauteur.

On consid\`ere 
$$\begin{array}{l}
\a' =\max \{q\ge \a\ |\ \forall q'\in [\a,q[,\ \winv(q'+1)<\winv(q')\}\\
\d' =\min \{q\le \d\ |\ \forall q'\in\, ]q,\d],\
\winv(q'-1)>\winv(q')\}.\end{array}$$

Soit $I=\{s_{\d'},\ldots,s_{\a'-1}\},$ et pour $i=1,\ldots, h,$
soient $k_i=\d'+\a'-\a+i-1$ et $J_i=I\setminus\{s_{k_i}\}.$ Notant
$w_I$ et $w_{J_i}$ les permutations maximales des sous-groupes
paraboliques de $\mathfrak{S}_n$ correspondants, on d\'efinit
$w_i=w_{J_i}w_Iw.$ Comme $w$ est maximal dans 
sa classe $\mathfrak S_I w$, $w_i$ est maximal dans sa classe
$\mathfrak S_{J_i}w_i.$  

\begin{ex}\label{graphe w_i}

On consid\`ere la permutation $$w=(8,7,9,6,1,11,5,4,2,10,3)\in \mathfrak
S_{11}$$ dont le graphe est repr\'esent\'e sur le premier diagramme
ci-dessous. Les points d'abscisses 2,3,9,11 forment une configuration
(3412) bien remplie et d'amplitude minimale (repr\'esent\'es par des
$\oplus$ sur le diagramme). Cette configuration est de hauteur 4 ; on
a $\d=3=\d', \a=7, \a'=8.$  

Le second diagramme est le graphe de la permutation $w_3$ associ\'ee ;
on a $k_3=6.$ 
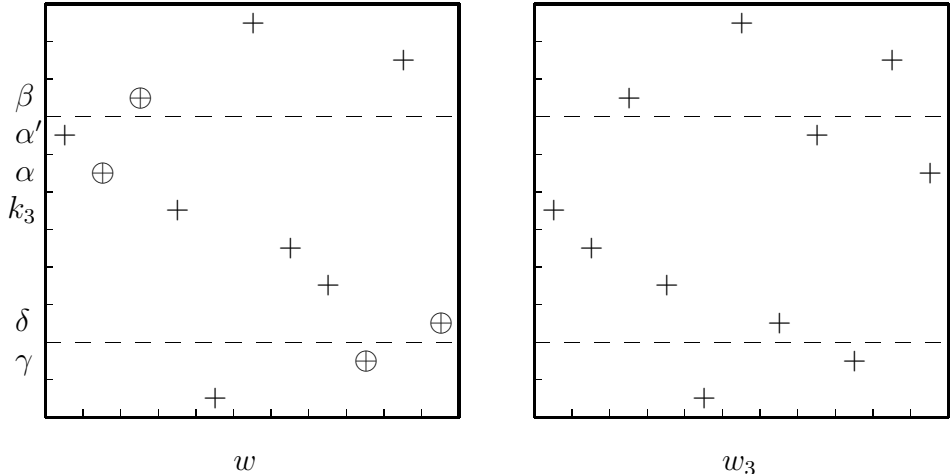
\begin{center}
\begin{picture}(13,6.5)

\thicklines
\put(0.5,0){\line(1,0){5.5}}
\put(0.5,5.5){\line(1,0){5.5}}
\put(0.5,0){\line(0,1){5.5}}
\put(6,0){\line(0,1){5.5}}

\put(7,0){\line(1,0){5.5}}
\put(7,5.5){\line(1,0){5.5}}
\put(7,0){\line(0,1){5.5}}
\put(12.5,0){\line(0,1){5.5}}

\linethickness{0.01mm}
\multiput(1,0)(0.5,0){10}{\line(0,1){0.1}}
\multiput(0.5,0.5)(0,0.5){10}{\line(1,0){0.1}}

\multiput(7.5,0)(0.5,0){10}{\line(0,1){0.1}}
\multiput(7,0.5)(0,0.5){10}{\line(1,0){0.1}}

\thinlines

\put(0.1,0.65){$\g$}
\put(0.1,1.15){$\d$}
\put(0,2.65){$k_3$}
\put(0.1,3.15){$\a$}
\put(0.1,3.65){$\a'$}
\put(0.1,4.15){$\b$}

\put(3,-0.7){$w$}
\put(9.5,-0.7){$w_3$}


\linethickness{0.01mm}
\multiput(0.5,1)(0.4,0){14}{\line(1,0){0.2}}
\multiput(0.5,4)(0.4,0){14}{\line(1,0){0.2}}

\multiput(7,1)(0.4,0){14}{\line(1,0){0.2}}
\multiput(7,4)(0.4,0){14}{\line(1,0){0.2}}
\thinlines
\multiput(2.6,0.15)(6.5,0){2}{$+$}
\multiput(3.1,5.15)(6.5,0){2}{$+$}
\multiput(5.1,4.65)(6.5,0){2}{$+$}

\put(0.6,3.65){$+$}
\put(1.1,3.15){$\oplus$}
\put(1.6,4.15){$\oplus$}
\put(2.1,2.65){$+$}
\put(3.6,2.15){$+$}
\put(4.1,1.65){$+$}
\put(4.6,0.65){$\oplus$}
\put(5.6,1.15){$\oplus$}

\put(7.1,2.65){$+$}
\put(7.6,2.15){$+$}
\put(8.1,4.15){$+$}
\put(8.6,1.65){$+$}
\put(10.1,1.15){$+$}
\put(10.6,3.65){$+$}
\put(11.1,0.65){$+$}
\put(12.1,3.15){$+$}

\end{picture}\end{center}
\vspace*{1cm}
\end{ex}

On d\'efinit enfin $$Z_i=P_I\times^{P_{J_i}}X_{w_i},$$ et l'on note
$\pi_i$ l'application naturelle de $Z_i$ dans $G/B,$ d\'efinie par
$\pi_i([p,x])=px,$ o\`u $[p,x]$ est la classe du couple $(p,x)$ de
$P_I\times X_{w_i}.$ 

\begin{prop}\label{prop-pi}
L'application $\pi_i$ a pour image $X_w,$ et la projection induite
$\pi_i\ :\ Z_i\longrightarrow X_w$ est birationnelle.
\end{prop}

\begin{proof}L'image de $\pi_i$ est $P_I X_{w_i},$ donc elle est
  contenue dans $X_w.$ De plus $\pi_i$ est un morphisme propre, car
  c'est le compos\'e de l'injection naturelle de 
  $P_I\times^{P_{J_i}}X_{w_i}$ dans $P_I\times^{P_{J_i}} G/B$ induite
  par l'inclusion de $X_{w_i}$ dans $G/B,$ suivi de l'isomorphisme
 $\begin{array}[t]{ccc}P_I \times^{P_{J_i}} G/B & \longrightarrow &
 P_I/P_{J_i}\times G/B,\\
{[}p,x{]} & \longmapsto & (p P_{J_i}, p x)
\end{array}$ et enfin de la deuxi\`eme projection. L'image de $\pi_i$
  est donc ferm\'ee. Or elle contient $P_I e_{w_i},$ dense dans $X_w,$
  on a donc $\pi_i(Z_i)=X_w.$ 

Soit maintenant $z_i=[w_Iw_{J_i},e_{w_i}]\in Z_i,$ l'orbite $Uz_i$ est
un ouvert dense de $Z_i,$ et il r\'esulte de la d\'ecomposition de Bruhat
que $\pi_i$ induit un isomorphisme $Uz_i \longrightarrow U e_w,$ ainsi
$\pi_i$ est birationnelle. 
\end{proof}

Remarquons n\'eanmoins que la vari\'et\'e $Z_i$ n'est en g\'en\'eral pas
lisse, on l'appelle donc quasi-r\'esolution de $X_w.$ On peut d\'ecrire
$Z_i$ de mani\`ere plus explicite, via l'isomorphisme 
$P_I \times^{P_{J_i}} G/B \longrightarrow 
 P_I/P_{J_i}\times G/B$ \'evoqu\'e plus haut ; l'image de $Z_i$ est 
$$\displaylines{\hspace*{3em}
\{(U^{k_i},V^{\bullet})\in Gr_{k_i}(n)\times X_w\,\vert\,
\K^{\d'-1}\subseteq U^{k_i} \subseteq \K^{\a'} \mbox{ et } \hfill\cr\hfill\dim
(V^p\cap U^{k_i})\ge r_{w_i}(p,k_i) \mbox{ pour tout } p\}\hspace*{3em}}$$
o\`u $Gr_{k_i}(n)$ d\'esigne la grassmannienne des $k_i$-plans dans
$\K^n.$

On va d\'ecrire l'image du lieu exceptionnel de chaque $\pi_i$ et \'etablir les
deux propositions suivantes :

\begin{prop}\label{inter}
L'intersection des images des lieux exceptionnels des $\pi_i,$ pour $i$
parcourant l'intervalle $[1,h],$ est contenue dans $\Sigma_w.$
\end{prop}

\begin{prop}\label{corr-config}
Pour toute configuration $\mathcal K$ de $w_i,$ de type I ou II,
param\'etrant une composante irr\'eductible $X_v$ du lieu singulier de
$X_{w_i},$ on a :\\ 
\hspace*{3ex}$\bullet$ ou bien $P_I\times^{P_{J_i}}X_v\subseteq
\Ex\,(\pi_i),$\\
\hspace*{3ex}$\bullet$ ou bien $w_Iw_{J_i}(\mathcal K)$ est une
configuration du m\^eme type de $w,$ et 
$\pi_i(P_I\times^{P_{J_i}}X_v) = X_{w_Iw_{J_i}v}$ est la composante
irr\'eductible du lieu singulier de $X_w$ associ\'ee  
(o\`u $w_Iw_{J_i}(\mathcal K)$ d\'esigne l'ensemble des points $(x,w(x))$ tels que
$(x,w_i(x))\in\mathcal K$).
\end{prop}

Nous renvoyons les d\'emonstrations de ces propositions \`a la section
suivante, pour donner d\`es \`a pr\'esent la preuve du 

\begin{thmppal}
Le lieu singulier de $X_w$ est la r\'eunion des composantes d\'ecrites
dans la section {\rm\ref{descomp}}. En particulier, les singularit\'es
g\'en\'eriques sont de type $S_1$ ou $S_2,$ c'est-\`a-dire soit un c\^one de
matrices de rang au plus 
$1,$ soit un c\^one quadratique non d\'eg\'en\'er\'e de dimension impaire $d\ge 5.$ 
\end{thmppal}

\begin{proof}
Nous allons proc\'eder par r\'ecurrence sur la dimension de $X_w$ pour
montrer que Sing $X_w=\Sigma_w.$ On peut supposer $n\ge 4$ puisque les
vari\'et\'es de Schubert de $GL_3/B$ sont lisses. L'\'enonc\'e a \'et\'e \'etabli
dans le cas des vari\'et\'es covexillaires dans \cite{covex}.  

On suppose maintenant l'\'egalit\'e d\'emontr\'ee pour toute vari\'et\'e de
Schubert de dimension $\le k-1,$ et on se donne $X_w$ de
dimension $k.$ On peut supposer $w$ non covexillaire. On fixe alors
une configuration (3412) de $w,$ bien remplie et d'amplitude minimale, et
on consid\`ere les quasi-r\'esolutions de $X_w$ associ\'ees. Etant donn\'e une
composante irr\'eductible $X_v$ de Sing $X_w,$  
ou bien $X_v\subseteq  \bigcap\limits_{i=1}^h \EE_{\pi_i},$ ou bien il existe
un entier $i$ tel que $X_v\not\subseteq \EE_{\pi_i}.$ Dans le premier
cas, on obtient $X_v\subseteq \Sigma_w,$ d'apr\`es la proposition \ref{inter}. 
  
Dans le second cas, consid\'erons l'ouvert $\Reg_{\pi_i}.$ Nous
utiliserons le fait suivant : si $\Omega$ 
est un ouvert d'une vari\'et\'e $X,$ les composantes irr\'eductibles de
Sing $X$ qui rencontrent $\Omega$ sont en bijection avec les
composantes irr\'eductibles de Sing $\Omega,$ par l'application
$Y\longmapsto Y\cap\Omega.$   

Ici, la composante irr\'eductible $X_v$ de Sing $X_w$ rencontre
$\Reg_{\pi_i},$ donc $X_v\cap \Reg_{\pi_i}$ est une composante
irr\'eductible de 
Sing $\Reg_{\pi_i}.$ Comme $\pi_i$ induit un isomorphisme de
$\Reg(\pi_i)$ sur $\Reg_{\pi_i}, \ \pi_i^{-1}(X_v\cap \Reg_{\pi_i})$
est une composante irr\'eductible de Sing $\Reg(\pi_i).$ Il existe
alors, d'apr\`es le fait \'evoqu\'e plus haut, une composante
irr\'eductible $Y$ de Sing $Z_i$ telle que $\pi_i^{-1}(X_v\cap \Reg_{\pi_i})=
Y\cap \Reg(\pi_i).$ Par ailleurs, comme la projection de $Z_i$ sur
$P_I/P_{J_i}$ est une fibration localement triviale de fibre
$X_{w_i},$ les composantes irr\'eductibles de Sing $Z_i$ sont les
$P_I\times^{P_{J_i}}X_{v_i},$ avec $X_{v_i}$ composante irr\'eductible de
Sing $X_{w_i}.$ Ainsi, il existe une composante irr\'eductible $X_{v_i}$ de
Sing $X_{w_i}$ telle que $Y=P_I\times^{P_{J_i}}X_{v_i}.$ Il vient
alors, comme $\pi_i$ est surjective et propre,
$X_v=\pi_i(P_I\times^{P_{J_i}}X_{v_i}).$ Par hypoth\`ese de r\'ecurrence,
on a Sing $X_{w_i}=\Sigma_{w_i},$ donc il existe une configuration
$\mathcal K$ de $w_i$ telle que $v_i=\gamma({\scriptstyle \mathcal K})
w_i.$ Il r\'esulte alors de la proposition \ref{corr-config} que
$X_v=X_{w_Iw_{J_i}v_i} \subseteq \Sigma_w.$ Le th\'eor\`eme est d\'emontr\'e.
\end{proof}

\begin{rem}
Le th\'eor\`eme pr\'ec\'edent, combin\'e avec \cite{BP}, §§3.3 et 4.6, permet de
d\'eterminer, pour chaque composante irr\'eductible $X_v$ de Sing $X_w,$ le
polyn\^ome de Kazhdan-Lusztig $P_{v,w}$ et la multiplicit\'e de $X_w$ en
$e_v,$ not\'ee $m_{v,w}$ (on renvoie par exemple \`a \cite{Mat}, chap.5, §14
pour cette notion). Si $\mathcal{N}_{v,w}\simeq \mathcal
C_{i+1,j+1}$ alors 
$$m_{v,w}=\left(\begin{array}{c}i+j\\i\end{array}\right)\ \mbox{ et }\ 
P_{v,w}= 1+q+\ldots +q^u,$$ o\`u $u=\mbox{Min }(i,j),$ tandis que si
$\mathcal{N}_{v,w}\simeq \mathcal K_{2k+1},$ alors 
$$m_{v,w}=2\ \mbox{ et }\ P_{v,w}= 1+q^k.$$ 
\end{rem}

\begin{ex}
On consid\`ere la permutation $$w=(5,10,7,2,9,8,1,6,3,4)$$ dans
$\mathfrak S_{10}.$ 

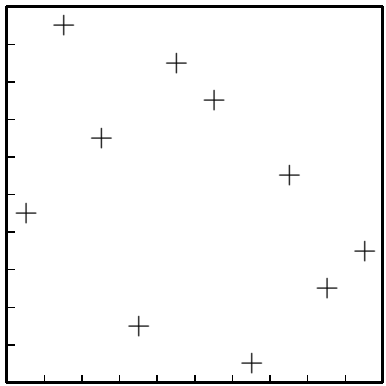
\begin{center}\begin{picture}(5,6)
\thicklines
\multiput(0,0.5)(0,5){2}{\line(1,0){5}}
\multiput(0,0.5)(5,0){2}{\line(0,1){5}}

\linethickness{0.01mm}
\multiput(0.5,0.5)(0.5,0){10}{\line(0,1){0.1}}
\multiput(0,1)(0,0.5){10}{\line(1,0){0.1}}
\thinlines
\put(0.1,2.65){$+$}
\put(0.6,5.15){$+$}
\put(1.1,3.65){$+$}
\put(1.6,1.15){$+$}
\put(2.1,4.65){$+$}
\put(2.6,4.15){$+$}
\put(3.1,0.65){$+$}
\put(3.6,3.15){$+$}
\put(4.1,1.65){$+$}
\put(4.6,2.15){$+$}

\end{picture}\end{center}

Il y a deux points coessentiels bien bord\'es $P_1=(5,5)$ et
$P_2=(5,7)$ ; $P_1$ admet trois bordages minimaux, donn\'es par les
couples de points d'abscisses $(3,7),\ (3,9)$ et $(3,10),$ tous d\'eg\'en\'er\'es
(l'un au Nord-Est, les deux autres au Sud-Ouest). Le point $P_2$ admet
deux bordages minimaux, donn\'es par les couples de points d'abscisses
$(2,7)$ et $(2,8),$ non d\'eg\'en\'er\'es, qui donnent deux
composantes de type $S_1.$ D'autre part, les configurations (3412)
incompressibles sont d'abscisses
$$\begin{array}{l}1,3,4,9\\
1,3,4,10,\\
1,6,7,9,\\
1,6,7,10,\\
1,8,9,10,\\
3,6,7,8,\\
3,6,9,10\ (\ast).\end{array}$$ Elles donnent toutes des configurations
II mixtes, sauf celle qui est marqu\'ee d'une $(\ast),$ qui donne une
configuration II pure. 

Le lieu singulier de $X_w$ a donc 9 composantes 
irr\'eductibles. Elles sont donn\'ees dans le tableau ci-apr\`es. Dans la
colonne de gauche, on a \'ecrit les configurations, donn\'ees par les
ordonn\'ees des points, dans l'ordre des abscisses croissantes, ainsi que le
type de la configuration. Pour les configurations I, les ordonn\'ees 
$\yinf$ et $\Yinf$ sont marqu\'ees en gras ; pour les configurations II, les
ordonn\'ees des points de la configuration (3412) sont indiqu\'ees en
gras. La composante irr\'eductible correspondante $v$ est donn\'ee dans la
deuxi\`eme colonne ; les points o\`u elle diff\`ere de $w$ sont en gras. Dans
les quatre colonnes suivantes, on donne la classe d'isomorphisme de la
transversale, sa dimension $d=\dim\mathcal{N}_{v,w}=\ell(w)-\ell(v),$
le polyn\^ome de Kazhdan-Lusztig $P_{v,w},$ et enfin la multiplicit\'e $m_{v,w}.$  

$$\begin{array}{llccclc}
\multicolumn{2}{c}{{\rm Configuration}} &v& \mathcal{N}_{v,w}& d & P_{v,w}&m_{v,w}\\ 
 & & & & & & \\
\mbox{I} & {\bf 10},7,2,9,8,{\bf 1} &

(5,{\bf 7},{\bf 2},{\bf 1},{\bf 10},{\bf 9},{\bf 8},6,3,4)
& \mathcal C_{3,3} & 5& 1+q+q^2 &6\\

\mbox{I} &{\bf 10},7,9,8,{\bf 6}
& (5,{\bf 7},{\bf 6},2,{\bf 10},{\bf 9},1,{\bf 8},3,4)

& \mathcal C_{2,3}&4 & 1+q&3\\
 & & \\

\mbox{II}_m&
{\bf 5},{\bf 7},{\bf 2},6,{\bf 3} &

({\bf 2},10,{\bf 5},{\bf 3},9,8,1,{\bf 7},{\bf 6},4)
& \mathcal C_{3,2}&4 & 1+q &3  \\

\mbox{II}_m& {\bf 5},{\bf 7},{\bf 2},6,{\bf 4}
&

({\bf 2},10,{\bf 5},{\bf 4},9,8,1,{\bf 7},3,{\bf 6})&
\mathcal C_{3,2}&4 & 1+q &3 \\

\mbox{II}_m& {\bf 5},{\bf 8},{\bf 1},6,{\bf 3}
&

({\bf 1},10,7,2,9,{\bf 5},{\bf 3},{\bf 8},{\bf 6},4)&
\mathcal C_{3,2}&4 & 1+q &3 \\

\mbox{II}_m&{\bf 5},{\bf 8},{\bf 1},6,{\bf 4}
&

({\bf 1},10,7,2,9,{\bf 5},{\bf 4},{\bf 8},3,{\bf 6})&\mathcal
C_{3,2} &4 & 1+q & 3  \\

\mbox{II}_m &{\bf 5},{\bf 6},{\bf 3},{\bf 4}
& 

({\bf 3},10,7,2,9,8,1,{\bf 5},{\bf 4},{\bf 6})&
\mathcal C_{2,2} &3 &1+q & 2\\

\mbox{II}_m&
{\bf 7},2,{\bf 8},{\bf 1},{\bf 6}
&

(5,10,{\bf 2},{\bf 1},9,{\bf 7},{\bf 6},{\bf 8},3,4)&
\mathcal C_{3,2} & 4 & 1+q & 3 \\ 

\mbox{II}_p& {\bf 7},{\bf 8},6,{\bf 3},{\bf 4}
& 

(5,10,{\bf 3},2,9,{\bf 7},1,6,{\bf 4},{\bf 8})&
\mathcal K_5 &5 & 1+q^2 & 2\\
\end{array}$$
\end{ex}

\section{Preuve des propositions \ref{inter} et \ref{corr-config}}

Dans cette derni\`ere section, nous d\'emontrons les propositions
\ref{inter} et \ref{corr-config}. Pour cela, on va d'abord d\'ecrire, pour chaque
entier $i,$ les composantes irr\'eductibles de l'image du lieu
exceptionnel de $\pi_i.$    

\subsection{Un lemme}

En vue de cette description, on
commence par \'etablir un lemme, d\^u \`a P. Polo. 
Soient $z\in \mathfrak S_n$ et $j\in [1, n-1]$ tels que $z<s_j z.$
A chaque point $(b',\b')$ de la fronti\`ere Sud-Est de $\Gamma_z\cap \{(p,q)\
|\ p<z^{-1}(j),\ q>j+1\},$ on associe une permutation $\t_{b'},$ dite de
type Nord-Ouest, d\'efinie comme suit : on note $y_1<\cdots<y_s$ les
ordonn\'ees des points de la fronti\`ere Sud-Ouest de l'ensemble
$\Gamma_z\cap \mathcal R^{(j+1,\b')}(z).$ La permutation $\t_{b'}$
associ\'ee est $$\t_{b'}=(j+1,y_1)(y_1,y_2)\cdots
(y_{s-1},y_s)(y_s,\b')z.$$ 
De fa\c con similaire, \`a chaque point $(c',\g')$ de la fronti\`ere
Nord-Ouest de $\Gamma_z\cap \{(p,q)\ |\ p>z^{-1}(j+1),\ q<j\},$ on
associe une permutation $\t_{c'},$ dite de type Sud-Est : on note
$y_{-t}<\cdots<y_{-1}$ les ordonn\'ees des points de la fronti\`ere
Nord-Est de l'ensemble $\Gamma_z\cap \mathcal R^{(\g',j)}(z),$ et la
permutation $\t_{c'}$ associ\'ee est
$$\t_{c'}=(j,y_{-1})(y_{-1},y_{-2})\cdots(y_{-t},\g')z.$$ 
Enfin, si $\big((b',\b'),(c',\g')\big)$ est un couple de points de
$\Gamma_z$ v\'erifiant $$z^{-1}(j)<b'<c'<z^{-1}(j+1),\ 
\g'<j\mbox{ et } j+1<\b',$$ et tels que le rectangle $\mathcal R_{(b',c')}(z)$
ne contienne pas de point du graphe de $z,$ on lui associe une
permutation $\t_{(b',c')},$ dite de type mixte : on d\'efinit les ordonn\'ees
$y_1<\cdots<y_s$ et $y_{-t}<\cdots<y_{-1}$ comme ci-dessus, et la
permutation $\t_{(b',c')}$ est alors
$$\t_{(b',c')}=\big[(j+1,y_1)(y_1,y_2)\cdots (y_{s-1},y_s)(y_s,\b')
\big]\big[(j,y_{-1})(y_{-1},y_{-2})\cdots(y_{-t},\g') \big] z.$$ 

On a alors le 
\begin{lem}\label{lemme-clef}
Soient $z\in \mathfrak S_n$ et $j\in [1, n-1]$ tels que $z<s_j z.$ Les
\'el\'ements maximaux de $\{\t\in \mathfrak S_n\ |\ \t< z \mbox{ et
  }s_j\t<\t\}$ sont pr\'ecis\'ement les permutations d\'ecrites ci-avant.
\end{lem}

\begin{proof} 
On construit pour commencer un ``ordre de r\'eflexion'' sur
l'ensemble des transpositions de $\mathfrak S_n$ ({\it cf.}
\cite{Dyer}). On consid\`ere la permutation $\sigma_j$ d\'efinie par 
la d\'ecomposition r\'eduite $\sigma_j=s_{j-1}\ldots s_1\,s_{j+1}\ldots
s_{n-1}.$ On consid\`ere une d\'ecomposition r\'eduite $s_{i_1}\cdots s_{i_N}$
de $w_0,$ l'\'el\'ement de plus grande longueur de $\mathfrak S_n$,
obtenue par concat\'enation de celle de $\sigma_j$ ci-dessus 
et d'une d\'ecomposition r\'eduite de $\sigma_j^{-1}w_0.$ D'apr\`es
\cite{Dyer}, prop. 2.13, en posant pour $k=1, \ldots, N,\
t_k=(s_{i_1}\cdots s_{i_{k-1}})\, s_{i_k}\, (s_{i_{k-1}}\cdots
s_{i_1}),$ l'ordre d\'efini par $t_1\prec\cdots\prec t_N$ est un ordre
de r\'eflexion. On a,  
$$\begin{array}{l} \mbox{pour } k=1, \ldots, j-1,\quad t_k=(j-k,j),\\
\mbox{pour } k=j,\ldots, n-2,\quad  t_k=(j+1,k+2).\end{array}$$ 

Cet ordre \'etant construit, on consid\`ere une permutation $\t$ telle
que $\t< z$ et $s_j\t<\t,$ maximale pour cette propri\'et\'e. Soit
$k=\ell(z)-\ell(\t),$ d'apr\`es 
\cite{Dyer}, prop. 4.3, il existe des transpositions $t_{l_1}\prec\cdots
\prec t_{l_k}$ telles que $\t=t_{l_k}\cdots t_{l_1}z,$ avec pour
tout $p=1, \ldots,k,\ \ell(t_{l_p}\cdots t_{l_1}z)=\ell(z)-p.$ Or on
observe que l'on a $t_{l_k}\in \{t_1, \ldots, t_{n-2}\}$ : 
c'est \'evident si $k=1,$ et cela r\'esulte de la maximalit\'e de $\t$
sinon. Alors, comme cet ensemble de r\'eflexions est un id\'eal pour
l'ordre $\prec,$ on en d\'eduit que  $t_{l_p}\in \{t_1, \ldots,
t_{n-2}\}$ pour tout $p=1, \ldots, k.$ Plus pr\'ecis\'ement, il existe des
entiers $q_1<\cdots< q_m$ et $q'_1<\cdots<q'_{m'},$ avec $m,m'\ge 0,$ et
$m+m'=k,$ tels que 
$$\left\{\begin{array}{lll}t_{l_p}&= (j-q_p,j) & \mbox{ pour
      }p=1,\ldots, m,\\   
t_{l_{m+p}}&=(j+1, j+1+q'_p) & \mbox{ pour } p=1,\ldots,
m'.\end{array}\right.$$ 

Supposons d'abord $m'=0.$ Alors le fait que $s_j\t<\t$ entra\^\i ne
$\zinv(j-q_k)>\zinv(j+1).$ Puis, par maximalit\'e de $\t,$ on obtient au
contraire $\zinv(j-q_{k-1})<\zinv(j+1).$ Enfin, comme la longueur
diminue de 1 \`a chaque \'etape \'el\'ementaire du produit $t_{l_k}\cdots
t_{l_1}z,$ on a 
$\zinv(j)<\zinv(j-q_1)<\cdots<\zinv(j-q_k),$ et de plus, posant
$q_0=0,$ on a 
$$\Gamma_z\cap \ \big(\,]\,\zinv(j-q_p),\zinv(j-q_{p+1})\,[\,\times\,
]\,j-q_{p+1}, j\,[\,\big)=\emptyset$$ pour $p=0,\ldots,k-1.$ Il en r\'esulte que le point
$(\zinv(j-q_k), j-q_k)$ est  dans la fronti\`ere Nord-Ouest de
$\Gamma_z\cap \{(p,q)\ |\ p>z^{-1}(j+1), q<j\},$ et que les points 
$(\zinv(j-q_p),j-q_p),$ pour $p=1,\ldots,k -1$ forment la
fronti\`ere Nord-Est de $\Gamma_z\cap \mathcal R^{(j-q_k, j)}(z).$ Ainsi
$\t$ est une permutation de type Nord-Ouest. On proc\`ede de m\^eme
lorsque $m=0$ pour montrer que $\t$ est une permutation de type Sud-Est. 

Supposons maintenant $mm'\not=0.$ Alors, le fait que $s_j\t<\t$
entra\^\i ne $\zinv(j+1+q'_{m'})<\zinv(j-q_m).$ Puis, par maximalit\'e de
$\t,$ on obtient $\zinv(j-q_{m-1})<\zinv(j+1+q'_{m'})$ et
$\zinv(j-q_m)<\zinv(j+1+q'_{m'-1}).$ Et, pour terminer, comme la
longueur diminue de 1  \`a chaque \'etape \'el\'ementaire du produit 
$t_{l_k}\cdots t_{l_1}z,$ on a d'abord
$\zinv(j)<\zinv(j-q_1)<\cdots<\zinv(j-q_m),$ et
$\zinv(j+1+q'_{m'})<\cdots <\zinv(j+1+q'_1)<\zinv(j+1),$ et de plus,
posant $q_0=q'_0=0,$ l'on a 
$$\displaylines{
\Gamma_z\cap\ \big(\, ]\,\zinv(j-q_p),\zinv(j-q_{p+1})\,[\,\times\,
]\,j-q_{p+1}, j\,[\,\big)=\emptyset\cr 
\hfill\mbox{ pour } p=0,\ldots,m-1,\hspace*{2em}\cr
\Gamma_z\cap \ \big(\,
]\,\zinv(j+1+q'_{p+1}),\zinv(j+1+q'_p)\,[\,\times\, ]\,j+1,
j+1+q'_{p+1}\,[\,\big)=\emptyset\cr 
\hfill\mbox{ pour } p=0,\ldots,m'-1.\hspace*{2em}}$$
On en d\'eduit que $\t$ est une permutation de type mixte.  

R\'eciproquement, il est clair que si $\t$ est une permutation de l'un des types
d\'ecrits avant l'\'enonc\'e, elle v\'erifie $s_j\t<\t.$ De plus, si $\t_{b'}$ est
de type Nord-Ouest, associ\'ee au point $(b',\b'),$ les points de $\Gamma_z$
d'ordonn\'ees 
$$j+1<y_1<\cdots<y_s<\b'$$ 
forment une configuration I d\'eg\'en\'er\'ee au Sud-Ouest, et $\t_{b'}$ est
la permutation associ\'ee. La r\'egion de $[1,n]^2$ associ\'ee comme en
\ref{D1} est ici not\'ee $D_{NO}.$ On a alors, d'apr\`es le lemme \ref{long}, 
$$r_{\t_{b'}}=r_z+\chi_{D_{NO}} \ \ (1),$$ 
d'o\`u en particulier $\t\le z.$ De m\^eme, si $\t_{c'}$ est de type
Sud-Est, associ\'ee au point $(c',\g'),$ les points de $\Gamma_z$ d'ordonn\'ees
$$\g'<y_{-t}<\cdots <y_{-1}<j$$ 
forment une configuration I d\'eg\'en\'er\'ee au Nord-Est, et $\t_{c'}$ est la
permutation associ\'ee. La r\'egion de $[1,n]^2$ associ\'ee comme en
\ref{D1} est ici not\'ee $D_{SE}.$ On a 
$$r_{\t_{c'}}=r_z+\chi_{D_{SE}} \ \ (2).$$ 
Enfin, si $\t_{(b',c')}$ est de type mixte, associ\'ee au couple
$\big((b',\b'),(c',\g')\big),$ alors d'une part les points de
$\Gamma_z$ d'ordonn\'ees 
$$j+1<y_1<\cdots<y_s<\b'$$ 
forment une configuration I d\'eg\'en\'er\'ee au Sud-Ouest, et l'on d\'efinit comme
ci-dessus la r\'egion $D_{NO}.$ D'autre part, les points de $\Gamma_z$
d'ordonn\'ees 
$$\g'<y_{-t}<\cdots <y_{-1}<j$$ 
forment une configuration I d\'eg\'en\'er\'ee au Nord-Est, et l'on d\'efinit
comme ci-dessus la r\'egion $D_{SE}.$ On pose alors $D_M=D_{NO}\cup
D_{SE},$ et l'on a  
$$r_{\t_{(b',c')}}=r_z+\chi_{D_M} \ \ (3).$$
Il r\'esulte alors de (1), (2) et (3) que les diff\'erents $\t$
sont deux \`a deux incomparables. Le lemme est d\'emontr\'e.
\end{proof}

\subsection{Le lieu exceptionnel $\Ex (\pi_i)$}

On revient maintenant \`a la situation de la section \ref{quasiresol}, avec
les m\^emes notations. Etant donn\'e $i\in[1,h],$ on note $\partial NO(i)$
la fronti\`ere Sud-Est de  
$$NO(i)=\Gamma_w\cap\{(p,q)\ |\ p<\winv(\a-i+1),\ q>\a'\},$$
et $\partial SE(i)$ la fronti\`ere Nord-Ouest
de $$SE(i)=\Gamma_w\cap \{(p,q)\ |\ p>\winv(\a-i),\ q<\d'\}.$$ 
A chaque point $(b',\b')$ de $\partial NO(i)$ on associe la permutation
$t^{i}(b'),$ \index{$t^{i}(b')$} dite de type Nord-Ouest, d\'ecrite
comme suit. Soit $\tilde{\a}_{b'},$ le plus grand 
entier de l'intervalle  $[\,\a-i+1,\a'\,]$ tel que
$b'<\winv(\tilde{\a}_{b'})$ ; si l'ensemble 
$\Gamma_w\cap\{(p,q)\ |\ \winv(\a-i+1)<p<\winv(\a-i),\ \a'<q<\b'\}$
est non vide, on note $x_s<\cdots<x_1$ les abscisses des points de sa
fronti\`ere Sud-Ouest. Alors les points d'abscisses $$b'<\winv(\tilde{\a}_{b'})
<\cdots<\winv(\a-i+1)<x_s<\cdots<x_1<\winv(\a-i)$$ forment une
configuration I de $w,$ \'eventuellement d\'eg\'en\'er\'ee au Nord-Est, not\'ee
$\mathcal I,$ 
et l'on pose $t^{i}(b')=\tau({\scriptstyle\mathcal I}).$ 

\begin{ex}

On reprend la permutation $w=(8,7,9,6,1,11,5,4,2,10,3)$ de l'exemple
\ref{graphe w_i}, et la configuration (3412) bien remplie et
d'amplitude minimale donn\'ee par les points d'abscisses 2,3,9,11.

\begin{minipage}{6cm}

Fixons $i=4$ : l'ensemble $\partial NO(4)$ est constitu\'e des deux
points du graphe d'abscisse respective 3 et 6. Prenons pour notre
exemple $b'=6.$ Les points de la configuration I qui permet de
d\'efinir  $t^4(6)$ sont repr\'esent\'es par des $\oplus$ sur le
diagramme.
\end{minipage}
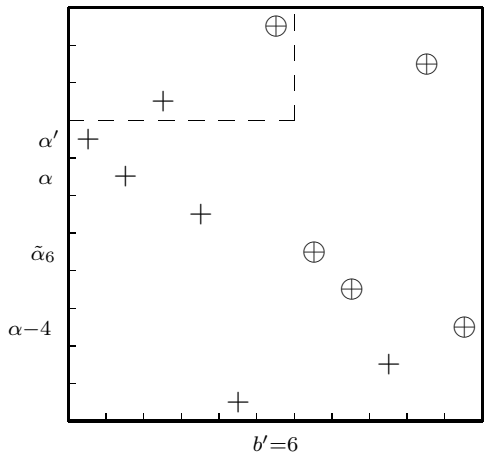
\begin{minipage}{6cm}
\begin{center}
\begin{picture}(5.5,6.5)

\thicklines
\put(1,0){\line(1,0){5.5}}
\put(1,5.5){\line(1,0){5.5}}
\put(1,0){\line(0,1){5.5}}
\put(6.5,0){\line(0,1){5.5}}

\linethickness{0.01mm}
\multiput(1.5,0)(0.5,0){10}{\line(0,1){0.1}}
\multiput(1,0.5)(0,0.5){10}{\line(1,0){0.1}}

\multiput(1,4)(0.4,0){8}{\line(1,0){0.2}}
\multiput(4,4)(0,0.4){4}{\line(0,1){0.2}}

\thinlines  

\put(1.1,3.65){$+$}
\put(1.6,3.15){$+$}
\put(2.1,4.15){$+$}
\put(2.6,2.65){$+$}
\put(3.1,0.15){$+$}
\put(3.6,5.15){$\oplus$}
\put(4.1,2.15){$\oplus$}
\put(4.6,1.65){$\oplus$}
\put(5.1,0.65){$+$}
\put(5.6,4.65){$\oplus$}
\put(6.1,1.15){$\oplus$}

\put(0.6,3.15){${\scriptstyle\a}$}
\put(0.6,3.65){${\scriptstyle\a'}$}
\put(0.5,2.15){${\scriptstyle\tilde{\a}_{6}}$}
\put(0.2,1.15){${\scriptstyle\a-4}$}
\put(3.45,-0.4){${\scriptstyle b'=6}$}

\end{picture}\end{center}
\end{minipage} 
\end{ex}

\vspace*{1cm}

De mani\`ere analogue, \`a chaque point $(c', \g')$ de $\partial SE(i),$ on
associe la permutation $t_i(c'),$ \index{$t_i(c')$} dite de type
Sud-Est, d\'ecrite comme suit. Soit 
$\tilde{\d}_{c'},$ le plus petit entier de $[\,\d',\a-i\,]$ tel que
$\winv(\tilde{\d}_{c'})<c'$ ; si l'ensemble $\Gamma_w\cap\{(p,q)\ |\
\winv(\a-i+1)<p<\winv(\a-i),\ \g'<q<\d'\}$ est non vide, on note
$x_{-1}<\cdots<x_{-t}$ les points de sa fronti\`ere Nord-Ouest. Alors
les points d'abscisses
$$\winv(\a-i+1)<x_{-1}<\cdots<x_{-t}<\winv(\a-i)
<\cdots<\winv(\tilde{\d}_{c'})<c'$$        
forment une configuration I de $w,$ \'eventuellement d\'eg\'en\'er\'ee au
Sud-Ouest, not\'ee 
$\mathcal I',$ et l'on pose $t_{i}(c')=\tau({\scriptstyle\mathcal I'}).$

Enfin, si $\winv(\a-i+1)<b'<c'<\winv(\a-i)$ est une
configuration (3412) incompressible de $w,$ not\'ee $\II,$ on note
$m_i(b',c')$ \index{$m_i(b',c')$} la permutation
$\sigma({\scriptstyle\II}),$ dite de type mixte. 

\begin{prop}\label{compex}
Soit $i\in [1,h].$ Les composantes irr\'eductibles de $\EE_{\pi_i},$
l'image du lieu exceptionnel $\Ex\,(\pi_i),$ correspondent exactement
aux permutations 
$t^i(b'),$ $t_i(c')$ et $m_i(b',c')$ d\'ecrites ci-dessus. 
\end{prop}

\begin{proof}
Puisque $\pi_i$ est propre, l'ensemble des $x\in X_w$ dont la fibre
$\pi_i^{-1}(\{x\})$ est finie est un ouvert de $X_w$ (voir, par exemple,
\cite{K}, Prop. 6.4.5). Puis, comme $\pi_i$ est birationnelle et que
$X_w$ est normale (voir, par exemple, \cite{Ram}), on d\'eduit du th\'eor\`eme
principal de Zariski que l'image du lieu exceptionnel de $\pi_i,\
\EE_{\pi_i},$ est \'egale \`a $$\{x\in X_w\, |\, \# 
\pi_i^{-1}(\{x\})>1 \}.$$ 

Le morphisme $\pi_i$ est $P_I$-\'equivariant, donc les composantes
irr\'eductibles de $\EE_{\pi_i}$ sont des vari\'et\'es de Schubert $X_v$ avec
$v\le w,\ v\in\, ^I\mathfrak S_{max}.$ On se donne une telle
composante $X_v,$ et l'on note $\mathcal V_1,\ldots, \mathcal V_m$ les
composantes irr\'eductibles de $\pi_i^{-1}(X_v).$ Elles sont toutes de la forme
$P_I\times^{P_{J_i}}X_{t_j},$ avec 
$t_j<w_i,\ t_j\in\, ^{J_i}\mathfrak S_{max}.$ Comme $\pi_i$ est
surjective et propre, on a  $$X_v=\bigcup\limits_{j=1}^m
\pi_i(\mathcal V_j)$$ et chaque $\pi_i(\mathcal V_j)$ est ferm\'e. Or
$X_v$ est irr\'eductible, il existe donc $j$ tel que $X_v=\pi_i(\mathcal
V_j)$ ; un tel $j$ est unique. En
effet, $t_j$ est une permutation maximale telle que $t_j\le w_i,
t_j\in \, ^{J_i}\mathfrak S_{max}$ et $t_j\in \mathfrak S_Iv.$ Comme
$v,w\in\, ^I\mathfrak S_{max}$ et $v\le w,$ d'apr\`es un lemme de  
Deodhar (voir \cite{lemme-d}, Lemma 11.1), il existe un unique $t\in
\mathfrak S_Iv$ tel que $t\le w_i,$ maximal pour cette
propri\'et\'e. Il vient alors $t\in \, ^{J_i}\mathfrak S_{max}.$ Ainsi
$t_j=t$ est unique ; on peut supposer $j=1.$ Soit
$\t_1=w_{J_i}t_1,$ le repr\'esentant minimal de $\mathfrak S_{J_i}t_1.$
On a ainsi 
$$\dim(P_I\times^{P_{J_i}}X_{t_1})= \dim (P_I/P_{J_i}) + \dim (X_{t_1})=
\ell(w_I)+\ell(\t_1),$$ et $$X_v =
\pi_i(P_I\times^{P_{J_i}}X_{t_1})= X_{w_I*\t_1}$$ 
(l'\'egalit\'e $\pi_i(P_I\times^{P_{J_i}}X_{t_1})= X_{w_I*\t_1}$ se
d\'emontre par le m\^eme argument que l'\'egalit\'e $\pi(Z_i)=X_w$ de la
proposition \ref{prop-pi}).
On en d\'eduit que $\t_1\not\in\,^I\mathfrak S_{min}.$ En effet,
supposons au contraire $\t_1\in\,
^I\mathfrak S_{min},$ il vient
$\dim(X_v)=\dim(P_I\times^{P_{J_i}}X_{t_1}).$ Alors si $z\in
\pi_i^{-1}(e_v),$ l'orbite $U\,z$ est un
ouvert dense de $P_I\times^{P_{J_i}}X_{t_1}.$ Il en r\'esulte que la
fibre $\pi_i^{-1}(e_v)$ est un singleton, et l'on a alors
$X_v\not\subseteq \EE_{\pi_i},$ une contradiction. On a donc $\t_1\not\in\,
^I\mathfrak S_{min},$ ce qui \'equivaut, comme $\t_1\in\, ^{J_i}\mathfrak
S_{min},$ \`a $s_{k_i}\t_1<\t_1.$ Finalement, on a montr\'e que $v=w_I*\t,$
avec $\t<w_Iw,\ \t\in\, ^{J_i}\mathfrak S_{min},$ et $s_{k_i}\t<\t.$ 

R\'eciproquement, si l'on se donne un tel $\t,$ on a $w_{J_i}\t<w_i,\ \pi_i( 
P_I\times^{P_{J_i}}X_{w_{J_i}\t})=X_{w_I*\t},$ et $\dim
(P_I\times^{P_{J_i}}X_{w_{J_i}\t})=\ell(w_I)+\ell(\t)>\ell(w_I*\t)$
puisque $\t\not\in\, ^I\mathfrak S_{min}.$ Il en r\'esulte que $X_{w_I*\t}$ est
contenu dans l'image du lieu exceptionnel de $\pi_i.$ 

Ainsi, $\EE_{\pi_i}$ est la r\'eunion des $X_{w_I*\t}$
pour $\t$ comme ci-dessus. Comme de plus on a $w_I*\t\le w_I*\t'$ si
$\t\le \t',$ on peut se restreindre aux tels $\t$ maximaux. 

En appliquant le lemme \ref{lemme-clef} \`a $z=w_Iw$ et $j=k_i,$ on
obtient la description de l'ensemble $\mbox{Max}\,_{k_i}(w_Iw)$ des
\'el\'ements maximaux de 
$\{u\in \mathfrak S_n\ |\ u<w_Iw \mbox{ et }s_{k_i}u<u\}.$ Comme
$w_Iw\in\, ^I\mathfrak S_{min},$ on a
$$(w_Iw)^{-1}(\d')<\cdots< (w_Iw)^{-1}(\a').$$ 
D'apr\`es le lemme \ref{lemme-clef}, tout \'el\'ement $u$ de
$\mbox{Max}\,_{k_i}(w_Iw)$ v\'erifie les conditions suivantes :
$$\left\{\begin{array}{lcl}u^{-1}(j)&=&(w_Iw)^{-1}(j) \mbox{ pour tout }
  j\in [\d',\a'],\, j\not=k_i,\,k_i+1,\\
u^{-1}(k_i)&\ge& (w_Iw)^{-1}(k_i),\\
u^{-1}(k_i+1)&\le& (w_Iw)^{-1}(k_i+1).\end{array}\right.$$
On en d\'eduit en particulier que $u^{-1}(\d')<\cdots<u^{-1}(k_i),$ et
$ u^{-1}(k_i+1)<\cdots<u^{-1}(\a'),$ c'est-\`a-dire que $u\in\, 
^{J_i}\mathfrak S_{min}.$ L'ensemble $\mbox{Max}\,_{k_i}(w_Iw)$
co\"\i ncide donc avec l'ensemble des $\t$ maximaux cherch\'es.

Il reste maintenant \`a d\'ecrire explicitement les $w_I*\t.$ Remarquons que ce
sont les repr\'esentants maximaux des classes $\mathfrak S_I\t,$ et
qu'ils sont ais\'ement d\'ecrits \`a partir 
des $\t$ : $(w_I*\t)^{-1}$ co\"\i ncide avec $\t^{-1}$ en dehors de
$[\d',\a'],$ et $(w_I*\t)^{-1}$ est d\'ecroissante sur $[\d',\a'].$

Voyons d'abord le cas des permutations $\t$ de type Nord-Ouest. On
rappelle que $k_i=\d'+\a'-\a+i-1,$ et l'on remarque que  
$$(w_Iw)^{-1}(k_i)=\winv(\a-i+1)\mbox{ et }(w_Iw)^{-1}(k_i+1)=\winv(\a-i).$$ 
De plus, on a $(w_Iw)^{-1}([k_i+1,\a'])\subseteq
[(w_Iw)^{-1}(k_i+1),n],$ donc il n'y a pas de point du graphe de
$w_Iw$ dans le rectangle $[1,(w_Iw)^{-1}(k_i)]\times[k_i+1,\a'].$ On a ainsi 
$$\displaylines{\hspace*{2em}\Gamma_{w_Iw}\cap \{(p,q)\ |\
  p<(w_Iw)^{-1}(k_i), q>k_i+1\} 
=\hfill\cr\hfill\Gamma_{w_Iw}\cap \{(p,q)\ |\ p<(w_Iw)^{-1}(k_i), q>
\a'\}.\hspace*{2em}}$$ Comme 
les graphes de $w$ et $w_Iw$ co\"\i ncident pour les ordonn\'ees de
$]\,\a',n],$ on obtient finalement $$\Gamma_{w_Iw}\cap \{(p,q)\ |\
p<(w_Iw)^{-1}(k_i), q> \a'\} =NO(i),$$ 
$NO(i)$ \'etant d\'efini avant l'\'enonc\'e de la proposition.  

Soit alors $(b',\b')$ un point de $\partial NO(i)$ ;
notons $(x_i, y_i),$ pour $i=1,\ldots,s,$ les coordonn\'ees des points de
la fronti\`ere Sud-Ouest de $\Gamma_{w_Iw}\cap \mathcal
R^{(k_i+1,\b')}(w_Iw),$ avec $x_s<\cdots<x_1.$ Comme $(w_Iw)^{-1}$ est 
croissante sur l'intervalle $[\d',\a'],$ il n'y a pas de point du
graphe de $w_Iw$ dans $]\,b',(w_Iw)^{-1}(k_i+1)\,[\,\times\,]\,k_i+1,\a'\,[,$
{\it i.e.} dans $]\,b',\winv(\a-i)\,[\,\times\,]\,k_i+1,\a'\,[.$ 
De plus, comme $(b',\b')$ est un point de la fronti\`ere Sud-Est de $NO(i),$ 
il n'y a pas non plus de point de $\Gamma_{w_Iw}$ dans
$]\,b',\winv(\a-i+1)\,[\,\times \,]\,\a',\b'\,[.$ Il vient donc 
$$\Gamma_{w_Iw}\cap \mathcal R^{(k_i+1,\b')}(w_Iw)\subseteq\ 
]\,\winv(\a-i+1),\winv(\a-i)\,[\,\times\, ]\,\a',\b'\,[,$$ cet
ensemble est donc form\'e de points du graphe de $w.$ 
Soit $\tilde{\a}_{b'}$ le plus grand entier de l'intervalle
$[\a-i+1,\a']$ tel que $b'<\winv(\tilde{\a}_{b'}).$ Alors on v\'erifie que
les points d'abscisses
$b'<\winv(\tilde{\a}_{b'})<\cdots<\winv(\a-i+1)<x_s<\cdots<x_1<\winv(\a-i)$
forment une configuration I de $w,$ not\'ee $\mathcal I,$ et que
$w_I*\t=\tau(\scriptstyle{\mathcal I}).$

Le cas des permutations de type Sud-Est se traite de mani\`ere
analogue. 

\vspace*{1em}

Soit maintenant $\t$ une permutation de type mixte, associ\'ee
\`a un couple $\big((b',\b'),(c',\g')\big)$ de points du graphe
de $w_Iw$ 
v\'erifiant $$(w_Iw)^{-1}(k_i)<b'<c'<(w_Iw)^{-1}(k_i+1),$$ c'est-\`a-dire
$$\winv(\a-i+1)<b'<c'<\winv(\a-i),$$ et $$\g'<k_i\ \mbox{ et }\ k_i+1<\b'.$$ On
suppose aussi que le rectangle $\mathcal R_{(b',c')}(w_Iw)$ ne contient
aucun point du graphe de $w_Iw.$ Les in\'egalit\'es
$b'<\winv(\a-i)$ et $k_i+1<\b'$ donnent $\b'>\a',$
donc $(b',\b')$ est un point du graphe de $w.$ On obtient de
m\^eme $\g'<\d',$ d'o\`u $(c',\g')\in \Gamma_w.$ Ainsi, les points
$\winv(\a-i+1)<b'<c'<\winv(\a-i)$ forment une configuration (3412)
de $w,$ not\'ee $\II.$ Cette configuration est
incompressible : en effet, elle est de hauteur 1, il suffit donc de
v\'erifier que les rectangles $MN$ et $MS$ associ\'es ne contiennent pas
de point du graphe de $w,$ ce qui r\'esulte imm\'ediatement du fait que  
$\mathcal R_{(b',c')}(w_Iw)\cap \Gamma_{w_Iw}=\emptyset.$ Ensuite, on
remarque que les points du graphe de $w_Iw$ contenus dans le rectangle
$\mathcal R^{(\g',k_i)}(w_Iw)$ sont d'ordonn\'ee $<\d',$ donc sont des
points du graphe de $w$ et que leur abscisse est en fait $<b'$ : ce sont donc
exactement les points de $SO_{II}\cap \Gamma_w.$ On voit de m\^eme que les points
du graphe de $w_Iw$ contenus dans le rectangle $\mathcal R^{(k_i+1,\b')}(w_Iw)$
sont des points du graphe de $w$, et que ce sont exactement les points
de $NE_{II}\cap \Gamma_w.$ On obtient alors que  
$w_I*\t=\sigma(\scriptstyle{\II}).$   

On remarque pour terminer que les permutations que l'on vient de d\'ecrire
sont deux \`a deux incomparables, \`a l'aide des propri\'et\'es de leurs
fonctions rang vues aux lemmes \ref{long} et \ref{longII}. 
Elles d\'ecrivent donc les composantes
irr\'eductibles de l'image du lieu exceptionnel de $\pi_i.$
\end{proof}

\begin{ex}\label{mixte}
On consid\`ere la permutation 
$$w=(6,4,2,7,1,5,3)$$ 
dans $\mathfrak S_7,$ dont le graphe est repr\'esent\'e sur le diagramme suivant.

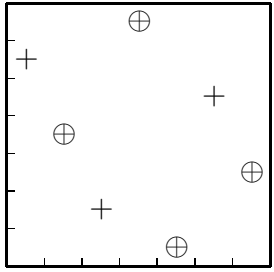
\begin{minipage}{4.5cm}

\begin{picture}(4,4.5)
\thicklines
\put(0,0.5){\line(1,0){3.5}}
\put(0,4){\line(1,0){3.5}}
\put(0,0.5){\line(0,1){3.5}}
\put(3.5,0.5){\line(0,1){3.5}}

\linethickness{0.01mm}
\multiput(0.5,0.5)(0.5,0){6}{\line(0,1){0.1}}
\multiput(0,1)(0,0.5){6}{\line(1,0){0.1}}

\thinlines

\put(0.6,2.15){$\oplus$}
\put(1.6,3.65){$\oplus$}
\put(2.1,0.65){$\oplus$}
\put(3.1,1.65){$\oplus$}

\put(0.1,3.15){+}
\put(1.1,1.15){+}
\put(2.6,2.65){+}

\end{picture}
\end{minipage}
\begin{minipage}{8.5cm}
La configuration (3412) donn\'ee par les points d'abscisses 2, 4, 5 et 7
(repr\'esent\'es par des $\oplus$) est bien remplie et d'amplitude minimale
; elle est incompressible. On a $h=1,\ \a'=\a=4,$ et $\d'=\d=3.$ La permutation
$w_1$ associ\'ee est $$w_1=(6,3,2,7,1,5,4).$$
\end{minipage}

Sur le diagramme suivant, on a repr\'esent\'e dans $\Gamma_w$ le quadrant
$NO(1)$ et la bande verticale d'abscisses $\in\, ]\,a,d\,[.$ 

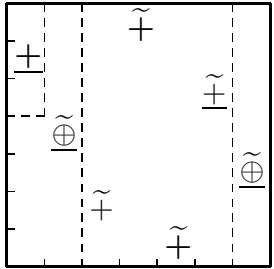
\begin{minipage}{4.5cm}

\begin{picture}(4,4.5)
\thicklines
\put(0,0.5){\line(1,0){3.5}}
\put(0,4){\line(1,0){3.5}}
\put(0,0.5){\line(0,1){3.5}}
\put(3.5,0.5){\line(0,1){3.5}}

\linethickness{0.01mm}
\multiput(0.5,0.5)(0.5,0){6}{\line(0,1){0.1}}
\multiput(0,1)(0,0.5){6}{\line(1,0){0.1}}
\linethickness{0.05mm}
\multiput(0.2,2.5)(0.2,0){2}{\line(1,0){0.1}}
\multiput(0.5,2.5)(0,0.2){8}{\line(0,1){0.1}}
\multiput(1,0.7)(0,0.2){17}{\line(0,1){0.1}}
\multiput(3,0.7)(0,0.2){17}{\line(0,1){0.1}}

\thinlines

\put(0.6,2.15){$\widetilde{\underline{\oplus}}$}
\boldmath
\put(1.6,3.55){$\widetilde{+}$}
\put(2.1,0.65){$\widetilde{+}$}
\unboldmath
\put(3.1,1.65){$\widetilde{\underline{\oplus}}$}

\boldmath\put(0.1,3.2){$\underline{+}$}\unboldmath
\put(1.1,1.15){$\widetilde{+}$}
\put(2.6,2.7){$\widetilde{\underline{+}}$}

\end{picture}
\end{minipage}
\begin{minipage}{8.5cm}
D'apr\`es la proposition qui pr\'ec\`ede, l'image du lieu exceptionnel de $\pi_1$ a
deux composantes irr\'eductibles, l'une de type Nord-Ouest associ\'ee au
point d'abscisse 1 (repr\'esent\'e par
\boldmath$\underline{+}$\unboldmath), l'autre de type mixte, associ\'ee
aux points d'abscisses 4 et 5 (repr\'esent\'es par \boldmath
$\widetilde{+}$\unboldmath). 
\end{minipage}

\noindent La configuration I qui d\'efinit la permutation $t^1(1)$
est non d\'eg\'en\'er\'ee (c'est un fait g\'en\'eral, qui sera \'etabli au lemme
\ref{1-r}), constitu\'ee des points d'abscisses 1, 2, 6, 7 (ce sont 
les points soulign\'es de la figure), et on a $$t^1(1)=(4,3,2,7,1,6,5).$$
La configuration II qui d\'efinit $m_1(4,5)$ est mixte, donn\'ee par les
abscisses 2, 3, 4, 5, 6, 7 (ce sont les points surmont\'es d'un $\
\widetilde{}\ $ sur la figure), et on a $$m_1(4,5)=(6,2,1,4,3,7,5).$$ 
\end{ex}


\subsection{Lieux exceptionnels et lieu singulier}
\subsubsection{Intersection des images des lieux exceptionnels}

Nous pouvons maintenant d\'emontrer la 

\begin{propA}
L'intersection des images des lieux  exceptionnels des $\pi_i,$ pour $i$
parcourant l'intervalle $[1,h],$ est contenue dans $\Sigma_w.$
\end{propA}

\noindent{\it Preuve.} On va montrer que pour toute famille $(X_{v_i})_{i\in
[1,h]},$ o\`u chaque $X_{v_i}$ est une composante irr\'eductible de
$\EE_{\pi_i},$ on a 
$\bigcap\limits_{i=1}^h X_{v_i} \subseteq \Sigma_w.$ Dans cette
perspective, il est utile de faire les remarques suivantes : d'abord, on peut
supposer que tous les $v_i$ sont de type Nord-Ouest ou Sud-Est, car
les composantes de type mixte sont de la forme
$\sigma(\scriptstyle{\II})$ pour une certaine configuration $\II$ de
type II mixte, donc telle que $X_{\sigma(\scriptscriptstyle{\II})} \subseteq
\Sigma_w.$ On peut aussi supposer que les configurations I qui d\'efinissent les
permutations $v_i$ sont d\'eg\'en\'er\'ees, puisque sinon $X_{v_i}\subseteq
\Sigma_w.$ Cela implique en particulier que la composante $X_{v_1}$
est de type Sud-Est, et que la composante $X_{v_h}$ est de type
Nord-Ouest, en vertu du lemme suivant.

\begin{lem}\label{1-r}
$($a$\,)$ Les composantes de type Nord-Ouest (resp. Sud-Est) de
$\EE_{\pi_1}$ (resp. $\EE_{\pi_h}$) sont associ\'ees \`a des configurations I
non d\'eg\'en\'er\'ees.\\
$($b$\,)$ Si $i>1,\ (b,\b)$ est le point le plus \`a l'Ouest de
$\partial NO(i),$ et de m\^eme, si $i<h,\ (c,\g)$ est le point le plus \`a
l'Est de $\partial SE(i).$ 
\end{lem}

\renewcommand\proofname{Preuve du lemme \rm{\ref{1-r}}}

\begin{proof}
Montrons par exemple que les composantes de type Nord-Ouest de
$\EE_{\pi_1}$ sont associ\'ees \`a des configurations I non d\'eg\'en\'er\'ees. On
remarque pour commencer que le fait que la configuration $a<b<c<d$
soit d'amplitude minimale parmi les configurations bien remplies
entra\^\i ne : 
$$w\big(\,]\,\winv(\a'),c\,[\,\big)\,\cap\, ]\,\a',\b\,[\, =\emptyset\
\ (\dagger).$$  
Consid\'erons alors un point $(b',\b')$ de $\partial NO(1)$ ; on a
$b'<a,$ et donc a fortiori $b'<b.$ On constate que si $\b'>\b,$ alors
la configuration I qui 
d\'efinit la composante irr\'eductible $X_{t^1(b')}$ de $\EE_{\pi_1}$
est non d\'eg\'en\'er\'ee, car sa suite NE contient le point $(b,\b).$

Voyons maintenant que si $h>1,$ alors on a $\b'>\b.$ En effet, supposons
$\b'<\b$ ; alors on obtient d'apr\`es $(\dagger),\ b'<\winv(\a').$ Il
r\'esulte aussi de $(\dagger)$ que $\winv(\a'+1)>c,$ et que pour tout
$j\in [\,\a'+1,\b'\,[,$ on a $\winv(j)<\winv(\a')$ ou $\winv(j)>c.$ Il
existe donc un entier $j\in [\,\a'+1,\b'\,[,$ tel que $\winv(j)>c$ et
$\winv(j+1)<\winv(\a').$ Comme de plus $h>1,$ on a $\winv(\a-1)\in
\,]\,b,c\,[,$ et donc les abscisses
$\winv(j+1)<b<\winv(\a-1)<\winv(j)$ forment une configuration (3412) de
$w,$ bien remplie (car de hauteur 1) et d'amplitude $\b-\a+1<\b-\g,$
exclu. On a donc $\b'>\b.$ 

Reste le cas $h=1$ et $\b'<\b$ : comme pr\'ec\'edemment, il existe $j\in
[\,\a'+1,\b'\,[,$ tel que $\winv(j)>c$ et $\winv(j+1)<\winv(\a').$ On a de
plus $\winv(j)\in\,]\,c,d\,[,$ car sinon, les abscisses
$\winv(j+1)<b<d<\winv(j)$ formeraient une configuration (3412) de
$w,$ bien remplie (car de hauteur 1) et d'amplitude $\b-\d<\b-\g,$
exclu. La configuration qui d\'efinit $t^1(b')$ est alors non d\'eg\'en\'er\'ee
car sa suite NE contient le point $(\winv(j),j).$ 
 
On montre de m\^eme que les composantes de type Sud-Est de $\EE_{\pi_h}$
sont associ\'ees \`a des configurations I non d\'eg\'en\'er\'ees. 

Montrons maintenant la deuxi\`eme assertion du lemme. Soit $i>1,$ on
rappelle que $$NO(i)=\Gamma_w\cap\{(p,q)\ |\ p<\winv(\a-i+1),\
q>\a'\}.$$ On a $(b,\b)\in NO(i)$ ; de plus, comme la configuration
$a<b<c<d$ est d'amplitude minimale parmi les configurations bien
remplies, $(b,\b)$ est dans la fronti\`ere Sud-Est, $\partial
NO(i).$ S'il y avait un point plus \`a l'Ouest dans $\partial NO(i),$
d'ordonn\'ee $\b',$ on aurait $\b'<\b,$ et on obtiendrait comme pr\'ec\'edemment
une configuration bien remplie d'amplitude plus petite, une contradiction. 

On montre de m\^eme que pour $i<h,\ (c,\g)$ est le point le plus \`a l'Est
de $\partial SE(i).$ 
\end{proof}

On a remarqu\'e au d\'ebut de la preuve de la proposition \ref{inter} que
les composantes de type mixte sont contenues dans $\Sigma_w.$ Donc,
dans le cas $h=1,$ le lemme \ref{1-r} ach\`eve la preuve de cette 
proposition.  

\begin{ex}
On renvoie \`a l'exemple \ref{mixte}, o\`u l'on a $h=1,$ et o\`u les
composantes irr\'eductibles de $\EE_{\pi_1}$ 
sont associ\'ees \`a la configuration I non d\'eg\'en\'er\'ee donn\'ee par les
abscisses 1, 2, 6, 7, et \`a la configuration II mixte donn\'ee par les
abscisses 2, 3, 4, 5, 6, 7. 
\end{ex}

On suppose dor\'enavant $h\ge 2,$ et on introduit la notion suivante.

\begin{defn}\label{good} Soit $(X_{v_i})_{i\in [1,h]}$ une famille de
  composantes irr\'eductibles des $\EE_{\pi_i},$ chacune de type
  Nord-Ouest ou Sud-Est, 
  et associ\'ee \`a une configuration I d\'eg\'en\'er\'ee. On dit que
  $(X_{v_i})_{i\in [1,h]}$ est une {\it bonne famille}, s'il
  existe des entiers $i<j$ tels que : \\
  \hspace*{2ex}($\,a\,$) la composante $v_i$ soit de type Sud-Est,
  associ\'ee au point $(c_i,\g_i),$\\
  \hspace*{2ex}($\,b\,$) la composante $v_j$ soit de type Nord-Ouest,
  associ\'ee au point $(b_j,\b_j),$\\
  \hspace*{2ex}($\,c\,$) on ait les in\'egalit\'es  
  $$\winv(\a-i+1)<b_j<c_i<\winv(\a-j).$$  
\end{defn}

Nous allons d\'emontrer le

\begin{lem}\label{bonne-inter} 
Si la famille $(X_{v_i})_{i\in [1,h]}$ est bonne, alors 
$\bigcap\limits_{i=1}^h X_{v_i}\subseteq \Sigma_w.$ \end{lem}

\renewcommand\proofname{Preuve du lemme \rm{\ref{bonne-inter}}}

\begin{proof}
Soient $i<j$ v\'erifiant les conditions $(a\,),\  (b\,)$ et $(c\,)$ de la
d\'efinition. On d\'efinit les entiers $\a_j$ et $\d_i$ dans
$[\,\d',\a'\,]$ par les conditions $\winv(\a_j)<b_j<\winv(\a_j-1)$ et 
$\winv(\d_i+1)<c_i<\winv(\d_i).$ 

On a $b\le b_j<\winv(\a-j+1)$ (la  premi\`ere in\'egalit\'e r\'esulte de
l'assertion $(b\,)$ du lemme \ref{1-r}, car $j\ge 2,$ 
et la seconde de la d\'efinition de $NO(j)$). De m\^eme, on a
$\winv(\a-i)<c_i\le c.$ On en d\'eduit que $\a_j\in
[\,\d+2,\a\,],$ et $\d_i\in [\,\d,\a-2\,].$ 

D'autre part, comme $b_j<c_i,$ on a $\winv(\a_j)\le \winv(\d_i+1),$ et donc
$\a_j\ge\d_i+1$ car $\winv$ est d\'ecroissante sur $[\,\d',\a'\,].$
Ainsi les points d'abscisses
$\winv(\a_j)<b_j<c_i<\winv(\d_i)$ forment une configuration (3412), que
l'on note $\II.$ Montrons que $\II$ est incompressible ; il s'agit
d'une configuration (3412) bien 
remplie, il suffit donc de montrer que
$MS = \,]\,b_j,c_i\,[\,\times\,]\,\g_i,\d_i\,[$ et
$MN = \,]\,b_j,c_i\,[\,\times\,]\,\a_j,\b_j\,[$ ne rencontrent pas
$\Gamma_w.$ D'une part, $(c_i,\g_i)\in \partial SE(i),$ on a donc 
$$\big(\,]\,\winv(\a-i),c_i\,[\,\times\,]\,\g_i,\d'\,[\,\big)\, \cap\, \Gamma_w
=\emptyset\hspace{3em} (1).$$ D'autre part, la configuration I
d\'efinissant $t_i(c_i)$ (d\'ecrite avant la proposition \ref{compex}) 
est suppos\'ee d\'eg\'en\'er\'ee, donc on a aussi 
$$\big(\,]\,\winv(\a-i+1),\winv(\a-i)\,[\,\times\,]\,\g_i,\d'\,[\,\big)
\,\cap\, \Gamma_w =\emptyset\hspace{3em}(2).$$ Par ailleurs, comme
$\winv$ est d\'ecroissante sur $[\,\d',\a'\,],$ on a
$\winv\big(\,[\,\d',\d_i\,]\,\big) \subseteq [\,\winv(\d_i),n\,],$ en
particulier $$\winv\big(\,[\,\d',\d_i\,]\,\big)\, \cap
\,]\,\winv(\a-i+1),c_i\,[\,=\emptyset\hspace{3em}(3).$$ On
d\'eduit de $(1), (2),$ et $(3)$ que 
$$\big(\,]\,\winv(\a-i+1),c_i\,[\,\times\,]\,\g_i,\d_i\,[\,\big)\,\cap\,
\Gamma_w =\emptyset\hspace{3em}(\diamond).$$ En particulier, il en r\'esulte que
$MS$ ne rencontre pas $\Gamma_w,$ puisque $b_j> \winv(\a-i+1).$ 

On d\'emontre de m\^eme
$$\,]b_j,\winv(\a-j)\,[\,\times\,]\,\a_j,\b_j\,[\,\cap\, \Gamma_w
=\emptyset\ \ (\diamond\diamond),$$ et l'on en d\'eduit que 
$MN$ ne rencontre pas $\Gamma_w.$ La configuration $\II$ est donc
incompressible. 

Notant simplement $\sigma$ la permutation associ\'ee \`a $\II,$ nous
allons montrer que $X_{v_i}\cap X_{v_j}\subseteq 
X_{\sigma}.$ 

La configuration I qui d\'efinit $v_j,$ d\'eg\'en\'er\'ee au
  Nord-Est, est donn\'ee par les points d'abscisses
  $b_j<\winv(\a_j-1)<\cdots<\winv(\a-j+1)<\winv(\a-j).$ Celle qui
  d\'efinit $v_i$ est d\'eg\'en\'er\'ee au Sud-Ouest, et donn\'ee par les points
  d'abscisses 
  $\winv(\a-i+1)<\winv(\a-i)<\cdots<\winv(\d_i+1)<c_i.$ 

  Il r\'esulte de $(\diamond)$ et $(\diamond\diamond)$ que les
  rectangles $SO_{II}$ et $NE_{II}$ de la configuration $\II$ ne
  rencontrent pas $\Gamma_w,$ 
  donc le graphe de $\sigma$ ne diff\`ere de celui de $w$ qu'en les
  points d'abscisses $\winv(\a_j),\ b_j,\ c_i$ et $\winv(\d_i).$  

On en d\'eduit que les graphes des permutations $w,\ v_j,\ v_i,$ et $\sigma$
ne diff\`erent qu'en les points d'abscisses
$\winv(\a-i+1)<\cdots<\winv(\a_j)<b_j<\winv(\a_j-1)<\cdots<\winv(\d_i+1)<c_i<
\winv(\d_i)<\cdots<\winv(\a-j).$      
Notons $\overline{v_j},\ \overline{v_i}$ et $\overline{\sigma}$ les
permutations obtenues en focalisant sur ces abscisses ; soient 
$$\begin{array}{lcl}m_1 & = & \#\,[\,\a_j,\a-i+1\,] \ge 1,\\
m_2 & = & \#\,[\,\d_i+1, \a_j-1\,] \ge 0,\\
m_3 & = & \#\,[\,\a-j, \d_i\,] \ge 1.\end{array}$$ et
soit $m=m_1+m_2+m_3+2.$ Les permutations $\overline{v_j},\
\overline{v_i}$ et $\overline{\sigma}$ 
sont dans $\mathfrak S_m$ ; $\overline{v_j}$ est la permutation
maximale telle que $\overline{v_j}(m_1+m_2+2)=1$ et 
$\overline{v_j}(m)=m,\ \overline{v_i}$ est la permutation maximale
telle que $\overline{v_i}(1) = 1$ et 
$\overline{v_i}(m_1)=m,$ et $\overline{\sigma}$ est la permutation
maximale de $\mathfrak S_m$ telle que $\overline{\sigma}(m_1)=1$ et
$\overline{\sigma}(m_1+m_2+3)=m.$ Soit $v$ la permutation maximale de
$\mathfrak S_m$ telle que $v(1)=1$ et $v(m)=m$ ; on voit facilement
que $v$ est le plus grand
\'el\'ement de $\Lambda(\overline{v_i},\overline{v_j}),$ et que $v\le 
\overline{\sigma}.$ On en d\'eduit, en vertu du lemme
\ref{focalisation}, que $X_{v_i}\cap X_{v_j}\subseteq X_{\sigma}.$ 

D'autre part, on a remarqu\'e que les rectangles $NE_{II}$ et $SO_{II}$ de la
configuration (3412) qui d\'efinit $\sigma$ ne rencontrent pas $\Gamma_w,$
on a donc, d'apr\`es le corollaire \ref{CompII},
$X_{\sigma}\subseteq \Sigma_w.$ Cela prouve le lemme
\ref{bonne-inter}. 
\end{proof}    

Pour terminer la preuve de la proposition \ref{inter}, nous allons
maintenant d\'emontrer le 

\begin{lem}\label{toutbon}
Toute famille $(X_{v_i})_{i\in [1,h]}$ de composantes
  irr\'eductibles des $\EE_{\pi_i},$ chacune de type Nord-Ouest ou Sud-Est,
  et associ\'ee \`a une configuration I d\'eg\'en\'er\'ee, est bonne.
\end{lem}

\renewcommand\proofname{Preuve du lemme \rm{\ref{toutbon}}}

\begin{proof}
On se donne une famille de composantes comme dans l'\'enonc\'e du
lemme, not\'ee $\mathcal F.$ Rappelons que, d'apr\`es le lemme
\ref{1-r}, la composante $v_1$ est de type Sud-Est, disons associ\'ee au
point $(c_1,\g_1),$ alors que la composante $v_h$ est de type Nord-Est,
disons associ\'ee au point $(b_h,\b_h).$ 
Si l'on a $b_h<c_1,$ alors $\mathcal F$ est bonne : en effet, d'apr\`es le
lemme \ref{1-r}, on a $b\le b_h$ et $c_1\le c,$ et il vient
$\winv(\a)=a<b_h<c_1<d=\winv(\d).$ 

On peut donc supposer $c_1<b_h.$
Posant $j_0=1$ et $i_0=h,$ supposons avoir construit des entiers
$$j_0<j_1<\cdots <j_{k-1}<i_{k-1}<\cdots<i_1<i_0$$ et des points
$(c_{j_l},\g_{j_l})\in \partial SE(j_l)$ et 
$(b_{i_l},\b_{i_l})\in \partial NO(i_l),$ pour 
$l=0,\ldots, k-1,$ tels que  
$$(\dagger)\hspace{2em}c_{j_0}<c_{j_1}<\cdots<c_{j_{k-1}}<b_{i_{k-1}}<\cdots<b_{i_0}$$ et 
que les entiers $j_{l+1}$ et $i_{l+1}$ soient d\'efinis par les
encadrements 
$$(\dagger\dagger)\hspace{2em}\left\{\begin{array}{l}\winv(\a-j_{l+1}+1)<c_{j_l}<\winv(\a-j_{l+1}),\\    
\winv(\a-i_{l+1}+1)<b_{i_l}<\winv(\a-i_{l+1}),\end{array}\right.$$ pour
$l=0,\ldots ,k-2.$   

On d\'efinit alors les 
entiers $j_k$ et $i_k$ par $$\winv(\a-j_k+1)<c_{j_{k-1}}<\winv(\a-j_k)$$ et
$$\winv(\a-i_k+1)<b_{i_{k-1}}<\winv(\a-i_k).$$ On a alors
$j_{k-1}<j_k\le i_k<i_{k-1}.$

Si $v_{j_k}$ est de type Nord-Ouest, d\'efinie par le point de
$\Gamma_w$ d'abscisse $b_{j_k},$ alors $b_{j_k}$ est dans
l'intervalle $]\,\winv(\a),\winv(\a-j_k+1)\,[,$ il existe donc $l\le
k-1,$ tel que 
$\winv(\a-{j_l}+1)<b_{j_k}<\winv(\a-{j_{l+1}}+1).$ On obtient alors,
en utilisant $(\dagger\dagger),$ que
$$\winv(\a-{j_l}+1)<b_{j_k}<c_{j_l}<\winv(\a-j_k),$$ et $\mathcal F$
est une bonne famille. 

Si $v_{j_k}$ est de type Sud-Est, d\'efinie par le point de
$\Gamma_w$ d'abscisse $c_{j_k},$ alors, d'apr\`es $(\dagger\dagger)$ et la
d\'efinition de $SE(j_k),$ on a
$c_{j_k}>\winv(\a-j_k)>c_{j_{k-1}}.$ Si de plus $c_{j_k}>b_{i_{k-1}},$ en
consid\'erant le plus petit entier $l$ tel que $b_{i_l}<c_{j_k},$ on a
$c_{j_k}<b_{i_{l-1}},$ et $b_{i_{l-1}}<\winv(\a-i_l)$ par
$(\dagger\dagger).$ On a aussi $b_{i_l}>\winv(\a-i_{l+1}+1)$ d'apr\`es
$(\dagger\dagger),$ mais $i_{l+1}\ge i_k\ge j_k,$ et comme $\winv$ est
d\'ecroissante sur $[\,\d',\a'\,],$ il vient $\winv(\a-i_{l+1}+1)\ge
\winv(\a-{j_k}+1).$ On a donc obtenu 
$$\winv(\a-{j_k}+1)<b_{j_l}<c_{j_k}<\winv(\a-i_l),$$ et
$\mathcal F$ est bonne.

On proc\`ede de m\^eme avec $v_{i_k}$ : on montre
d'abord que si $v_{i_k}$ est de type Sud-Est, $\mathcal F$ est une
bonne famille. Ensuite, si 
$v_{i_k}=t^{i_k}(b_{i_k})$ est de type Nord-Ouest, on remarque que
$b_{i_k}<b_{i_{k-1}}$ et l'on montre que si de plus 
$b_{i_k}<c_{j_{k-1}}$ alors la famille $\mathcal F$ est bonne. 

En regroupant ces r\'esultats, on constate que l'on a montr\'e que
$\mathcal F$ est une bonne famille, sauf, \'eventuellement,
si $j_k<i_k,\ v_{j_k}=t_{j_k}(c_{j_k})$ de type Sud-Est,
$v_{i_k}=t^{i_k}(b_{i_k})$ de type Nord-Ouest, et $b_{i_k},c_{j_k}\in\, 
]\,c_{j_{k-1}},b_{i_{k-1}}\,[.$ Si l'on a de plus
$b_{i_k}<c_{j_k},$ il vient d'une part $\winv(\a-j_k+1)<b_{i_k}$ car
$b_{i_k}>c_{j_{k-1}},$ par hypoth\`ese et $c_{j_{k-1}}>\winv(\a-j_k+1)$
par $(\dagger\dagger).$ On obtient de m\^eme $c_{j_k}<\winv(\a-i_k),$
d'o\`u finalement 
$$\winv(\a-j_k+1)<b_{i_k}<c_{j_k}<\winv(\a-i_k),$$ et la famille 
$\mathcal F$ est bonne. 

On peut donc supposer que
$c_{j_k}<b_{i_k}$ et l'on est \`a nouveau dans les hypoth\`eses de la
r\'ecurrence. Ce processus s'arr\^ete, puisque les suites d'entiers $b_{i_l}$
et $c_{j_l}$ sont strictement monotones et born\'ees, on en d\'eduit donc que
$\mathcal F$ est une bonne famille. Cela termine
la preuve du lemme \ref{toutbon}, et par cons\'equent aussi celle de la
proposition \ref{inter}
\end{proof}

\begin{ex}
On consid\`ere la permutation $$w=(6,7,5,1,8,4,2,3)$$
de $\mathfrak S_8.$

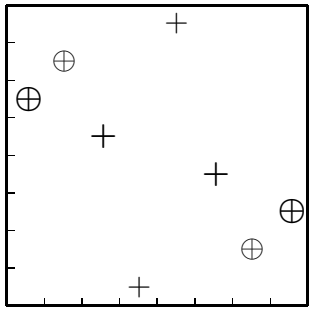
\begin{minipage}{4.5cm}

\begin{picture}(4.5,5)
\thicklines
\multiput(0,0.5)(0,4){2}{\line(1,0){4}}
\multiput(0,0.5)(4,0){2}{\line(0,1){4}}

\linethickness{0.01mm}
\multiput(0.5,0.5)(0.5,0){8}{\line(0,1){0.1}}
\multiput(0,1)(0,0.5){8}{\line(1,0){0.1}}
\thinlines
\put(0.6,3.65){$\oplus$}
\put(1.6,0.65){$+$}
\put(2.1,4.15){$+$}
\put(3.1,1.15){$\oplus$}

\boldmath
\put(0.1,3.15){$\oplus$}
\put(1.1,2.65){$+$}
\put(2.6,2.15){$+$}
\put(3.6,1.65){$\oplus$}
\unboldmath

\end{picture}
\end{minipage}
\begin{minipage}{8.5cm}
La configuration (3412) form\'ee par les points d'abscisses 1, 2, 7, 8 est
bien remplie et d'amplitude minimale. On
a $h=3,\ \a'=\a=6$ et $\d'=\d=3.$ Les permutations associ\'ees sont
$$\begin{array}{c}w_1=(3,7,6,1,8,5,2,4),\\
w_2=(4,7,3,1,8,6,2,5),\\
w_3=(5,7,4,1,8,3,2,6),\end{array}$$\end{minipage}
dont les graphes sont repr\'esent\'es ci-dessous.
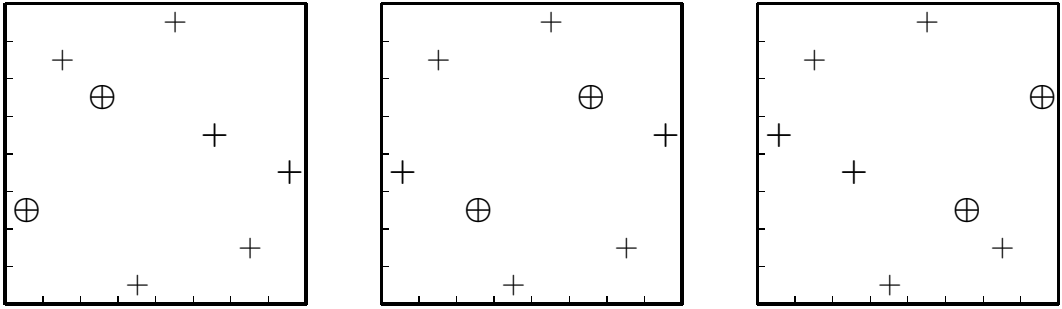
\begin{center}  
\begin{picture}(14,5)
\thicklines
\multiput(0,0)(5,0){3}{\line(1,0){4}}
\multiput(0,4)(5,0){3}{\line(1,0){4}}
\multiput(0,0)(5,0){3}{\line(0,1){4}}
\multiput(4,0)(5,0){3}{\line(0,1){4}}

\linethickness{0.01mm}
\multiput(0.5,0)(0.5,0){8}{\line(0,1){0.1}}
\multiput(0,0.5)(0,0.5){8}{\line(1,0){0.1}}

\multiput(5.5,0)(0.5,0){8}{\line(0,1){0.1}}
\multiput(5,0.5)(0,0.5){8}{\line(1,0){0.1}}

\multiput(10.5,0)(0.5,0){8}{\line(0,1){0.1}}
\multiput(10,0.5)(0,0.5){8}{\line(1,0){0.1}}

\thinlines

\multiput(0.6,3.15)(5,0){3}{$+$}
\multiput(1.6,0.15)(5,0){3}{$+$}
\multiput(2.1,3.65)(5,0){3}{$+$}
\multiput(3.1,0.65)(5,0){3}{$+$}
\boldmath
\put(0.1,1.15){$\oplus$}
\put(1.1,2.65){$\oplus$}
\put(2.6,2.15){$+$}
\put(3.6,1.65){$+$}

\put(5.1,1.65){$+$}
\put(6.1,1.15){$\oplus$}
\put(7.6,2.65){$\oplus$}
\put(8.6,2.15){$+$}

\put(10.1,2.15){$+$}
\put(11.1,1.65){$+$}
\put(12.6,1.15){$\oplus$}
\put(13.6,2.65){$\oplus$}

\end{picture}

\end{center}

Sur les trois diagrammes suivants, on a repr\'esent\'e les quadrants
associ\'es respectivement \`a chacune des trois quasi-r\'esolutions, ainsi
que les bandes verticales d'abscisses $\in\,]\,\winv_i(\a'),\winv_i(\d')\,[$ :

\begin{center}
\begin{picture}(14,5)

\thicklines
\multiput(0,0.5)(0,4){2}{\line(1,0){4}}
\multiput(0,0.5)(4,0){2}{\line(0,1){4}}

\multiput(5,0.5)(0,4){2}{\line(1,0){4}}
\multiput(5,0.5)(4,0){2}{\line(0,1){4}}

\multiput(10,0.5)(0,4){2}{\line(1,0){4}}
\multiput(10,0.5)(4,0){2}{\line(0,1){4}}

\linethickness{0.01mm}
\multiput(0.5,0.5)(0.5,0){8}{\line(0,1){0.1}}
\multiput(0,1)(0,0.5){8}{\line(1,0){0.1}}

\multiput(5.5,0.5)(0.5,0){8}{\line(0,1){0.1}}
\multiput(5,1)(0,0.5){8}{\line(1,0){0.1}}

\multiput(10.5,0.5)(0.5,0){8}{\line(0,1){0.1}}
\multiput(10,1)(0,0.5){8}{\line(1,0){0.1}}

\thinlines
\multiput(0.6,3.65)(5,0){3}{$\oplus$}
\multiput(1.6,0.65)(5,0){3}{$+$}
\multiput(2.1,4.15)(5,0){3}{$+$}
\multiput(3.1,1.15)(5,0){3}{$\oplus$}

\boldmath
\multiput(0.1,3.15)(5,0){3}{$\oplus$}
\multiput(1.1,2.65)(5,0){3}{$+$}
\multiput(2.6,2.15)(5,0){3}{$+$}
\multiput(3.6,1.65)(5,0){3}{$\oplus$}
\unboldmath

\linethickness{0.05mm}
\multiput(1.5,1.5)(0.2,0){13}{\line(1,0){0.1}}
\multiput(1.5,0.5)(0,0.2){5}{\line(0,1){0.1}}
\multiput(0.5,0.5)(0,0.2){20}{\line(0,1){0.1}}
\multiput(1,0.5)(0,0.2){20}{\line(0,1){0.1}}

\multiput(5,3.5)(0.2,0){5}{\line(1,0){0.1}}
\multiput(6,3.5)(0,0.2){5}{\line(0,1){0.1}}
\multiput(6.5,0.5)(0,0.2){20}{\line(0,1){0.1}}
\multiput(7.5,0.5)(0,0.2){20}{\line(0,1){0.1}}
\multiput(8,1.5)(0.2,0){5}{\line(1,0){0.1}}
\multiput(8,0.5)(0,0.2){5}{\line(0,1){0.1}}

\multiput(10,3.5)(0.2,0){13}{\line(1,0){0.1}}
\multiput(12.5,3.5)(0,0.2){5}{\line(0,1){0.1}}
\multiput(13,0.5)(0,0.2){20}{\line(0,1){0.1}}
\multiput(13.5,0.5)(0,0.2){20}{\line(0,1){0.1}}
\end{picture}
\end{center}

Ainsi, d'apr\`es la proposition \ref{compex}, chacun des $\EE_{\pi_i}$ a
deux composantes irr\'eductibles. Les deux composantes de $\EE_{\pi_1}$
sont de type Sud-Est, $t_1(4)=(1,7,6,5,8,4,2,3)$ et
$t_1(7)=(2,7,6,1,8,5,4,3).$ Le lieu $\EE_{\pi_2}$ a une composante de
type Nord-Ouest et une de type Sud-Est, $t^2(2)=(6,5,4,1,8,7,2,3)$ et
$t_2(7)=(6,7,2,1,8,5,4,3).$ Enfin, $\EE_{\pi_3}$ a deux
composantes de type Nord-Ouest, $t^3(2)=(6,5,4,1,8,3,2,7)$ et
$t^3(5)=(6,7,5,1,4,3,2,8).$ 

Aucune de ces composantes n'est contenue dans $\Sigma_w$ : elles
correspondent toutes \`a des configurations I d\'eg\'en\'er\'ees, qui donnent
des points lisses d'apr\`es le th\'eor\`eme \ref{gencomp}. En revanche,
d'apr\`es le lemme \ref{toutbon}, toute intersection $X_{v_1}\cap
X_{v_2}\cap X_{v_3},$ o\`u chaque $X_{v_i}$ est une composante de
$\EE_{\pi_i},$ est contenue dans $\Sigma_w.$ 

Voyons par exemple
l'intersection $X_{t_1(4)} \cap X_{t_2(7)}\cap X_{t^3(5)}.$ Les
entiers $i=2$ et $j=3$ remplissent les conditions de la d\'efinition
\ref{good} : ($a$) $t_2(7)$ est de type Sud-Est, associ\'ee au point
$(7,2)$ ; ($b$) $t^3(5)$ est de type Nord-Ouest, associ\'ee au point
$(5,9)$ ; ($c$) on a $\a=6$ et $\winv(5)=3<b_3=5<c_2=7<\winv(3)=8.$ 

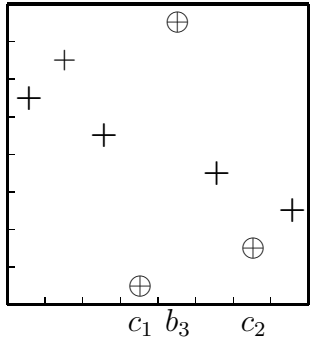
\begin{center}
\begin{picture}(4.5,5)
\thicklines
\multiput(0,0.5)(0,4){2}{\line(1,0){4}}
\multiput(0,0.5)(4,0){2}{\line(0,1){4}}

\linethickness{0.01mm}
\multiput(0.5,0.5)(0.5,0){8}{\line(0,1){0.1}}
\multiput(0,1)(0,0.5){8}{\line(1,0){0.1}}
\thinlines
\put(0.6,3.65){$+$}
\put(1.6,0.65){$\oplus$}
\put(2.1,4.15){$\oplus$}
\put(3.1,1.15){$\oplus$}

\put(1.6,0.15){$c_1$}
\put(2.1,0.15){$b_3$}
\put(3.1,0.15){$c_2$}
\boldmath
\put(0.1,3.15){$+$}
\put(1.1,2.65){$+$}
\put(2.6,2.15){$+$}
\put(3.6,1.65){$+$}
\unboldmath

\end{picture}
\end{center}

D'apr\`es le lemme \ref{bonne-inter},
l'intersection $X_{t_2(7)}\cap X_{t^3(5)}$ est contenue dans
$X_{\sigma},$ o\`u $\sigma$ est la permutation associ\'ee \`a la
configuration II pure form\'ee par les points d'abscisses 3, 5, 6, 7, 8,
c.-\`a-d. $\sigma=(6,7,2,1,5,4,3,8).$
\end{ex}

\renewcommand\proofname{Preuve}
%
\subsubsection{Correspondance entre configurations}

Nous terminons par la preuve de la proposition \ref{corr-config}.

Pour toute configuration $\mathcal K$ de $w_i,$ de type I ou II, on
note $w_Iw_{J_i}(\mathcal K)$ l'ensemble des points $(x,w(x))$ tels que
$(x,w_i(x))\in\mathcal K.$  

\begin{propB}
Pour toute configuration $\mathcal K$ de $w_i,$ de type I ou II,
param\'etrant une composante irr\'eductible $X_v$ du lieu singulier de
$X_{w_i},$ on a :\\ 
\hspace*{3ex}$\bullet$ ou bien $P_I\times^{P_{J_i}}X_v\subseteq
\Ex(\pi_i),$\\
\hspace*{3ex}$\bullet$ ou bien $w_Iw_{J_i}(\mathcal K)$ est une
configuration du m\^eme type de $w,$ et 
$\pi_i(P_I\times^{P_{J_i}}X_v) = X_{w_Iw_{J_i}v}$ est la composante
irr\'eductible du lieu singulier de $X_w$ associ\'ee.  
\end{propB}

Au cours de la d\'emonstration, nous utiliserons la notion
suivante 

\begin{defn}
Soit $z$ une permutation et $\mathcal F$ une famille de points du
graphe de $z,\ \mathcal F =\{(x_1,z(x_1)),\ldots, (x_m,z(x_m))\},$
avec $x_1<\cdots< x_m.$ 
Pour $j=1,\ldots, m-1,$ on appelle {\it successeur de $(x_j,z(x_j))$ dans
$\mathcal F$} le point $(x_{j+1},z(x_{j+1})).$ 
\end{defn}

\begin{proof}
On consid\`ere une configuration $\mathcal K$ de $w_i,$ param\'etrant une
composante irr\'eductible $X_v$ de Sing $X_{w_i}.$ On suppose de plus que  
$P_I\times^{P_{J_i}}X_v\not\subseteq\Ex(\pi_i).$ 

D'apr\`es la discussion au d\'ebut de la preuve de la proposition \ref{compex}, il
vient que $w_{J_i}v\in\, ^I\mathfrak S_{min},$ de sorte que
$w_I*v=w_Iw_{J_i}v.$  
Il suffit alors de montrer que $w_Iw_{J_i}(\mathcal K)$ est une
configuration du m\^eme type de $w.$ En effet, la derni\`ere assertion en d\'ecoule :
on a $v=\gamma ({\scriptstyle\mathcal K}) w_i,$ o\`u $\gamma
({\scriptstyle\mathcal K})$ est un certain cycle 
sur les ordonn\'ees de la configuration $\mathcal K,$ et de m\^eme, la
composante irr\'eductible de Sing $X_w$ associ\'ee \`a $w_I w_{J_i}(\mathcal K)$
de $w$ est donn\'ee par $v' =\gamma ({\scriptstyle w_I w_{J_i}(\mathcal
  K)}) w,$ et comme $\gamma ({\scriptstyle z(\mathcal
  K)})=z\,\gamma ({\scriptstyle\mathcal K})\, \zinv,$ pour tout $z,$
il vient $v'=w_I w_{J_i}\gamma ({\scriptstyle \mathcal K}) w_i=w_I w_{J_i}v.$ 

On remarque de plus que l'on a $(\d',\a')\,v\not\le w_i.$ En effet, notant
$v'=(\d',\a')\,v,$ on a $v'=w_{J_i}s_{k_i}\t,$ et comme $\t\in\,
^I\mathfrak S_{min},$ il vient $v'\in \,^{J_i}\mathfrak S_{max}.$ Si
l'on avait $v'\le w_i,$ alors $P_I\times^{P_{J_i}}X_{v'}\subseteq
Z_i$ aurait m\^eme image que $P_I\times^{P_{J_i}}X_v,$ d'o\`u
$P_I\times^{P_{J_i}}X_v\subseteq \Ex(\pi_i),$ une 
contradiction. 

On note $\mathcal S_i^S$ (resp. $\mathcal S_i^N$) l'ensemble des
points du graphe de $w_i$ dont l'ordonn\'ee est dans $[\d',k_i]$
(resp. $[k_i+1,\a']$), et $\mathcal S_i= \mathcal S_i^S\cup \mathcal
S_i^N.$ On note aussi $\mathcal S=w_Iw_{J_i}(\mathcal S_i).$

\vspace*{1em}
{\bf 1.} Supposons pour commencer que $\mathcal K$ est une configuration
I, n\'ecessairement non d\'eg\'en\'er\'ee puisqu'elle param\`etre un point singulier de
$X_{w_i}.$ On conserve les notations des d\'efinitions \ref{confI} et
\ref{PtCoBB}. 
On note $\mathcal K^N$ la r\'eunion de $\{P_+\}$ et de
la suite NE, et $\mathcal 
K^S$ la r\'eunion de $\{P_-\}$ et de la suite SO. Si la configuration
$\mathcal K$ ne rencontre pas $\mathcal S_i,$ le r\'esultat est clair. 
Supposons au contraire que $\mathcal K\cap \mathcal
S_i\not=\emptyset.$ Soit $p_2$ la deuxi\`eme projection de $[1,n]^2$ sur
$[1,n].$ Comme $\winv_i$ est d\'ecroissante sur chacun des
$p_2(E),$ pour $E=\mathcal S_i^S,\mathcal S_i^N, \mathcal K^S,
\mathcal K^N,$ et que $\mathcal S_i^S$ est au Sud-Ouest de $\mathcal
S_i^N,$ chacune des deux parties $\mathcal K^S$ et $\mathcal 
K^N$ ne peut rencontrer qu'au plus l'une de $\mathcal S_i^S$ et $\mathcal
S_i^N.$

Montrons de plus que $\mathcal S_i$ ne peut rencontrer simultan\'ement $\mathcal
K^N$ et $\mathcal K^S.$ En effet, supposons que ce soit le cas. Alors,
n\'ecessairement, $\mathcal K^N$ rencontre $\mathcal S_i^N$ et $\mathcal
K^S$ rencontre $\mathcal S_i^S,$ et cela entra\^\i ne que $\mathcal S_i$
ne contient ni $P_+$ ni $P_-.$ Ainsi $\mathcal K^N\cap \mathcal S_i^N$
est contenu dans la suite NE, et $\mathcal K^S\cap \mathcal S_i^S$ dans
la suite SO. Comme $P_-$ est au Sud-Est de la suite NE, on en d\'eduit
que $\xinf>\winv_i(\a'),$ d'o\`u $\xinf>\winv_i(\d').$ 

Soit $y_m=\inf\ \{y\ |\ (\winv_i(y),y)\in\mathcal K^S\cap\mathcal
S_i^S\}$ ; on a $y_m=\d',$ car sinon le point $(\winv_i(\d'),\d')$
serait contenu dans $\Gamma_{w_i}\cap \mathcal R_{(\Xinf,\xinf)}$ mais
pas dans $\bigcup_{j=1}^t SO(x_{-j},y_{-j}) \cup \bigcup_{j=1}^s
NE(x_j-1,y_j-1),$ ce qui contredirait $(\bigtriangleup).$ On montre de
m\^eme que $\a'$ est l'ordonn\'ee d'un point de la suite NE. Mais on voit alors que
$v':=(\d',\a')\,v<w_i.$ En effet, les graphes $\Gamma_{v'}$ et
$\Gamma_{w_i}$ co\"\i ncident en dehors des ordonn\'ees de $\mathcal Y$
({\it cf.} 3.1), il suffit donc de comparer les permutations obtenues
en focalisant sur ces ordonn\'ees. Posant 
$n'=s+t+2,$ on a alors 
$$\begin{array}{l}\overline{w_i}=(n',\  t+1, \ldots , 2,\ s+t+1,
  \ldots,t+2,\ 1)\\
\overline{v}=(t+1, \ldots ,1\ n', \ldots,t+2)\end{array}$$ 
et $\overline{v'}=(j_1,j_2)\overline{v}$ avec $j_1\in
[2,t+1]$ et $j_2\in [t+2,s+t+1].$ On a donc
$\overline{v'}<\overline{w_i},$ d'o\`u $v'<w_i,$ exclu. On a ainsi
d\'emontr\'e que $\mathcal K\cap \mathcal S_i$ est contenu dans $\mathcal
K^N$ ou bien dans $\mathcal K^S,$ donc \'egal \`a l'un des quatre
$\mathcal K^{\ast}\cap\mathcal S_i^{\ast'}$ o\`u $\ast,\ast'\in \{N,S\}.$

Supposons par exemple que $\mathcal K \cap\mathcal S_i =
\mathcal K^N\cap\mathcal S_i^S.$ Alors il
est clair que les coordonn\'ees des points de $w_Iw_{J_i}(\mathcal K)$
v\'erifient les in\'egalit\'es requises ({\it cf.} 3.1), et il suffit donc
de v\'erifier l'inclusion $(\bigtriangleup).$ Notons 
$$\begin{array}{rcl}\mathcal K &=& \left\{(\Xinf,\Yinf),\
    (\xinf,\yinf)\}\cup\{(x_j,y_j),\ j\in 
  [-t,-1]\cup [1,s]\right\},\\ 
w_Iw_{J_i}( \mathcal K) & =& \left\{(\Xinf,\Yinf'),\
  (\xinf,\yinf')\}\cup\{(x_j,y'_j),\ j\in [-t,-1]\cup
  [1,s]\right\}.\end{array}$$ 
Comme les graphes $\Gamma_{w_i}$ et $\Gamma_w$ co\"\i ncident sur
les ordonn\'ees hors de $[\,\d',\a'\,],$ il suffit  
de montrer que  
$$\mathcal S\cap \mathcal R_{(\Xinf,\xinf)}(w)\subseteq
  \bigcup_{j=1}^t SO(x_{-j},y'_{-j}) \cup \bigcup_{j=1}^s NE(x_j-1,y'_j-1)\hspace*{3em}(\bigtriangleup').$$   

Soit $y_m=\inf\ \{y \ |\ (\winv_i(y),y)\in\mathcal K^N\cap\mathcal
S_i^S\}$ ; on a $m\in [1,s]\cup \{\infty\}.$ Soit $m'$ l'entier tel
que $x_{m'}=v^{-1}(y_m)$ ; on a $m'\in [1,s]\cup
\{-\infty\}.$ Le point $(x_{m'},y_{m'})$ est \`a l'Ouest, au sens
large, du successeur de $(\winv_i(y_m),y_m)$ dans $\mathcal S_i.$ En effet, si
$y_m\not=\d',$ alors le successeur de $(\winv_i(y_m),y_m)$ dans
$\mathcal S_i$ est $(\winv_i(y_m-1),y_m-1).$ Si l'on avait
$\winv_i(y_m-1)<x_{m'},$ alors le point $(\winv_i(y_m-1),y_m-1)$
serait contenu dans $\Gamma_{w_i}\cap \mathcal
R_{(\Xinf,\xinf)}(w_i),$ et pas dans $\bigcup_{j=1}^t
SO(x_{-j},y_{-j}) \cup \bigcup_{j=1}^s NE(x_j-1,y_j-1),$ ce qui
contredit $(\bigtriangleup).$ D'autre part, si $y_m=\d',$ alors le
successeur de $(\winv_i(y_m),y_m)$ dans $\mathcal S_i$ est
$(\winv_i(\a'),\a').$ Si l'on avait $\winv_i(\a')<x_{m'},$ alors on
aurait $(\d',\a')\,v<v,$ d'o\`u $(\d',\a')\,v<w_i,$ exclu.

Il vient alors : si $m'=-\infty,\ \mathcal S\cap \mathcal
R_{(\Xinf,\xinf)}(w)\subseteq w_Iw_{J_i}(\mathcal K),$
d'o\`u $(\bigtriangleup'),$ et si $m'\not=-\infty,$ comme on a 
$w_i(x_{m'})=y_{m'}<\d',$ alors 
$w(x_{m'})=w_i(x_{m'}),$ ou encore $y'_{m'}=y_{m'}$ et alors  
$$\mathcal S \cap \mathcal R_{(\Xinf,\xinf)}(w) \subseteq
w_Iw_{J_i}(\mathcal K) \cup NE(x_{m'}-1,y'_{m'}-1),$$ d'o\`u l'on d\'eduit 
$(\bigtriangleup').$

Le cas $\mathcal K \cap\mathcal S_i = \mathcal K^S\cap\mathcal
S_i^N$ est semblable, et les cas o\`u $\mathcal K \cap\mathcal S_i$ est
\'egal \`a $\mathcal K^S\cap\mathcal S_i^S$ ou $\mathcal K^N\cap\mathcal
S_i^N$ sont similaires, mais plus simples.

On a donc d\'emontr\'e le r\'esultat lorsque $\mathcal K$ est une configuration
I.  
\vspace*{1em}

{\bf 2.} Supposons maintenant que $\mathcal K$ est une configuration II : comme
$X_v$ est une composante irr\'eductible du lieu singulier, on a, d'apr\`es
le corollaire \ref{CompII}, $s=t=0$ ou $r=0.$ 

{\bf 2.1.} Traitons pour commencer le cas d'une configuration
pure. Notons $A,B,C$ et $D$ les 
quatre points de la configuration (3412) qui d\'etermine $\mathcal K,$ et
$(x_A,y_A),$ {\it etc.} leurs coordonn\'ees. On va montrer que 
$w_Iw_{J_i}(\{A,B,C,D\})$ est une configuration (3412) de $w,$
incompressible, et telle que les zones $NE_{II}(w)$ et $SO_{II}(w)$
associ\'ees ne contiennent pas de point de $\Gamma_w.$ Comme
pr\'ec\'edemment, il suffit de montrer que   
$$\big(MN(w)\cup NE_{II}(w)\cup ME(w) \cup MS(w)\cup SO_{II}(w)\cup
MO(w) \big)\cap 
 \mathcal S =\emptyset\ \ \ (\bigtriangledown).$$

Comme l'ensemble des ordonn\'ees des points de $\mathcal S_i$ est un
intervalle, le r\'esultat est clair si $\mathcal K\cap \mathcal
S_i=\emptyset.$ On suppose donc $\mathcal K\cap \mathcal
S_i\not=\emptyset.$ 

{\bf 2.1.1.} Supposons pour commencer $B\in \mathcal S_i.$ Alors le
seul autre point de $\mathcal K$ qui pourrait appartenir \`a $\mathcal
S$ est $A,$ mais on aurait alors $(\d',\a')\,v<v,$ exclu. On a donc
$\mathcal K\cap \mathcal S_i=\{B\}.$ Alors il est clair
que $w_Iw_{J_i}(\{A,B,C,D\})$ est une configuration (3412) de $w,$ not\'ee
$\II.$ Si $B\in \mathcal S_i^N,$ alors comme on a vu que $A\not\in
\mathcal S_i,$ on a $p_2(\mathcal S_i^S)\subseteq\, ]\,y_A,y_B\,[,$ et
l'on obtient que $\II$ est de la forme voulue. D'autre part, si $B\in
\mathcal S_i^S,$ on a d'abord $y_A<\d',$ et donc $\big(ME(w) \cup
MS(w)\cup SO_{II}(w)\cup MO(w) \big)\cap \mathcal S =\emptyset.$ 
Ensuite, si $y_B\not=\d',$ le successeur de $B$
dans $\mathcal S_i$ a pour ordonn\'ee $y_B-1\in \ ]\,y_A, y_B\,[,$ et comme
$\mathcal K$ est incompressible et que le 
rectangle $NE_{II}(\mathcal K)$ ne rencontre pas $\Gamma_{w_i},$ on a
n\'ecessairement $x_D<\winv_i(y_B-1).$ D'autre part, si l'on a
$y_B=\d',$ alors le successeur de $B$ dans $\mathcal S_i$ est le point
d'ordonn\'ee $\a',$ et l'on a $x_D<\winv(\a'),$ car sinon on aurait
$(\d',\a')\,v<v.$ On en d\'eduit, dans les deux cas,
$\big(MN(w)\cup NE_{II}(w)\big) \cap 
 \mathcal S =\emptyset.$ On a ainsi obtenu $(\bigtriangledown)$
 lorsque $B\in \mathcal S_i.$ Le cas o\`u $C\in\mathcal S_i$ se traite
 de fa\c con semblable.  

{\bf 2.1.2.} Il reste donc \`a traiter le
cas $\mathcal K\cap \mathcal S_i\subseteq \mathcal K^M,$ o\`u $\mathcal
K^M$ d\'esigne la r\'eunion de $A,\ D$ et de la suite
centrale. Dans ce cas, $\mathcal K\cap \mathcal S_i$ est  
contenue dans l'une des deux zones $\mathcal S_i^N$ ou $\mathcal
S_i^S,$ disons $\mathcal S_i^S.$ 

Si l'on a $A\in \mathcal S_i^S$ et $D\not\in\mathcal S_i^S$ :
alors on a d'une part $y_D<\d',$ donc $\big(MS(w)\cup SO_{II}(w)\big)\cap
\mathcal S =\emptyset.$ D'autre part, on a $y_B>\a',$ donc, comme la
configuration $\mathcal K$ est incompressible et que $NE_{II}(\mathcal
K)$ ne rencontre pas $\Gamma_{w_i},$ les abscisses de $B$ et $D$ sont
contenues dans un m\^eme intervalle de la subdivision donn\'ee par les
abscisses des points de $\mathcal S_i.$ Alors les prori\'et\'es voulues
concernant les rectangles $MN(w),\ NE_{II}(w),\ MO (w)$ et $ME(w)$
r\'esultent de leur analogue dans $\Gamma_{w_i}.$ On a donc \'etabli
$(\bigtriangledown)$ dans ce cas. Les
autres cas, plus simples, sont laiss\'es au lecteur. On a donc d\'emontr\'e
le r\'esultat pour $\mathcal K$ de type II telle que $s=t=0.$

\vspace*{1em}

{\bf 2.2.} Soit maintenant $\mathcal K$ de type II mixte. On
d\'ecompose alors $\mathcal K$ en trois parties : on note $\mathcal K^N$ la
r\'eunion de $\{B\}$ et de la suite NE, $\mathcal K^M=\{A,D\}$ et
$\mathcal K^S$ est la r\'eunion de $\{C\}$ et de la suite SO. On va
montrer que $w_Iw_{J_i}(\{A,B,C,D\})$ est une configuration (3412) de 
$w,$ incompressible et de zone centrale ne contenant pas de point de
$\Gamma_w,$ et que 
$$\displaylines{\hspace*{3em}\Gamma_w\cap\big(SO_{II}(w)\cup NE_{II}(w)\big)
  \subseteq \hfill\cr 
\hfill\bigcup_{j=1}^t SO(x_{-j},y'_{-j}) \cup \bigcup_{j=1}^s
NE(x_j-1,y'_j-1)\ \ \ (\Diamond)\hspace*{3em}}$$ 
avec les m\^emes conventions d'\'ecriture que pr\'ec\'edemment. A nouveau,
on peut remplacer $(\Diamond)$ par  
$$\displaylines{\hspace*{3em}\mathcal S\cap \big(SO_{II}(w)\cup
  NE_{II}(w)\big)
  \subseteq\hfill\cr 
\hfill\bigcup_{j=1}^t SO(x_{-j},y'_{-j}) \cup \bigcup_{j=1}^s
NE(x_j-1,y'_j-1)\ \ \ (\Diamond').\hspace*{3em}}$$ 
Le r\'esultat est clair si $\mathcal K\cap \mathcal
S_i=\emptyset,$ on suppose donc $\mathcal K\cap \mathcal
S_i\not=\emptyset.$ On remarque que l'intersection de chacune des
parties $\mathcal S_i^S$ et $\mathcal S_i^N$ avec $\mathcal K$ est
contenue dans l'une de $\mathcal K^N,\ \mathcal K^M$ et $\mathcal K^S.$

{\bf 2.2.1.} Supposons pour commencer que $\mathcal K$ ne rencontre
qu'une seule de $\mathcal S_i^S$ et $\mathcal S_i^N.$ Les deux cas
sont sym\'etriques, il suffit donc de traiter par exemple $\mathcal
K\cap \mathcal S_i\subseteq \mathcal K\cap \mathcal S_i^S.$ On suppose
d'abord $\mathcal K\cap \mathcal S_i^S =\mathcal K^N \cap \mathcal S_i^S.$ 
Alors il est clair que $w_Iw_{J_i}(\{A,B,C,D\})$
est une configuration (3412) de $w.$ Par ailleurs, on a $y_A<\d',$ donc
les propri\'et\'es voulues concernant les rectangles $MO(w),\ C(w),$
$ME(w),\ MS(w),$ et $SO_{II}(w)$ r\'esultent de leur analogue dans
$\Gamma_{w_i}.$ Les abscisses de $B$ et $C$ sont n\'ecessairement
contenues dans un m\^eme intervalle de la subdivision donn\'ee par les
abscisses des points de $\mathcal S_i$ (sinon $MN(w_i)$ contiendrait
des points de $\mathcal S_i$). Il en r\'esulte que $MN(w)\cap
\Gamma_w=\emptyset.$ Concernant $NE_{II}(w),$ il suffit de montrer que 
$$NE_{II}(w)\cap w_Iw_{J_i}(\mathcal S_i^N)\subseteq  \bigcup_{j=1}^s
NE(x_j-1,y'_j-1)\ \ \ (\Diamond'').$$ Or on montre que l'on a
$x_D<\winv_i(\a')$ ou bien 
qu'il existe $j\in [1,s]$ tel que $x_j<\winv_i(\a')$ et $y_j<\d'.$
Dans le premier cas, il vient $$NE_{II}(w)\cap w_Iw_{J_i}(\mathcal
S_i^N)=\emptyset,$$ et dans le second, on a $y'_j=y_j$ et on obtient
$$NE_{II}(w)\cap w_Iw_{J_i}(\mathcal S_i^N)\subseteq NE(x_j-1,y'_j-1).$$ 
On a donc \'etabli $(\Diamond'').$ Cela prouve le r\'esultat voulu dans le
cas $\mathcal K\cap \mathcal S_i=\mathcal K^N \cap \mathcal S_i^S.$ 
Le cas $\mathcal K\cap \mathcal S_i = \mathcal
K^S\cap\mathcal S_i^S,$ plus simple que celui que
nous venons de d\'etailler, est laiss\'e au lecteur.

Le cas $\mathcal K\cap \mathcal S_i \subseteq \mathcal K\cap\mathcal
S_i^S\subseteq \mathcal K^M$ se traite de la m\^eme mani\`ere que son
analogue trait\'e en 2.1.2.

{\bf 2.2.2.} Il reste pour terminer \`a montrer que $\mathcal K$ ne peut pas
rencontrer \`a la fois $\mathcal S_i^N$ et $\mathcal S_i^S.$ On remarque
d'abord que si $\mathcal S_i^S$ rencontre $\mathcal K^N$ ou contient
$D,$ alors $\mathcal S_i^N\cap \mathcal K =\emptyset.$ De m\^eme, si 
$\mathcal S_i^N$ rencontre $\mathcal K^S$ ou contient $A,$ alors 
$\mathcal S_i^S\cap \mathcal K =\emptyset.$ Ainsi, si $\mathcal K$
rencontrait \`a la fois $\mathcal S_i^N$ et $\mathcal S_i^S,$ on aurait 
$$\mathcal S_i^S\cap \mathcal K\subseteq \{A\}\cup \mathcal K^S,\mbox{ et }
\mathcal S_i^N\cap \mathcal K\subseteq \{D\}\cup \mathcal K^N.$$ De
plus, on n'a pas simultan\'ement $\mathcal S_i^S\cap \mathcal
K\subseteq\mathcal K^S$ et $\mathcal S_i^N\cap \mathcal
K\subseteq\mathcal K^N,$ car les ordonn\'ees de $A$ et $D$
sont entre celles des points de $\mathcal K^S$ et $\mathcal K^N.$ Les
possibilit\'es qui restent sont donc 
$$A\in \mathcal S_i^S  \mbox{ et } \mathcal S_i^N\cap \mathcal
K^N\not=\emptyset$$ ou, sym\'etriquement,  
$$D\in \mathcal S_i^N  \mbox{ et } \mathcal S_i^S\cap \mathcal
K^S\not=\emptyset.$$
Mais on montre alors que l'on aurait $(\d',\a')\,v<v,$ exclu.
Cela termine la preuve de la proposition \ref{corr-config}.
\end{proof}

%


\begin{thebibliography}{00}
%
\bibitem{Bia}A. Bialynicki-Birula, {\it On induced actions of
    algebraic groups}, Ann. Inst. Fourier
  (Grenoble) {\bf 43} (2) (1993), 365-368.
%
\bibitem{BW}S. Billey, G. Warrington, {\it Maximal singular loci
    of Schubert varieties in $SL(n)/B$}, math.AG/0102168.
%
\bibitem{BP}M. Brion, P. Polo, {\it Generic singularities of
    certain Schubert varieties}, Math. Z. {\bf 231} (1999), 301-324.  
%
\bibitem{BV}S. Buoncristiano, A. B. Veit, {\it The intrinsic
    stratification of a Schubert variety}, Adv. Math. {\bf 91} (1)
    (1992), 1-26. 
%
\bibitem{covex}A. Cortez, {\it Singularit\'es g\'en\'eriques des
  vari\'et\'es de Schubert covexillaires}, Ann. Inst. Fourier
  (Grenoble) {\bf 51} (2) (2001), 375-393.
%
\bibitem{Note}A. Cortez, {\it Singularit\'es g\'en\'eriques et
    quasi-r\'esolutions des vari\'et\'es de Schubert pour le groupe
    lin\'eaire}, C. R. Acad. Sci. Paris {\bf 333} (2001), 561-566.
%
\bibitem{propZ} V. Deodhar, {\it Some characterizations of Bruhat
    ordering on a Coxeter group and determination of the relative 
Möbius function}, Invent. Math. {\bf 39} (1977), no. 2, 187-198.
%
\bibitem{Dyer}M. J. Dyer, {\it Hecke algebras and shellings of Bruhat
    intervals}, Compositio Math. {\bf 89} (1993), no. 1, 91-115.
%
\bibitem{F2}W. Fulton, {\it Flags, Schubert polynomials , degeneracy
    loci, and determinantal formulas}, Duke Math. J. {\bf 65} (1992),
    381-420. 
%
\bibitem{Gash}V. Gasharov, {\it Sufficiency of Lakshmibai-Sandhya
singularity conditions for Schubert varieties}, Compositio Math. {\bf 126} (2001), 47-56. 
%
\bibitem{Gon}N. Gonciulea, {\it Singular loci of varieties of
    complexes, II}, J. Algebra {\bf 235} (2) (2001), 547-558. 
%
\bibitem{KLR}C. Kassel, A. Lascoux, C. Reutenauer, {\it The singular
    locus of a Schubert variety}, preprint, IRMA n$^{\circ}$2001-004.
%
\bibitem{K-L}D. Kazhdan, G. Lusztig, {\it Representations of Coxeter
  groups and Hecke algebras}, Inv. Math. {\bf 53} (1979), 165-184.
%
\bibitem{K}G. Kempf, {\it Algebraic varieties}, London Mathematical
  Society Lecture Note Series {\bf 172}, Cambridge University Press,
  Cambridge, 1993. 
%
\bibitem{lemme-d} V. Lakshmibai, C. Musili, C. S. Seshadri, {\it
    Geometry of $G/P$. IV. Standard monomial theory for classical
    types}, Proc. Indian Acad. Sci. Sect. A (Math. Sci.) {\bf 88}
    (1979), no. 4, 279-362. 
%
\bibitem{L-Sa}V. Lakshmibai, B. Sandhya, {\it Criterion for
  smoothness of Schubert varieties in $SL(n)/B$}, Proc. Indian
  Acad. Sci. (Math. Sci.) {\bf 100} (1990), 45-52.
%
\bibitem{L-Se}V. Lakshmibai, C. S. Seshadri, {\it Singular locus of a
    Schubert variety}, Bull. Amer. Math. Soc. {\bf 11} (1984),
    363-366.
%
\bibitem{LS-2}A. Lascoux, M.-P. Sch\"utzenberger, {\it
  Treillis et bases des groupes de Coxeter}, Electron. J. Combin.
  {\bf 3} (2) (1996) (The Foata Festschrift volume). 
%
\bibitem{Man}L. Manivel, {\it Fonctions sym\'etriques, polyn\^omes de
    Schubert et lieux de d\'eg\'en\'erescence}, Cours sp\'ecialis\'es
    {\bf 3}, Soci\'et\'e Math\'ematique de France, Paris, 1998. 
%
\bibitem{Man2}L. Manivel, {\it Le lieu singulier des vari\'et\'es de
    Schubert}, Int. Math. Res. Notices {\bf 16} (2001), 849-871.
%
\bibitem{Man3}L. Manivel, {\it Generic singularities of Schubert
    varieties}, math.AG/\linebreak 0105239.
%
\bibitem{Mat}H. Matsumura, {\it Commutative ring theory}, Cambridge
  studies in advanced mathematics {\bf 8}, Cambridge University Press,
  1986.
%
\bibitem{Mum}D. Mumford, {\it The red book of varieties and schemes},
  Lecture notes in mathematics {\bf 1358}, Springer, 1999.
%
\bibitem{polo}P. Polo, {\it Construction of arbitrary Kazhdan-Lusztig
    polynomials in symmetric groups}, Represent. Theory {\bf 3} (1999),
    90-104.  
%
%
\bibitem{Ram}S. Ramanan, A. Ramanathan, {\it Projective normality of flag 
varieties and Schubert varieties}, Invent. Math. {\bf 79} (1985),
217-224. 
%
\bibitem{Ry}K. Ryan, {\it On Schubert varieties in the flag manifold
    of SL$(n,\C)$}, Math. Ann. {\bf 276} (1987),205-224.
%
\bibitem{Ser}J.-P. Serre, {\it Espaces fibr\'es alg\'ebriques}, Anneaux
  de Chow et applications, S\'eminaire Chevalley, E.N.S. Paris, 1958.

\end{thebibliography}
\end{document}